\documentclass[10pt]{article}

\usepackage{a4wide}
\usepackage{amssymb}
\usepackage{amsfonts}
\usepackage{amsmath}
\input xy
\xyoption{arrow} \xyoption{matrix}

\date{}

\newtheorem{proposition}{Proposition}[section]
\newtheorem{theorem}[proposition]{Theorem}
\newtheorem{lemma}[proposition]{Lemma}

\newtheorem{definition}[proposition]{Definition}
\newtheorem{corollary}[proposition]{Corollary}

\def\der{\partial }

\def\nFM0{{\nu }_{F,M_0}}
\def\nFN0{{\nu }_{F,N_0}}
\def\nGN0{{\nu }_{G,N_0}}

\def\N0{ {\bf N}_0 }

\def\t{\otimes}

\def\v{\varphi}
\def\ra{\rightarrow}

\def\Xpm{X^{\pm }}

\def\s{\sigma}
\def\Z{{\bf Z }}

\def\l1{{\lambda}_1}

\def\a{\alpha}
\def\a0{ {\alpha }_0}
\def\a1{ {\alpha }_1}

\def\l{\lambda}

\def\nFGM0{{\nu }_{F,G,M_0}}

\def\nFN0{{\nu}_{F,N_0}}

\def\sm{{\sigma}^m}

\def\sm1{{\sigma}^{-1}}

\def\smtp1{{\sigma}^{-t+1}}

\def\S1{S^{-1}}

\def\Xpm1{X^{\pm 1}_1}

\def\sPM1{{\sigma }^{\pm 1}}
\def\sMP1{{\sigma }^{\mp 1 }}

\def\di{{\rm d.ind}}

\def\L{\Lambda}

\def\G{\Gamma}

\def\CA{{\cal A}}

\def\CD{{\cal D}}

\def\Ytm1{Y^{t-1}}
\def\Yim1{Y^{i-1}}

\def\CL{{\cal L}}
\def\CM{{\cal M}}
\def\CN{{\cal N}}
\def\CS{{\cal S}}
\def\CF{{\cal F}}
\def\CG{{\cal G}}

\def\ass{{\rm ass}}
\def\bT{\overline{T}}

\def\supp{{\rm supp }}

\def\dim{{\rm dim }}

\def\ker{ {\rm ker } }

\def\SL2Z{ {\rm SL}_2({\bf Z}) }

\def\th{ \theta }

\def\CL{{\cal L}}

\def\Gp1{ G^{1 , 1 } }
\def\P11{ P^{-1 , 1 } }
\def\Pp1{ P^{1 , 1 } }

\def\th{\theta}

\def\nCLsr{{}^\nu\kern-2pt {\cal L}^{\sigma , \rho  }}
\def\nP{{}^\nu \kern-2pt P}
\def\nL{{}^\nu\kern-2pt L}
\def\nLL{{}^\nu\kern-2pt \Lambda}
\def\nPsr{{}^\nu\kern-2pt P^{\sigma , \rho  }}
\def\nLsr{{}^\nu\kern-2pt L^{\sigma , \rho  }}
\def\nuCL{{}^\nu\kern-2pt  {\cal L}}
\def\nCLsr{{}^\nu\kern-2pt {\cal L}^{\sigma , \rho  }}
\def\nCL1m{{}^\nu\kern-2pt {\cal L}^{-1 , 1  }}
\def\x1nu{x^\frac{1}{\nu}}
\def\xm1nu{x^{-\frac{1}{\nu}}}

\def\rad{{\rm rad}}

\def\CN{{\cal N}}
\def\ra{\rightarrow }

\def\CB{{\cal B}}

\def\CI{{\cal I}}

\def\CT{{\cal T}}

\def\CC{ {\cal C}}

\def\CP{ {\cal P}}

\def\nAM0{{\nu }_{{\cal A},M_0}}
\def\nAN0{{\nu }_{{\cal A},N_0}}

\def\End{ {\rm End }}

\def\CP{ {\cal P }}

\def\bR{\overline{R}}

\def\ga{\mathfrak{a}}
\def\gb{\mathfrak{b}}
\def\gc{\mathfrak{c}}

\def\gn{\mathfrak{n}}
\def\gm{\mathfrak{m}}
\def\gp{\mathfrak{p}}
\def\gq{\mathfrak{q}}

\def\Spec{{\rm Spec}}

\def\di!{\frac{\der^i}{i!}}
\def\dik!{\frac{\der^k_i}{k!}}

\def\CC{{\cal C}}

\def\CT{{\cal T}}

\def\CU{{\cal U}}

\def\id{{\rm id}}

\def\N{\mathbb{N}}

\def\0{\overline{0}}
\def\1{\overline{1}}

\def\Ln1{\L_{n,\overline{1}}}

\def\a1{a_{\overline{1}}}

\def\bs{\overline{s}}

\def\S{\Sigma}

\def\CU{{\cal U}}

\def\vn1{\overrightarrow{n-1}}

\def\im{{\rm im}}

\def\mL{\mathbb{L}}

\def\gf{\mathfrak{f}}

\def\Sub{{\rm Sub}}

\def\mJ{\mathbb{J}}
\def\mI{\mathbb{I}}

\def\ann{{\rm ann}}

\def\bM{\overline{M}}

\def\K1{{\rm K}_1}

\def\hmI1{\widehat{\mI_1}}
\def\tmI1{\widetilde{\mI_1}}
\def\tmJ1{\widetilde{\mJ_1}}
\def\hB1{\widehat{B_1}}
\def\hCB1{\widehat{\CB_1}}

\def\bS{\overline{S}}

\def\Den{{\rm Den}}

\def\Ore{{\rm Ore}}

\def\Den{{\rm Den}}

\def\maxDen{{\rm max.Den}}

\def\br{\overline{r}}

\def\gt{\mathfrak{t}}
\def\bs{\overline{s}}

\def\ga{\mathfrak{a}}

\def\udim{{\rm udim}}

\def\gll{\mathfrak{l}}

\def\pCCR{{}'\CC_R}

\def\RSm{ R\langle S^{-1}\rangle}
\def\mLR{ \mathbb{L}(R)}
\def\mLlR{ \mathbb{L}_l(R)}
\def\mLrR{ \mathbb{L}_r(R)}
\def\mLsR{ \mathbb{L}_*(R)}

\def\mL{\mathbb{L}}

\def\pCCR{{}'\CC_R}

\def\Z{\mathbb{Z}}

\def\mL{\mathbb{L}}

\def\CS{{\cal S}}
\def\CL{{\cal L}}

\def \mrL{\mathrm{L}} 
\def\mrLR{ \mathrm{L}(R)}
\def\mrLRa{ \mathrm{L}(R, \ga )}

\def\mrLoc{\mathrm{L}{\rm oc}}

\begin{document}

\author{V. V. \  Bavula %(GenLocSets.tex)
}

\title{General theory of localizations of rings and modules}

\maketitle
\begin{abstract}
 The aim of the paper is to start  to develop the most general theory of localizations/inversion. Several new concepts are introduced and studied. \\

% The concepts of localizable set, localization of a ring and a module at a localizable set are introduced  and studied. Localizable sets are generalization of Ore sets and denominator sets, and the localization of a ring/module at a localizable set is a generalization of localization of a ring/module at a denominator set. 

%***  For a semiprime left Goldie ring, it is proven that the set of  maximal left localizable sets that contain all regular elements is equal to the set of maximal  left denominator sets (and they are explicitly described). For a semiprime Goldie ring, it is  proven that the following five sets coincide: the maximal   Ore sets,  the maximal  denominator sets,  the maximal left or right or two-sided  localizable sets that contain all regular elements (and they are explicitly described).  ***.\\



$${\bf Contents}$$
\begin{enumerate}
\item Introduction.
\item Localization of a ring at a set. 
\item Maximal localizable sets, the localization radical, the sets of localizable and non-localizable elements.
\item Localizations of commutative  rings.
\item Localizations of semiprime left Goldie rings and direct products of division rings. 
\item  Localization of a module at a localizable set.
\end{enumerate}

\end{abstract}

%%%%%%%%%%%%%%%%%% SECTION 1 %%%%%%%%%%%%%%%%%%%%%%%%

\section{Introduction}

Throughout, all rings are associative, with a unit element 1, which is inherited
by subrings, preserved by homomorphisms and which acts unitally on modules. Let $K$ be a commutative ring with 1. In the paper, a ring $R$ means a $K$-{\em algebra} with fixed  ring homomorphism $\nu_R: K\ra R$ where the image $\im (\nu_R)$ belongs to the centre $Z(R)$ of $R$. A ring  homomorphism $f:R\ra R'$ is a $K$-{\em homomorphism}, that is $f\nu_R=\nu_{R'}$.  Ring  homomorphism means a $K$-homomorphism. For a ring $R$, let $R^\times$ be its group of units and  $\gn_R$ be its prime radical.

Inversion of elements in a ring  is an  important and difficult operation that  `simplifies'  the situation as a rule.   Ore's method of localization is an example of a theory of {\em one-sided fractions} where by definition only the elements of a {\em denominator set} can be inverted and the result is a ring of one-sided fractions which always exsists, i.e. it is not equal to zero, \cite{Stenstrom-RingQuot,Jategaonkar-LocNRings,MR}. In \cite[Theorem 4.15]{Bav-intdifline}, it is proven that the elements of an {\em arbitrary} (left and right) {\em Ore set} can be inverted and the result  is also a ring of one-sided fractions (which always exists)  but in general Ore set is not a denominator set. This was the starting point in \cite{LocSets} and \cite{Loc-groupUnits} where the  most general theory of one-sided fractions was presented. The goal of the paper is  to consider the general situation and to start to develop the most general theory of localization/inversion where the result of localization/inversion of elements does not necessarily yields one-sided fractions.

%***************** 

% For that the following new concepts are introduced: the  almost Ore set, the  localizable set and the  localizable perfect set. Their relations are given by the chain of inclusions:
%$$\{ {\rm Denominator \; sets}  \} \subseteq \{ {\rm Ore \; sets}  \} \subseteq \{ {\rm almost\;  Ore\;  sets}  \} \subseteq \{ {\rm perfect\;  localizable\;  sets }  \} \subseteq \{ {\rm  localizable\;  sets}  \} . $$

%*****************

%{\bf The ring $\RSm$ and the ideal $\ass_R(S)$.}  
Let $R\langle X_S\rangle$ be a ring freely generated by the ring $R$ and a set $X_S=\{ x_s\, | \, s\in S\}$ of free noncommutative indeterminates (indexed by the elements of the set $S$). Let $I_S$ be the  ideal of $R\langle X_S\rangle$  generated by the set  $\{ sx_s-1, x_ss-1 \, | \, s\in S\}$. The factor  ring  
%\marginpar{IRbSbm}
\begin{equation}\label{IRbSbm}
\RSm := R\langle X_S\rangle/ I_S.
\end{equation}
 is called the {\em localization of $R$ at} $S$. 
Let $\ass (S) = \ass_R(S)$ be the  kernel of the ring homomorphism
%\marginpar{IRbSbm1}
\begin{equation}\label{IRbSbm1}
\s_S: R\ra \RSm , \;\; r\mapsto r+ I_S.
\end{equation}
  The factor ring $\bR :=R/\ass_R(S)$ is a subring of $\RSm$ and the map $\pi_S:R\ra \bR$, $r\mapsto \br := r+\ass_R(S)$ is   an epimorphism.  The  ideal $\ass_R(S)$ of $R$ has a complex structure, its description is given in \cite[Proposition 2.12]{LocSets} when  $$\RSm =\{ \bs^{-1}\br\, | \, s\in S, r\in R\}$$ is a ring of left fractions.   
 There is an example of a domain $R$ and a finite set $S$ such that $\ass (S)=R$, i.e. $\RSm =\{ 0\}$  \cite[Exercises 9.5]{Lam-Exbook}. 
    Proposition \ref{A12Jan19}.(1) and  its proof  describe an explicit ideal  $\ga (S)$ of $R$ that is contained in $\ass_R(S)$. The ideal $\ga(S)$ is the least ideal such that the  elements of the set $\{ s+\ga (S)\in R/\ga (S)\, | \, s\in S\}$ are non-zero-divisors in the factor ring $R/\ga (S)$, \cite[Proposition 1.1]{LocSets}. The proof of Proposition \ref{A12Jan19} contains an explicit description of the ideal $\ga (S)$. 
 % The ideal $\ga (S)$ is the key part in the definition of perfect localizable sets.
 If $S$ is a left denominator set then the set  $$\ass_l(S):=\{ r\in R\, | \, sr=0\;\;{\rm  for\; some}\;\;s\in S\}$$ is  an   ideal of $R$ and  $\ass_l(S)=\ga (S)=\ass_R(S)$. In general, for a multiplicative set $S$, $\ass_l(S)\subseteq \ga (S)\subseteq \ass_R(S)$ (in general, the inclusions are strict).

 %%%%%%%%%%%%%%%%%%%  Sec 2   %%%%%%%%%%%%%%%%%%%%%
 
In Section \ref{LRLSGRU},  we recall Ore's method of localization. The aim of this section is to introduce several new concepts that are associated with  localization of a ring at a set and to prove some general properties for them (Corollary \ref{a17Dec22}, Lemma \ref{b17Dec22}, Proposition \ref{A18Dec22}, Proposition \ref{c17Dec22}, Corollary \ref{xA12Jan19}): A subset $S$ of $R$ is called  a {\em localizable set} if $\RSm\neq \{ 0\}$; $\mrLR$ is the set of all localizable  sets of $R$ and $ \mrLRa :=\{ S\in \mrLR\, | \, \ass_R(S)=\ga\}$; 
$\mrLoc (R)$ and $\mrLoc (R,\ga )$ are the sets of $R$-isomorphism classes in the sets $$\{ \RSm\, | \, S\in\mrLR\} \;\; {\rm  and}\;\;\{ \RSm\, | \, S\in\mrLRa\},$$ respectively;  $\mrL (R, \ga , A ):=\{ S\in \mrL (R, \ga )\, | \, \RSm \simeq A$,  an $R$-isomorphism (necessarily unique)$\}$. All the above sets are partially ordered sets.

 Proposition \ref{A5Dec22} is the universal property of localization. Proposition \ref{a4Dec22} is a criterion for a subset of a ring to be a localizable set.  Theorem \ref{S10Jan19} and Theorem \ref{B16Jan19} demonstrates some of the results of this section for Ore sets. Direct limits are essential  in proving that the partially ordered sets above   admit maximal elements. Several results on direct limits of rings are proven that are used later in the paper (Theorem \ref{A14Dec22}, Corollary \ref{a7Dec22}, and  Lemma \ref{a21Dec22}). For a ring $R$, the concept of  the {\em absolute quotient ring} $$Q_a(R):=\varinjlim_{S\in \mrLR}\RSm$$ is introduced. Theorem \ref{XA30Dec22} is a description  of the ring $Q_a(R)$ in terms of maximals localization sets of $R$. Corollary \ref{ya30Dec22} presents a sufficient condition for a ring $R$ to have at least two maximal localization sets.

%%%%%%%%%%%%%%%%%%%  Sec 3   %%%%%%%%%%%%%%%%%%%%%

In Section \ref{MAXLOCSETS},    for a ring $R$, the following concepts are introduced: the maximal localizable sets, the localization  radical $\mrL\rad (R)$, the set of localizable elements $\CL\mrL (R)$, the set of non-localizable elements $\CN\CL\mrL (R)$, the set of completely  localizable elements $\CC\mrL (R)$, and the complete localization $$Q_c(R):=R\langle \CC\mrL (R)^{-1}\rangle$$ of $R$. These sets  are  invariant under the action of the automorphism group of $R$. 
 Proposition \ref{A29Dec22} shows that there are tight connections between the sets $\CL \mrLR$, $\CN \CL\mrLR$,  and $ \mrL\rad (R)$. For an arbitrary ring, 
Theorem \ref{SC12Jan19} states that  the set of maximal localizable set is a non-empty set and every localizable set  is a subset of a maximal localizable set. 
 
Proposition \ref{SA19Aug21}.(1), is an explicit description of    the largest element  $ \CS (R, \ga , A)$ of the poset $(\mrL(R, \ga,  A), \subseteq )$. Corollary \ref{a20Dec22} is a criterion for a set  $S\in \mrL (R, \ga, A )$ to be the largest element of the poset $\mrL (R, \ga, A )$.  
Theorem \ref{22Dec22} describes the sets  $\max \mrL (R)$ and $\max \mrLoc (R)$, presents a bijection between them, and shows that $\max \mrLoc (R)\neq \emptyset$.  Similarly, for each ideal $\ga\in \ass \, \mrLR$, Theorem \ref{A23Dec22} describes the sets  $\max \mrLRa$ and $\max \mrLoc (R, \ga)$, presents a bijection between them, and shows that they are nonempty sets. The $\ga$-absolute quotient ring $$Q_a(R, \ga ):=\varinjlim_{S\in  \mrLRa}\RSm$$ is introduced for each ideal $\ga \in \ass \, \mrLR$. Theorem \ref{aXA30Dec22} gives an explicit description of the ring $Q_a(R,\ga)$, a criterion for $Q_a(R,\ga)\neq \{0\}$, and a criterion for $\ga_{a,\ga}=\ga$.
  
  Theorem \ref{20Dec22}.(1) is an $R$-isomorphism  criterion for  localizations of the ring $R$ and Theorem \ref{20Dec22}.(2) is a criterion for existing of an $R$-homomorphism between localizations of $R$. Proposition \ref{B29Dec22}  contains some properties of the set of completely localizable elements $\CC\mrLR$ and the ideals $\gc_R:= \ass_R(\CC\mrLR)$. 
Theorem \ref{A17Jan23} is a criterion for the ring $\RSm$ to be an $R$-isomorphic to a localization of $R$ at a denominator set and  it also  describes all such denominator sets.

For a denominators set $T$ of $R$, Proposition \ref{B17Jan23}  describes all localizable sets $S$ of the ring $R$ such that the ring $\RSm$ is an $R$-isomorphic to the localization of $R$ at $T$. For a localizable  multiplicative   set $S$ of the ring $R$, Proposition \ref{C17Jan23} is a criterion for the set $S$ to be a left denominator set of $R$.

%%%%%%%%%%%%%%%%%%%  Sec 4   %%%%%%%%%%%%%%%%%%%%%

In Section \ref{COMRINGS}, for an arbitrary commutative ring $R$, descriptions of the following sets are obtained: the sets of localizable and non-localizable elements (Lemma \ref{a29Dec22}), the maximal localizable sets, maximal localization rings, the set of completely localizable elements, and the ideal $\gc_R$  (Theorem \ref{C29Dec22}). Theorem \ref{C29Dec22}  also describes relations between the  ideals $\mrL\rad (R)$, $\gn_R$, and $\gc_R$.

For a commutative ring $R$ with finitely many minimal primes  and   nilpotent prime radical (eg, $R$ is a commutative Noetherian ring), Proposition \ref{XC29Dec22} describes the spectrum  $\Spec (Q_c(R))$ and  the rings  $Q_c (R)$ and $Q_c(R)/\gn_{Q_c(R)}$.

%%%%%%%%%%%%%%%%%%%  Sec 5  %%%%%%%%%%%%%%%%%%%%%

In Section \ref{LEFTGOLDIE}, 
 Proposition \ref{A26Dec22}.(4) describes the maximal localizable sets of finite direct product of rings and  their localizations. Proposition \ref{A26Dec22}.(1) describes localizations of the direct product of rings via localizations of its components. For a direct product of simple rings 
 $$A=\prod_{i=1}^sA_i\;\; {\rm such \;that} \;\;\CL\mrL (A_i)=A_i^\times\;\;{\rm  for}\;\;i=1, \ldots , s,$$  
 Theorem  \ref{b26Dec22} describes the sets $\max \mrL (A)$, $\ass_A(S)$ for all $S\in \max \mrL (A)$, the sets of localizable, non-localizable, and completely localizable elements of $A$. It also shows that $\mrL\rad (A)=0$. Every semisimple Artinian ring satisfies the assumptions of Theorem  \ref{b26Dec22}, see Corollary \ref{26Dec22} for detail. For a semiprime left Goldie ring, Lemma \ref{a27Dec22} describes all the  maximal localizable sets that contain the set of regular elements of the ring. Let 
 $$D=\prod_{i\in I}D_i$$ be  a direct product of  division rings $D_i$ where $I$ is an arbitrary set. Proposition \ref{B28Dec22} describes all the localizable sets in the ring $D$ and all the localizations of $D$. Theorem \ref{A28Dec22} describes  the following sets: $\max \mrL (D)$, $\max \mrLoc (D)$, $\mrL\rad (D)$, $\CC\mrL (D)$, $Q_c(D)$, $\CL\mrL (D)$, and 
 $\CN\CL\mrL (D)$. Theorem \ref{A28Dec22}.(2) shows that for every $\CS\in \max \mrL (D)$, the ring $D\langle \CS^{-1}\rangle$ is a division ring.

%%%%%%%%%%%%%%%%%%%  Sec 6  %%%%%%%%%%%%%%%%%%%%%

In  Section \ref{LOCMOD},   basic properties of the  localization of a module at a localizable set  are considered. Proposition \ref{SA20Jan19}.(2) is a universal property of localization of a module. Theorem \ref{S20Jan19} is a criterion for the  localization functor $M\ra S^{-1}M$ to be an exact functor.

%%%%%%%%%%%%%%%%%   Section   2  %%%%%%%%%%%%%%%%%%%

\section{Localization of a ring at a set}\label{LRLSGRU}%\marginpar{LRLSGRU}

At the beginning of the section we recall Ore's method of localization, i.e. the localization of a ring at a left or right denominator set. The aim of this section is to introduce several new concepts that are associated with  localization of a ring at a set and to prove some general properties for them (Corollary \ref{a17Dec22}, Lemma \ref{b17Dec22}, Proposition \ref{A18Dec22}, Proposition \ref{c17Dec22}, Corollary \ref{xA12Jan19}). Proposition \ref{A5Dec22} is the universal property of localization. Proposition \ref{a4Dec22} is a criterion for a subset of a ring to be a localizable set. Every Ore set is a localizable set,  \cite[Theorem 4.15]{Bav-intdifline}. Theorem \ref{S10Jan19} and Theorem \ref{B16Jan19} demonstrates some of the results of this section for Ore sets.

Direct limits play a prominent role in proving that certain  natural posets (partially ordered sets) that are associated with localizations admit maximal elements. Several results on direct limits of rings are proven that are used later in the paper (Theorem \ref{A14Dec22}, Corollary \ref{a7Dec22}, and  Lemma \ref{a21Dec22}).

For a ring $R$, the concept of the  absolute quotient ring $Q_a(R)$ is introduced. Theorem \ref{XA30Dec22} is a description  of the ring $Q_a(R)$ in terms of maximals localization sets of $R$. Corollary \ref{ya30Dec22} presents a sufficient condition for a ring $R$ to have at least two maximal localization sets. \\

{\bf Ore and denominator sets, localization of a ring at a denominator set.} Let $R$ be a ring. A subset $S$ of $R$ is called a {\em multiplicative set} if  $SS\subseteq S$, $1\in S$ and $0\not\in S$. A 
multiplicative subset $S$ of $R$   is called  a {\em
left Ore set} if it satisfies the {\em left Ore condition}: for
each $r\in R$ and
 $s\in S$, $$ Sr\bigcap Rs\neq \emptyset .$$
Let $\Ore_l(R)$ be the set of all left Ore sets of $R$.
  For  $S\in \Ore_l(R)$, $\ass_l (S) :=\{ r\in
R\, | \, sr=0 \;\; {\rm for\;  some}\;\; s\in S\}$  is an ideal of
the ring $R$. 

%$\noindent $

A left Ore set $S$ is called a {\em left denominator set} of the
ring $R$ if $rs=0$ for some elements $ r\in R$ and $s\in S$ implies
$tr=0$ for some element $t\in S$, i.e., $r\in \ass_l (S)$. Let
$\Den_l(R)$ (resp., $\Den_l(R, \ga )$) be the set of all left denominator sets of $R$ (resp., such that $\ass_l(S)=\ga$). For
$S\in \Den_l(R)$, let $$S^{-1}R=\{ s^{-1}r\, | \, s\in S, r\in R\}$$
be the {\em left localization} of the ring $R$ at $S$ (the {\em
left quotient ring} of $R$ at $S$). Let us stress that in Ore's method of localization one can localize {\em precisely} at left denominator sets.
 In a similar way, right Ore and right denominator sets are defined. 
Let $\Ore_r(R)$ and 
$\Den_r(R)$ be the set of all right  Ore and  right   denominator sets of $R$, respectively.  For $S\in \Ore_r(R)$, the set  $\ass_r(S):=\{ r\in R\, | \, rs=0$ for some $s\in S\}$ is an ideal of $R$. For
$S\in \Den_r(R)$,  $$RS^{-1}=\{ rs^{-1}\, | \, s\in S, r\in R\}$$   is the {\em right localization} of the ring $R$ at $S$. 

Given ring homomorphisms $\nu_A: R\ra A$ and $\nu_B :R\ra B$. A ring homomorphism $f:A\ra B$ is called an $R$-{\em homomorphism} if $\nu_B= f\nu_A$.  A left and right  set is called an {\em Ore set}.  Let $\Ore (R)$ and 
$\Den (R)$ be the set of all   Ore and    denominator sets of $R$, respectively. For
$S\in \Den (R)$, $$S^{-1}R\simeq RS^{-1}$$ (an $R$-isomorphism)
 is  the {\em  localization} of the ring $R$ at $S$, and $\ass (S):=\ass_l(S) = \ass_r(S)$. \\ 

{\bf The ring $\RSm$ and the ideal $\ass_R(S)$.} Let $R$ be a ring and $S$ be a  subset of  $R$. Let $R\langle X_S\rangle$ be a ring freely generated by the ring $R$ and a set $X_S=\{ x_s\, | \, s\in S\}$ of free noncommutative indeterminates (indexed by the elements of the set $S$). Let $I_S$ be the  ideal of $R\langle X_S\rangle$  generated by the set  $\{ sx_s-1, x_ss-1 \, | \, s\in S\}$ and 
%\marginpar{RbSbm}
\begin{equation}\label{RbSbm}
\RSm := R\langle X_S\rangle/ I_S.
\end{equation}
The ring $\RSm$ is called the {\em localization of $R$ at} $S$. 
Let $\ass (S) = \ass_R(S)$ be the  kernel of the ring homomorphism
%\marginpar{RbSbm1}
\begin{equation}\label{RbSbm1}
\s_S: R\ra \RSm , \;\; r\mapsto r+ I_S.
\end{equation}
 The map $\pi_S:R\ra \bR:= R/ \ass_R(S)$, $r\mapsto \br := r+\ass_R(S)$ is   an epimorphism.  The  ideal $\ass_R(S)$ of $R$ has a complex structure, its description is given in \cite[Proposition 2.12]{LocSets} when  $\RSm =\{ \bs^{-1}\br\, | \, s\in S, r\in R\}$ is a ring of left fractions.  We identify the factor ring $\bR $ with its isomorphic copy in the ring $\RSm$  via the monomorphism
%\marginpar{bRS}
\begin{equation}\label{bRS}
\overline{\s_S}: \bR\ra \RSm , \;\; r+\ass_R(S)\mapsto r+ I_S. 
\end{equation}
Clearly, $\bS := (S+\ass_R(S) ) /\ass_R(S)\subseteq \CC_{\RSm}$. Corollary \ref{a17Dec22}.(2) shows that the rings $\RSm$ and $\bR\langle \bS^{-1}\rangle$ are $R$-isomorphic. For $S=\emptyset$,  $R\langle \emptyset^{-1}\rangle := R$ and $\ass_R(\emptyset ):=0$. 

{\it Definition.} A subset $S$ of a ring $R$ is called a {\em localizable set} of $R$ if  
$\RSm \neq \{ 0\}$. Let $\mrLR$ be the set of localizable sets of $R$ and 
%\marginpar{SxLzRa}
\begin{equation}\label{SxLzRa}
 \ass \, \mrL (R) :=\{ \ass_R (S) \, | \, S\in \mrL (R) \}.
\end{equation}

Proposition \ref{A5Dec22} is the {\em universal property of  localization}. 

\begin{proposition}\label{A5Dec22}%\marginpar{A5Dec22}
Let $R$ be a ring,  $S\in \mrLR$, and $\s_S : R\ra  R\langle S^{-1}\rangle$, $r\mapsto r+\ass_R(S)$. Let $f: R\ra A$ be a ring homomorphism such that $f(S)\subseteq A^\times$. The there is a unique $R$-homomorphism $f':R\langle S^{-1}\rangle \ra A$ such  $f=f'\s_S$, i.e. the diagram below is commutative 
$$\begin{array}[c]{ccc}
R&\stackrel{\s_S}{\ra} &\RSm\\
&\stackrel{f}{\searrow}&{\downarrow}^{\exists !\, f'}\\
&  & A
\end{array}$$
\end{proposition}
{\it Proof}. The unique $R$-homomorphism $f'$   is defined by the rule 
$$f':R\langle S^{-1}\rangle =R\langle X_S \rangle/I_X \ra A, \;\; r+I_S\mapsto f(r), \;\; x_s +I_S\mapsto f(s)^{-1}$$ for all $r\in R$ and $s\in  S$. 
 $\Box $

 \begin{corollary}\label{a17Dec22}%\marginpar{a17Dec22}
 Suppose that $S\in \mrL (R, \ga )$, $\gb$ is an ideal of $R$ such that $\gb \subseteq \ga$, $\bS:=\pi_\gb (S)$ where $\pi_\gb : R\ra R/\gb$, $r\mapsto r+\gb$. Then
\begin{enumerate}
\item  $\bS \in \mrL (R/\gb , \ga / \gb)$ and $(R/\gb ) \langle \bS^{-1}\rangle\simeq R\langle S^{-1}\rangle$, an $R$-isomorphism.
\item In particular, for $\ga = \gb$,  $\bS \in \mrL (R/\ga , 0)$ and $(R/\ga ) \langle \bS^{-1}\rangle\simeq R\langle S^{-1}\rangle$, an $R$-isomorphism.
\end{enumerate}
\end{corollary}

{\it Proof}. 1.  By the universal property of the localization, there is a unique $R$-homomorphism $R\langle S^{-1}\rangle\ra (R/\gb ) \langle \bS^{-1}\rangle$ (since $\gb \subseteq \ga$). Conversely, since $\gb \subseteq \ga$,  there is a unique $R$-homomorphism $(R/\gb )\langle \bS^{-1}\rangle\ra  R\langle S^{-1}\rangle$, by the universal property of the localization, and the lemma follows.

2. Statement 2 is a particular case of statement 1 where $\gb = \ga$.   $\Box $

Corollary \ref{a4Jan23} is a criterion for two localizations of a ring $R$ to be $R$-isomorphic.

\begin{corollary}\label{a4Jan23}%\marginpar{a4Jan23}
 Let $S\in \mrLRa$, $T\in \mrL (R, \gb)$, $\bS =\pi (S)$ and $\bT=\pi (T)$ where $\pi :R\ra \bR =R/\ga$, $r\mapsto r+\ga$. Then the rings 
 $\RSm$ and $R\langle T^{-1}\rangle$ are $R$-isomorphic iff $\ga = \gb$, $\bS\subseteq  R\langle T^{-1}\rangle^\times$, and $\bT\subseteq \RSm^\times$.
\end{corollary}

{\it Proof}. $(\Rightarrow )$ The implication is obvious.

$(\Leftarrow )$ Suppose that  $\ga = \gb$, $\bS\subseteq  R\langle T^{-1}\rangle^\times$, and $\bT\subseteq \RSm^\times$. By Proposition \ref{A5Dec22}, we have  commutative diagrams of $R$-homomorphisms, 
$$\begin{array}[c]{ccc}
R&\stackrel{\s_S}{\ra} &\RSm\\
&\stackrel{\s_T}{\searrow}&{\downarrow}^{\exists !\, f}\\
&  & R\langle T^{-1}\rangle
\end{array}\;\;\; {\rm and}\;\;\;
\begin{array}[c]{ccc}
R&\stackrel{\s_T}{\ra}& R\langle T^{-1}\rangle\\
&\stackrel{\s_S}{\searrow}&{\downarrow}^{\exists !\, g}\\
&  & \RSm .
\end{array}
$$ 
By Proposition \ref{A5Dec22}, $fg=\id_{R\langle T^{-1}\rangle}$ and $gt=\id_{\RSm}$, and so $f$ is an $R$-isomorphism.  $\Box $

%****************   dor-t    *********************

%Corollary \ref{b4Jan23} is an isomorphism  criterion  for  localizations of a ring.

%\begin{corollary}\label{b4Jan23}%\marginpar{b4Jan23}
% Let $S\in \mrLRa$, $T\in \mrL (R, \gb)$, $\bS =\pi (S)$ and $\bT=\pi (T)$ where $\pi :R\ra \bR =R/\ga$, $r\mapsto r+\ga$. Then the rings 
% $\RSm$ and $R\langle T^{-1}\rangle$ are isomorphic iff $\ga = \gb$ and there are homomorphisms $f:\RSm\ra R\langle T^{-1}\rangle$ and $g: R\langle T^{-1}\rangle\ra \RSm$ such that $fg(\bR)=\bR$, $gf(\bR)=\bR$, $f(\bS)\subseteq R\langle T^{-1}\rangle^\times$, and $g(\bR)\subseteq \RSm^\times$.
%\end{corollary}

%{\it Proof}. $(\Rightarrow )$ Suppose that $f:\RSm\ra R\langle T^{-1}\rangle$ is an isomorphism. Then clearly $\ga =\gb$, and  it suffices to take $g=f^{-1}$. 

%$(\Leftarrow )$ 

%$\Box $

%********************************************

 If $S\subseteq T\subseteq R$ then $\ass_R(S)\subseteq \ass_R(T)$ and, by Proposition \ref{A5Dec22},  there is  a {\em unique} $R$-homomorphism
%\marginpar{RSmRT}
\begin{equation}\label{RSmRT}
\phi_{TS}: R\langle S^{-1}\rangle \ra R\langle T^{-1}\rangle,
\end{equation}
i.e. $\phi_{TS}(ra)=r\phi_{TS}(a)$ for all elements $a\in  R\langle S^{-1}\rangle $.  Clearly, $\phi_{SS}=\id_{R\langle S^{-1}\rangle}$ and 
%\marginpar{RSmRT2}
\begin{equation}\label{RSmRT2}
\ker (\phi_{TS})\supseteq  R\langle S^{-1}\rangle \Big( \ass_R(T)/\ass_R(S) \Big)R\langle S^{-1}\rangle.
\end{equation}
If $S\subseteq T\subseteq  U\subseteq R$ then
%\marginpar{RSmRT1}
\begin{equation}\label{RSmRT1}
\phi_{US}=\phi_{UT}\phi_{TS}.
\end{equation}

{\bf Criterion for a set to be a localizable set.} Proposition \ref{a4Dec22} is a   criterion for a set to be a localizable set.

\begin{proposition}\label{a4Dec22}%\marginpar{a4Dec22}
Let $S$ be a subset of $R$. The following statements are equivalent:
\begin{enumerate}
\item $S\in \mrLR$. 
\item $\ass_R(S)\neq R$.
\item There is an ring homomorphism $f: R\ra A$ such that $f(S)\subseteq A^\times$ where $A^\times$ is a group of units of the ring $A$. 
\end{enumerate} 
\end{proposition}

{\it Proof}. $(1\Leftrightarrow 2)$ $S\in \mrLR $ iff $\RSm \neq \{ 0\}$ iff the image of the ring $R$  under the $R$-homomorphism $R\ra \RSm$ is not equal to $0$ iff $\ass_R(S)\neq R$.

$(1\Rightarrow 3)$ If $S\in \mrLR$ then the image of the set $S$ under the homomorphism $\s_S: R\ra \RSm$ consists of units of the ring $ \RSm$. 

$(3\Rightarrow 2)$ Suppose that there is a ring homomorphism  $f: R\ra A$ such that $f(S)\subseteq A^\times$. By the universal property of  localization, there is a unique $R$-homomorphism $\RSm\ra A$. Hence, $\RSm\neq \{0\}$ (since $f(S)\subseteq A^\times$), and so $\ass_R(S)\neq R$.   $\Box $ 

 Proposition \ref{D17Jan23} gives  a sufficient condition for a subset of a ring to be a localizable set.

\begin{proposition}\label{D17Jan23}%\marginpar{D17Jan23}
Let $I$ be an ideal of a ring $R$, $\pi : R\ra R/I$, $r\mapsto r+I$,  and $S$ be a subset of $R$ such that $\pi (S)\subseteq \CC_{R/I}$ and $\ass_R(S)=\ga (S)$. Then $S\in \mrL (R, \ga (S))$ and $\ga (S)\subseteq I$.
\end{proposition}

{\it Proof}. By Corollary \ref{xA12Jan19}.(1), $\ga (S) \subseteq I$ since by induction $\ga_\l \subseteq I$ for all $\l \in \N$.
% By ***  Theorem \ref{5Jan23}  ***, 
Since  $\ga (S)=\ass_R(S)$, the set $S$ is a localizable set of $R$ since $\ass_R(S)\neq R$ (Proposition \ref{a4Dec22}) as $\ass_R(S)\subseteq I$.  $\Box $

Let $S$ be a subset of the ring $R$ and $S_{mon}$ be the multiplicative submonoid of $(R, \cdot )$ generated by the set $S$. Lemma \ref{b17Dec22} shows that the localization of the ring $R$ at $S$ is the same as the localization of the ring $R$ at $S_{mon}$.

\begin{lemma}\label{b17Dec22}%\marginpar{b17Dec22}
Let $S$ be a subset of the ring $R$.
% and $S_{mon}$ be the multiplicative submonoid of $(R, \cdot )$ generated by the set $S$. 
Then
\begin{enumerate}
\item $S\in \mrLRa$ iff $S_{mon}\in \mrLRa$. 
\item If $S\in \mrLRa$ then $\RSm\simeq R\langle S_{mon}^{-1}\rangle$, an $R$-isomorphism  
\end{enumerate}
\end{lemma}

{\it Proof}. 1. If $S\in \mrLR$ then the image of the monoid $S_{mon}$ under the $R$-homomorphism $\s_S: R\ra \RSm$ consists of units, and so $S_{mon}\in \mrLR$, by Proposition \ref{a4Dec22}. Similarly, if $S_{mon}\in \mrLR$ then the image of the set $S$ under the $R$-homomorphism $\s_{S_{mon}}: R\ra R\langle S_{mon}^{-1}\rangle$ consists of units, and so $S\in \mrLR$, by Proposition  \ref{a4Dec22}.

2.  Since  $S\in \mrLR$,   $S_{mon}\in \mrLR$,  by  statement 1. By the universal property of localization,  there is a unique $R$-homomorphism $\RSm\ra R\langle S_{mon}^{-1}\rangle$.  By the universal property of localization,  there is a unique $R$-homomorphism  $ R\langle S_{mon}^{-1}\rangle\ra \RSm$. So, the rings  $\RSm$ and $ R\langle S_{mon}^{-1}\rangle$ are  an $R$-isomorphic, and statement 2 follows.  $\Box $\\

 {\bf The set $\CM (S,R)$ of ideals of $R$  and its minimal element.} For a set $S\in  \mrLR$, let $\CM (S, R)$ be the set of ideal $\ga$ of the ring $R$ such that there is a ring homomorphism $\v :R/\ga\ra A$ for some ring $A$ such that $\v (S+\ga ) \subseteq A^\times$. Proposition \ref{A18Dec22} shows that $\ass_R(S)$ is the least element of the set $\CM (S, R)$ w.r.t. inclusion. 
 
 \begin{proposition}\label{A18Dec22}%\marginpar{A18Dec22}
For every set $S\in  \mrLR$, $\min \CM (S, R)=\{ \ass_R(S)\}$. 
\end{proposition}

{\it Proof}. Let $\ga \in \CM (S, R)$ and  $$f: R\ra R/\ga \stackrel{\v}{\ra} A, \;\;r\mapsto r+\ga \mapsto \v(r+\ga )$$ for some ring $A$ such that $\v (S+\ga ) \subseteq A^\times$. Clearly, $\ker (f)=\ga$. By Proposition \ref{A5Dec22}, there is                                                                                                                                                                                                                                                                                                                                                                                                                                                                                                                                                                                                                                                                                                                                                                                                                                                                                                                                                                                                                                                                                                                                                                                                                                                                                                                                                                                                                                                                                                                                                                                                                                                                                                                                                                                                                                                                                                                                                                                                                                                                                                                                                                                                                                                                                                                                                                                                                                                                                                                                    an $R$-homomorphism $f':\RSm\ra A$ such that $f=f'\s_S$. Therefore, $$ \ass_R(S)=\ker (\s_S)\subseteq \ker (f'\s_S)=\ker (f)=\ga,  $$ 
 and the proposition follows. $\Box$\\
 
 {\bf General results on localizations.} Some useful general results on localizations are collected below that are used in the proofs of this paper.

\begin{proposition}\label{A7Dec22}%\marginpar{A7Dec22}

\begin{enumerate}
\item If $S\subseteq R^\times$ then $\ass_R(S)=0$ and $\RSm =R$.
\item If $T\in  \mrLR$ and $S\subseteq T$. Then $S\in  \mrLR$, $\ass_R(S)\subseteq \ass_R(T)$, 
 $\s_S(T)\in \mrL (\RSm)$,  $$R \langle T^{-1}\rangle \simeq \RSm \langle \s_S(T)^{-1}\rangle,$$ an $R$-isomorphism, and $\ass_{R \langle S^{-1}\rangle}(\s_S(T))\supseteq R \langle S^{-1}\rangle\Big( \ass_R(T)/\ass_R(S)\Big) R \langle S^{-1}\rangle$. 
 \item If $S,T\in  \mrLR$ such that $S\cup T\in \mrLR$ then $\s_S(T)\in \mrL (\RSm)$ and $R \langle (S\cup T)^{-1}\rangle \simeq \RSm \langle \s_S(T)^{-1}\rangle$, an $R$-isomorphism. 
\end{enumerate}
\end{proposition}

{\it Proof}. 1. Both rings  $\RSm$ and $R$  satisfy the universal property of localization,  Proposition \ref{A5Dec22}. Therefore, $\RSm =R$. Hence, $\ass_R(S)=0$. 

2. Since   $T\in  \mrLR$ and $S\subseteq T$, there a unique $R$-homomorphism $\RSm\ra R\langle T^{-1}\rangle$. Therefore,  $S\in  \mrLR$. Both rings $R \langle T^{-1}\rangle $ and $\RSm \langle \s_S(T)^{-1}\rangle$ satisfy the universal property of localization for the set $T$ (Proposition \ref{A5Dec22}), and so they are  $R$-isomorphic. It follows that $\s_S(T)\in \mrL (\RSm)$, by statement 2 (since $T\subseteq S\cup T$).

3.  Statement 3 is a particular case of statement 2 (since $S\subseteq S\cup T$ and $S\cup T\in \mrL (R)$. $\Box $

\begin{lemma}\label{c17Dec22}%\marginpar{c17Dec22}
Suppose that $S\in \mrLRa$, $T\subseteq R$,  and $\s_S :R\ra\RSm$.
\begin{enumerate}
\item If  $\s_S (T)\in \mrL (\RSm )$ then $S\cup T\in \mrLR$ and $R\langle (S\cup T)^{-1}\rangle \simeq \RSm \langle \s_S (T)^{-1}\rangle$, an $R$-isomorphism.
\item If  $\s_S (T)\subseteq \RSm^\times $ then $S\cup T\in \mrLRa$ and $R\langle (S\cup T)^{-1}\rangle \simeq \RSm $, an $R$-isomorphism.
\item Let $\CT =\s^{-1}(\RSm^\times)$. Then $\CT\in \mrLRa$ and $R\langle \CT^{-1}\rangle \simeq \RSm $, an $R$-isomorphism.
\end{enumerate}
 \end{lemma}

{\it Proof}. 1. The image of the set $S\cup T$ in $\RSm \langle \s_S (T)^{-1}\rangle$ consists of units, hence $S\cup T\in \mrLR$, by Proposition \ref{a4Dec22}. By Proposition \ref{A7Dec22}.(3), the rings $R\langle (S\cup T)^{-1}\rangle $ and $ \RSm \langle\s_S (T)^{-1}\rangle$  are $R$-isomorphic. 

2.  Since  $\s_S (T)\subseteq  \RSm^\times $,  $\s_S (T)\in \mrL (\RSm )$. Now, by statement 1, $S\cup T\in \mrLR$ and $$R\langle (S\cup T)^{-1}\rangle \simeq \RSm \langle \s_S (T)^{-1}\rangle\simeq \RSm ,  $$ $R$-isomorphisms where the second $R$-isomorphism is due to Proposition \ref{A7Dec22}.(1).  Hence, $\ass_R(S\cup T)=\ga$.

3. Statement 3 follows from statement 2 where $T=\CT$ since $S\cup \CT=\CT$.  $\Box$ \\

 {\bf The ideals $\ga (S)$, ${}'\ga (S)$ and $\ga'(S)$.} The ideals $\ga (S)$, ${}'\ga (S)$ and $\ga'(S)$ are contained in the ideal $\ass_R(S)$. They are defined iteratively and it is difficult to compute  them in general. They reveal a complex structure of the ideal $\ass_R(S)$.
 
 For each element $r\in R$, let $r\cdot : R\ra R$, $x\mapsto rx$ and $\cdot r : R\ra R$, $x\mapsto xr$. The sets $${}'\mathcal{C}_R := \{ r\in R\, | \, \ker (\cdot r)=0\}\;\; {\rm  and}\;\;\mathcal{C}_R' := \{ r\in R\, | \, \ker (r\cdot )=0\}$$ are called the {\em sets of left and right regular elements} of $R$, respectively.  Their intersection $$\mathcal{C}_R={}'\mathcal{C}_R \cap \mathcal{C}_R'$$ is the  {\em set of regular elements}  of $R$. The rings $$Q_{l,cl}(R):= \mathcal{C}_R^{-1}R\;\; {\rm  and}\;\;Q_{r,cl}(R):= R\mathcal{C}_R^{-1}$$ are called the {\em classical left and right quotient rings} of $R$, respectively. If both rings exist then they are isomorphic and the ring
 $$Q_{cl}(R):= Q_{l, cl}(R)\simeq  Q_{r, cl}(R)$$
 is called the {\em classical quotient ring} of $R$.   
 Goldie's Theorem states that the ring $Q_{l, cl}(R)$ is  a semisimple Artinian ring iff the ring $R$ is a semiprime left Goldie ring (the ring $R$ is called a {\em left Goldie ring} if   $\udim (R)<\infty$ and the ring $R$ satisfies the a.c.c. on left annihilators where $\udim$ stands for the {\em left uniform dimension}). In a similar way a {\em right Goldie ring} is defined. A left and right Goldie ring is called a {\em Goldie ring}.

\begin{proposition}\label{A12Jan19}%\marginpar{A12Jan19}
(\cite[Proposition 1.1]{LocSets}) Let $R$ be a ring and $S$ be a non-empty subset of $R$.
\begin{enumerate}
\item  There is  the least ideal, $\ga= \ga (S)$, such that  $(S+\ga ) / \ga \subseteq \mathcal{C}_{R/ \gb}$.
\item There is  the least ideal,  ${}'\ga= {}'\ga (S)$, such that  $(S+{}'\ga ) / {}'\ga \subseteq {}'\mathcal{C}_{R/ \gb}$.
\item There is  the least ideal, $\ga'= \ga' (S)$, such that  $(S+\ga' ) / \ga' \subseteq \mathcal{C}_{R/ \gb}'$.
\end{enumerate}
\end{proposition}
 
Notice that ${}'\ga (S) \subseteq \ga (S)$ and $\ga' (S)\subseteq \ga (S)$.

Recall that   $S_{mon}$ is the multiplicative submonoid of $(R, \cdot)$ which is generated by the set $S$. If $0\in  S_{mon}$ then $\ga (S)= {}'\ga (S)= \ga' (S)=R$ and $\RSm=0$. So, we can assume that $0\not\in S_{mon}$, i.e. the set $S$ is a multiplicative set. Notice that (Lemma \ref{b17Dec22}) 
%\marginpar{SSmon}
\begin{equation}\label{SSmon}
\ga (S)=\ga (S_{mon}), \;\; \ass (S)= \ass (S_{mon}), \;\; {\rm and}\;\; \RSm=R\langle S_{mon}^{-1}\rangle. 
\end{equation}

  The proof of Proposition \ref{A12Jan19} (\cite[Proposition 1.1]{LocSets}) contains an explicit description of the ideals $\ga= \ga (S)$,  ${}'\ga= {}'\ga (S)$, and $\ga'= \ga' (S)$. In more detail,  let $\G$ be the set of ordinals. The ideal  $\ga$ (resp., ${}'\ga$, $\ga'$) is the union 
%\marginpar{aappa}
\begin{equation}\label{aappa}
\ga= \bigcup_{\l \in \G}\ga_\l\;\;\;   ({\rm resp.,}\;\;  {}'\ga = \bigcup_{\l \in \G}{}'\ga_\l, \;\; \ga' = \bigcup_{\l \in \G}\ga_\l') 
\end{equation}
of ascending chain of ideals $\{ \ga_\l \}_{\l \in \G}$ (resp., $\{ {}'\ga_\l  \}_{\l \in \G}$, $\{ \ga'_\l \}_{\l \in \G}$), where $\l \leq \mu$ in $\G$ implies $\ga_\l \subseteq \ga_\mu$ (resp., ${}'\ga_\l \subseteq {}'\ga_\mu$, $\ga_\l' \subseteq \ga_\mu'$). The ideals $\ga_\l$ (resp., ${}'\ga_\l$, $\ga_\l'$) are defined inductively as follows: the ideal $\ga_0$ (resp.,  ${}'\ga_0$,  $\ga_0'$) is generated by the set $\{ r\in R\,  | \, sr=0$ or $ rt=0$ for some elements $s,t\in S\}$ (resp., $\{ r\in R\,  | \,  rt=0$ for some element $t\in S\}$, $\{ r\in R\,  | \, sr=0$ for some element $s\in S\}$), and  for $\l \in \G$ such that $\l >0$ (where below $\Big( \{ \ldots \} \Big)$ means the ideal of $R$ generated by the set  $ \{ \ldots \} $), 
%\marginpar{ala}
\begin{equation}\label{ala}
\ga_\l = \begin{cases}
 \bigcup_{\mu <\l \in \G}\ga_\mu & \text{if $\l$ is a limit ordinal},\\
\Big( \{ r\in R\, | \, sr\in \ga_{\l -1}\; {\rm or}\; rt\in \ga_{\l -1} \; {\rm for \; some} \; s,t\in S\}\Big) & \text{if if $\l$ is not a limit ordinal}.\\
\end{cases}
\end{equation}
(resp., 
%\marginpar{ala1}
\begin{equation}\label{ala1}
{}'\ga_\l = \begin{cases}
 \bigcup_{\mu <\l \in \G}{}'\ga_\mu & \text{if $\l$ is a limit ordinal},\\
\Big( \{ r\in R\, | \, rt\in {}'\ga_{\l -1}\; {\rm for \; some} \; t\in S\}\Big) & \text{if if $\l$ is not a limit ordinal}.\\
\end{cases}
\end{equation}

%\marginpar{ala2}
\begin{equation}\label{ala2}
\ga_\l' = \begin{cases}
 \bigcup_{\mu <\l \in \G}\ga_\mu' & \text{if $\l$ is a limit ordinal},\\
\Big( \{ r\in R\, | \, sr\in \ga_{\l -1}'\; {\rm for \; some} \; s\in S\}\Big)  & \text{if if $\l$ is not a limit ordinal}).
\end{cases}
\end{equation}
Let us define an   ideal $\gb$ of $R$ which  is the  union $\gb =\bigcup_{\l \in \G}\gb_\l$ of an ascending chain of ideals $\gb_\l$ that are defined  inductively 
%\marginpar{bla}
\begin{equation}\label{bla}
\gb_\l = \begin{cases}
 \bigcup_{\mu <\l \in \G}\gb_\mu & \text{if $\l$ is a limit ordinal},\\
\Big( \{ r\in R\, | \, sr\in \gb_{\l -1}\; {\rm or}\; rt\in \gb_{\l -1}  \; {\rm or}\; srt\in \gb_{\l -1} \;{\rm for \; some} \; s,t\in S\}\Big) & \text{if if $\l$ is not a limit ordinal}.\\
\end{cases}
\end{equation}
where  the ideal $\gb_0$  is generated by the set $\{ r\in R\,  | \, sr=0$ or $ rt=0$ or $srt=0$ for some elements $s,t\in S\}$. 
\begin{corollary}\label{xA12Jan19}%\marginpar{xA12Jan19}
 Let $R$ be a ring and $S$ be a non-empty subset of $R$. Then
\begin{enumerate}
\item  $\ga (S)=\bigcup_{\l\in \N}\ga_\l =\bigcup_{\l\in \N}\gb_\l $.
\item  ${}'\ga (S)=\bigcup_{\l\in \N}{}'\ga_\l$. 
\item $\ga' (S)=\bigcup_{\l\in \N}\ga'_\l$. 
\end{enumerate}
\end{corollary} 

{\it Proof}. The corollary follows from the definitions of the ideals $\ga_\N$, $\gb_\N$, ${}'\ga_\N$, and $\ga_\N'$. 
 In more detail, let us show that  $\ga (S)=\gb_\N :=\bigcup_{\l\in \N}\gb_\l $. The other cases can be treated in a similar way.  The inclusion  $\ga (S)\supseteq \gb_\N $ is obvious. To prove that the reverse inclusion holds it suffices to show that the maps $$s\cdot : R/\gb_\N\ra R/\gb_\N, \;\; r\mapsto sr+\gb_\N\;\; {\rm  and}\;\;\cdot s: R/\gb_\N\ra R/\gb_\N, \;\;r\mapsto rs+\gb_\N$$  are injections  for all $s\in S$. Suppose that $ r+\gb_\N\in \ker (s\cdot)$ or $ r+\gb_\N\in \ker (\cdot s)$, that is $sr\in \gb_\N$ or $rs\in \gb_\N$. Then $sr\in \gb_\l$ or $rs\in \gb_\l$ for some $\l \in \N$, and so $r\in \gb_{\l+1}$ in both cases. This means that $r\in \gb_\N$, as required.  $\Box$\

 {\it Definition.}  The first natural number such that the ascending chain of ideals in (\ref{bla})  (resp., (\ref{aappa})) stabilizes is called the $\G$-{\em length} of the ideal denoted $l_\G (\ga )$ (resp., $l_\G ({}'\ga )$, $l_\G (\ga' )$). If there is no such a number then $l_\G (\ga ):=\infty $ (resp., $l_\G ({}'\ga ):=\infty$, $l_\G (\ga' ):=\infty$).

Lemma \ref{a6Jan23} shows that the $\G$-length $l_\G (\ga )$ can be any natural number or $\infty$.

Let $K$ be a field. Consider the following rings
\begin{enumerate}
\item $R_0=K\langle x,y_0\, | \, xy_0=0\rangle$,
\item $R_n=K\langle x,y_i, u_j,v_j\, | \,  xy_0=0,xy_j=u_jy_{j-1}v_j, i=0, \ldots , n, j=1, \ldots , n\rangle$,
\item $R_\infty =K\langle x,y_i, u_j,v_j\, | \,  xy_0=0,xy_j=u_jy_{j-1}v_j, i\geq 0, j\geq 1\rangle$.
\end{enumerate}
 
Notice that $R_0/(y_0)\simeq K[x]$ and $R_n/(y_0, \ldots , y_n)\simeq K\langle x, u_j,v_j\, | \,  j=1, \ldots , n\rangle$ is a free algebra for all $n\geq 1$ and $n=\infty$.

\begin{lemma}\label{a6Jan23}%\marginpar{a6Jan23}
Let $R_n$, $n\in \N \cup \{ \infty\}$, be the rings above,  $S=\{x^i\, | \, i\geq 0\}$, and $\ga (S, R_n) :=\ga (S)$ where $S\subseteq R_n$. Then  $\ga (S, R_n)=\ass_{R_n}(S)=(y_0, \ldots , y_n)$, $l_\G (\ga (S, R_n))=n$, and
$$R_n\langle S^{-1}\rangle \simeq \begin{cases}
K[x,x^{-1}]& \text{if }n=0,\\
K\langle x^{\pm 1}, u_j,v_j\, | \,  j=1, \ldots , n\rangle& \text{if } n=1, \ldots , \infty.\\
\end{cases}$$
\end{lemma}

{\it Proof}. Suppose that $n=0$. Then $y_0\in \gb_0\subseteq \gb \subseteq \ga (S, R_0)\subseteq \ass_{R_0}(S)$, $R_0/(y_0)=K[x]$, and   so $R_0\langle S^{-1}\rangle \simeq K[x,x^{-1}]$ and  $\ga (S, R_0)=\ass_{R_0}(S)=(y_0)=\gb_0$. 
 Hence,  $l_\G (\ga (S, R_0))=0$.

 Suppose that $n\geq 1$. Then $$R_n=\bigoplus_{w\in W_n}Kw$$ where $W_n$ is the set of all words in the alphabet $\{x,y_i, u_j,v_j\, | \,  i=0, \ldots , n, j=1, \ldots , n\}$ that do not contained   the subwords $xy_i$, where $i=0, \ldots , n$.
 Then $\gb_0 =(y_0)$ since $\ker_{R_n}(x\cdot)=y_0R_n$ and $\ker_{R_n}(\cdot x)=0$. Using the same argument in the case of the  factor ring
 $$R_n/(\gb_0)=K\langle x,y_i, u_j,v_j\, | \,  xy_1=0,xy_k=u_ky_{k-1}v_k, k=2, \ldots , n\rangle \;\; {\rm where}\;\;  i=1, \ldots , n, j=1, \ldots , n,$$
 we conclude that $\gb_1=(y_0, y_1)$ and 
 $$R_n/(\gb_1)=K\langle x,y_i, u_j,v_j\, | \,  xy_2=0,xy_k=u_ky_{k-1}v_k, k=3, \ldots , n\rangle\;\; {\rm where}\;\;  i=2, \ldots , n, j=1, \ldots , n.$$
 Repeating the same argument again and again (or use induction), we prove that for all $m=1, \ldots , n$, 
 $\gb_m=(y_0,\ldots, y_m)$ and 
 $$R_n/(\gb_m)=K\langle x,y_i, u_j,v_j\, | \,  xy_{m+1}=0,xy_k=u_ky_{k-1}v_k, k=m+2, \ldots , n\rangle , i=m+1, \ldots , n, j=1, \ldots , n.$$
Hence, $$R_n\langle S^{-1}\rangle \simeq (R_n/\gb_n)\langle S^{-1}\rangle \simeq K\langle x^{\pm 1}, u_j,v_j\, | \,  j=1, \ldots , n\rangle$$ and  $\ga (S, R_n)=\ass_{R_n}(S)=\gb_n$. 
 It follows that  $l_\G (\ga (S, R_n))=n$.  $\Box $\\

{\bf Perfect localization set and perfect localizations.} 
 By \cite[Lemma 1.2]{LocSets}, 
%\marginpar{gasas}
\begin{equation}\label{gasas}
\ga (S)\subseteq \ass_R (S). 
\end{equation}

{\it Definition, \cite{LocSets}.} If the equality holds, $\ga (S)= \ass_R (S)$, then set $S$ is called a {\em perfect localization  set} and the localization $\RSm$ is called a {\em perfect localization}.

 Theorem \ref{8May23} shows that there are plenty of perfect localizations (especially in applications). 
 
 \begin{theorem}\label{8May23}%\marginpar{8May23}
Let $R$ be a $K$-algebra over a field $K$, $\ga :=\ga (S)$, $\pi_\ga :R\ra \bR :=R/\ga$, $r\mapsto \br := r+\ga$, $S$ be a multiplicative set in $R$. Suppose that $\bR$ is a domain and for each  pair of elements $s\in S$ and $r\in R$, the subalgebra  $K\langle \bs, \br\rangle$ of $\bR$ is not a free algebra. Then 
\begin{enumerate}
\item $\ass_R(S)=\ga$. 
\item $\bS :=\pi_\ga (S)\in \Den (\bR , 0)$ and $\RSm\simeq \bS^{-1}\bR \simeq \bR\bS^{-1}$.
\end{enumerate}
\end{theorem}

{\it Proof}. By the assumption,  for each  pair of elements $s\in S$ and $r\in R$, the subalgebra  $K\langle \bs, \br\rangle$ of $\bR$ is not a free algebra. Then  $\bS :=\pi_\ga (S)\in \Den (\bR , 0)$ (by Jategaonkar's Lemma, see \cite[Lemma 9.21]{Lam-ModRingbook}). Therefore,  $\ass_R(S)=\ga$ and $\RSm\simeq \bS^{-1}\bR \simeq \bR\bS^{-1}$. $\Box$ \\

{\bf The localizable sets of a ring and the posets $(\mrL (R), \subseteq )$ and $(\mrLoc (R),\ra )$.}  For an ideal $\ga $ of $R$, let $ \mrLRa :=\{ S\in \mrLR \, | \, \ass_R(S) = \ga \}$. Then 

%\marginpar{SLzRa}
\begin{equation}\label{SLzRa}
\mrLR= \coprod_{\ga \in \ass\,  \mrLR } \mrLRa
\end{equation}
is a disjoint union of non-empty sets. A set with binary operation, $(I, \leq )$,  is a {\em partially ordered set} (a poset, for short)  if $i\leq i$ for all $i\in I$, and if $i\leq j$ and $j \leq k$ for some elements $i,j,k\in I$ then $i\leq k$. The set $(\mrLR , \subseteq )$ is a poset (w.r.t. inclusion $\subseteq $), and $(\mrLRa, \subseteq )$ is a sub-poset of $(\mrLR , \subseteq )$ for every $\ga \in \ass \mrLR$. 

Let $\mrLoc (R)$ be the set of $R$-isomorphism classes of  all the localizations of the ring $R$, i.e. $R$-isomorphism classes in the set $\{ \RSm \, | \, S\in \mrLR \}$. If the rings $\RSm$ and $R\langle T^{-1}\rangle$ are $R$-isomorphic for some $S,T\in \mrLR$ then $\ass_R(S)=\ass_R(T)$. Therefore, 
%\marginpar{LocR}
\begin{equation}\label{LocR}
\mrLoc (R)= \coprod_{\ga \in \ass\,  \mrLR } \mrLoc (R,\ga )
\end{equation}
is a disjoint union of non-empty sets where $\mrLoc (R,\ga)$ is the  set of $R$-isomorphism classes in the set $\{ \RSm\, | \,  S\in \mrLRa\}$.

 The set
$(\mrLoc (R, \ga ),  \ra )$ is a poset where $A\ra B$ if there is a necessarily unique $R$-homomorphism from $A$ to $B$ which is denoted by 
$$\phi_{BA}:A\ra B.$$
In particular, the unique $R$-homomorphism $\phi _{TS}: R\langle S^{-1}\rangle \ra R\langle T^{-1}\rangle$ is also denoted by $\phi_{ R\langle S^{-1}\rangle  R\langle T^{-1}\rangle}$. 
 If  $A\ra B$ for some rings $A\in \mrLoc (R,\ga )$ and $B\in \mrLoc (R, \gb)$ then 
 %\marginpar{gasgb}
\begin{equation}\label{gasgb}
 \ga \subseteq \gb.
\end{equation}
  If, in addition,    $A =R\langle 
S^{-1}\rangle$ and  $B= R\langle 
T^{-1}\rangle$ for some localizable  sets $S,
T \in \mrLR $ then  if necessary we can assume that  $S\subseteq T$ (for example, by replacing  $T$ by $ \s_T^{-1}(B^\times)$ where $\s_T: R\ra B$,  see Proposition  \ref{SA20Aug21}.(1)).

 The map 
%\marginpar{Sm1RD}
\begin{equation}\label{Sm1RD}
 R\langle (\cdot )^{-1}\rangle: \, \mrLR\ra \mrLoc (R), \;\; S\mapsto \RSm,
\end{equation}
is an epimorphism of  posets. For each ideal $\ga \in  \ass \, \mrL (R)$, it
induces the epimorphism of   posets,  
 %\marginpar{Sm1RD1}
\begin{equation}\label{Sm1RD1}
 R\langle (\cdot )^{-1}\rangle: \, \mrL (R, \ga) \ra \mrLoc (R, \ga ), \;\; S\mapsto \RSm .
\end{equation}
For each ideal $\ga \in \ass \, \mrL (R)$, 
%\marginpar{SDen2}
\begin{equation}\label{SDen2}
\mrL (R, \ga )=\bigsqcup_{A \in \mrLoc (R, \ga )}\mrL (R, \ga ,
A )
\end{equation}
where $\mrL (R, \ga , A ):=\{ S\in \mrL (R, \ga )\, | \, \RSm \simeq A$,  an $R$-isomorphism (necessarily unique)$\}$.\\

{\bf The poset of localizable multiplicative sets $(\mrLR_m, \subseteq )$.} An Ore set or a denominator set of a ring is a multiplicative set. This fact suggests the following definitions. Let $\mrLR_m$ be the set of all localizable  multiplicative sets of $R$. By Lemma \ref{b17Dec22}, the localization of a ring at a localizable set is the same as  the localization of the ring at the localizable multiplicative set that the set generates.  For an ideal $\ga $ of $R$, let $ \mrLRa_m :=\{ S\in \mrLR_m \, | \, \ass_R(S) = \ga \}$. By Lemma \ref{b17Dec22}, $\ass\,  \mrLR =\ass\,  \mrLR_m:=\{ \ass_R(S)\, | \, S\in \mrLR_m\}$,  and so 
%\marginpar{MSLzRa}
\begin{equation}\label{MSLzRa}
\mrLR_m= \coprod_{\ga \in \ass\,  \mrLR } \mrLRa_m
\end{equation}
is a disjoint union of non-empty sets. The map 
%\marginpar{SSLmon}
\begin{equation}\label{SSLmon}
\mrLR\ra \mrLR_m, \;\; S\mapsto S_{mon}
\end{equation}
 is an epimorphism of posets with section 
 %\marginpar{SSLmon1}
\begin{equation}\label{SSLmon1}
\mrLR_m\ra \mrLR, \;\; S\mapsto S.
\end{equation}
 
 Similarly, for each $\ga \in \ass\, \mrLR$,  the map 
%\marginpar{aSSLmon}
\begin{equation}\label{aSSLmon}
\mrLRa\ra \mrLRa_m, \;\; S\mapsto S_{mon}
\end{equation}
 is an epimorphism of posets with section 
 %\marginpar{aSSLmon1}
\begin{equation}\label{aSSLmon1}
\mrLRa_m\ra \mrLRa, \;\; S\mapsto S.
\end{equation}
 
  The map 
%\marginpar{MSm1RD}
\begin{equation}\label{MSm1RD}
 R\langle (\cdot )^{-1}\rangle: \, \mrLR_m\ra \mrLoc (R), \;\; S\mapsto \RSm,
\end{equation}
is an epimorphism of  posets. For each ideal $\ga \in  \ass \, \mrL (R)$, it
induces the epimorphism of   posets,  
 %\marginpar{MSm1RD1}
\begin{equation}\label{MSm1RD1}
 R\langle (\cdot )^{-1}\rangle: \, \mrL (R, \ga)_m \ra \mrLoc (R, \ga ), \;\; S\mapsto \RSm .
\end{equation}
For each ideal $\ga \in \ass \, \mrL (R)$, 
%\marginpar{MSDen2}
\begin{equation}\label{MSDen2}
\mrL (R, \ga )_m=\bigsqcup_{A \in \mrLoc (R, \ga )}\mrL (R, \ga ,
A )_m
\end{equation}
where $\mrL (R, \ga , A )_m:=\{ S\in \mrL (R, \ga )_m\, | \, \RSm \simeq A$,  an $R$-isomorphism (necessarily unique)$\}$.\\

{\bf Direct limits of rings.} Let $(I, \leq )$ be a poset, $\{ R_i\}_{i\in I}$ be rings,  and for each pair of elements $i$ and $j$ of $I$ such that $i\leq j$ there is a ring homomorphism $ f_{ji}: R_i\ra R_j$ such that $f_{ii}$ is the identity map for all $i\in I$, and $$f_{ki}=f_{kj}f_{ji}$$ for all elements $i,j,k\in I$ such that  $i\leq j \leq k$. The system  $(R_i,f_{ij})$ is called  a {\em direct system} of rings on $I$. A ring 
$\varinjlim_{i\in I}  R_i$ together with ring homomorphisms $g_i: R_i\ra \varinjlim_{i\in I}  R_i$ 
such that $g_i=g_jf_{ji}$  for all $i\leq j$
is called the {\em direct limit} of the direct system $(R_i, f_{ji})$ if it satisfies  the following {\em universal property}: For each system of ring homomorphisms $g_i': R_i\ra A$ such that $g_i'=g_j'f_{ji}$  for all $i\leq j$, there is a unique homomorphism  
$$\alpha :\varinjlim_{i\in I}  R_i\ra A$$ such that $g_i'=\alpha g_i$ for all $i\in I$. It follows from the universal property that the direct limit is unique (up to isomorphism).

 Suppose that  the poset $(I,\leq )$ is a {\em directed set}, that is for each pair of  elements $i,j\in I$ there is an element $k\in I$ such that $i\leq k$ and $j\leq k$.   On the disjoint union $\coprod_{i\in I}R_i$, we say that elements $a_i\in R_i$ and $a_j\in R_j$ are {\em equivalent}, $a_i\sim a_k$,   if $$f_{ki}(a_i)=f_{kj}(a_j)$$
 for some  element $k\in I$ such that $i,j\leq k$. Then 
 the set of equivalence classes,  $\coprod_{i\in I}R_i/\sim, $
 is isomorphic to $\varinjlim_{i\in I}  R_i$ where for each $i\in I$, the map $$g_i: R_i\ra \varinjlim_{i\in I}  R_i, \;\; r\mapsto  \widetilde{r}$$ is a ring homomorphism and  $\widetilde{r}$ is the equivalence class of the element $r$, such that for all $i\leq j$, $g_i=g_jf_{ji}$.

Let $(I, \leq )$ be a poset,   $(R_i,f_{ij})$ be  a  direct system of rings on $I$. Let $\CI $ be an additive subgroup of the the direct sum of $K$-modules  $\bigoplus_{i\in I}R_i$ generated by all the elements 
%\marginpar{thopi1}
\begin{equation}\label{thopi1}
\{ a_i-a_j\, | \, f_{ki}(a_i)=f_{kj}(a_j),  a_i\in R_i, a_j\in R_j, i\leq k , j\leq k\}.
\end{equation} 
By the definition, $\CI$ is a $K$-submodule of $\bigoplus_{i\in I}R_i$.  Let $*_{i\in I}R_i$ be the free product of  $K$-algebras  $R_i$. For each $i\in I$, there is a natural ring homomorphism $\theta_i: R_i\ra *_{i\in I}R_i$, and the ring $ *_{i\in I}R_i$ is generated by the images of the rings $R_i$. Let
%\marginpar{thopi}
\begin{equation}\label{thopi}
\th : \bigoplus_{i\in I}R_i\ra  *_{i\in I}R_i, \;\; r_i\mapsto \th_i(r_i)\;\; {\rm for} \;\; r_i\in R_i.
\end{equation} 
An element $i\in I$ is the largest element of the poset $I$ if for all elements $j\in I$, $j\leq i$.
 
 \begin{theorem}\label{A14Dec22}%\marginpar{A14Dec22}
Let $(I, \leq )$ be a poset and    $(R_i,f_{ij})$ be  a  direct system of rings on $I$.  Then
\begin{enumerate}
\item $\varinjlim_{i\in I}  R_i\simeq *_{i\in I}R_i/ (\th (\CI ))$  where $(\th (\CI))$ is the ideal of the free product of $K$-algebras $*_{i\in I}R_i$ generated by the set $\th (\CI)$. 
\item Suppose that for each element $i\in I$ there is an element $j\in \max (I)$ such that $i\leq j$.  Then $$\varinjlim_{i\in I}  R_i\simeq *_{\mu\in \max (I)}R_\mu/ \gf$$  where $*_{\mu\in \max (I)}R_\mu$ is a free product of the set of $K$-algebras $\{ R_\mu \, | \, \mu \in \max (I)\}$ and $ \gf$ is an ideal of $*_{\mu\in \max (I)}R_\mu$ which is generated by  the set $\th  \bigg( \bigg(\bigoplus_{\mu \in \max (I)}R_\mu\bigg)\bigcap \CI \bigg)$. 
\item Suppose that an element $i\in I$ is the largest element of the poset $I$. Then $\varinjlim_{j\in I}  R_j\simeq R_i$. 
\end{enumerate}
\end{theorem}

{\it Proof}. 1.  To prove statement 1 it suffices to show that  the ring  $*_{i\in I}R_i/ (\th (\CI))$ satisfies the universal property of the ring $\varinjlim_{i\in I}  R_i$.  Let $g_i': R_i\ra A$, $i\in I$, be  a set of ring homomorphisms such that $g_i'=g_j'f_{ji}$ for all $i\leq j$.  Then there is a unique ring homomorphism $\rho : *_{i\in I}R_i \ra A$ such that $g_i'=\rho \th_i$ for all $i\in I$. Then 
 $$\rho \th (\CI )=0$$ since for all elements $a_i\in R_i$ and $a_j\in R_j$ such that $f_{ki}(a_i)=f_{kj}(a_j)$ for some $k\in I$ such that $i,j\leq k$, 
 $$\rho (\th_i(a_i)-\th_j(a_j))=g_i'(a_i)-g_j'(a_j)=g_k'f_{ki}(a_i)-g_k'f_{kj}(a_j)=g_k'(f_{ki}(a_i)-f_{kj}(a_j))=g_k'(0)=0.$$
 
2.  Let $\varinjlim^{K}_{i\in I} R_i $ be the direct limit of the  $K$-modules $R_i$. 

(i) $\varinjlim^{K}_{i\in I} R_i =\bigg( \bigoplus_{i\in I}R_i\bigg)/\CI$: By the definition, the $K$-module $ \bigg( \bigoplus_{i\in I}R_i\bigg)/\CI$ satisfies the universal property of the direct limit of $K$-modules $R_i$  and the statement (i) follows from  the uniqueness of the direct limit of $K$-modules.

(ii) $\bigg( \bigoplus_{i\in I}R_i\bigg)/\CI\simeq \bigg( \bigoplus_{\mu \in \max (I)}R_\mu\bigg)/\bigg(\bigoplus_{\mu \in \max (I)}R_\mu\bigg)\bigcap \CI$: 
 By the assumption, for each element $i\in I$ there is an element $j\in \max (I)$ such that $i\leq j$. Then it follows from the  inclusion $\bigoplus_{\mu \in \max (I)}R_\mu\subseteq \bigoplus_{i\in I}R_i$ that  
 $$ \bigg( \bigoplus_{i\in I}R_i\bigg)/\CI= \bigg( \bigoplus_{\mu\in \max (I)}R_i+\CI\bigg)/\CI\simeq \bigg( \bigoplus_{\mu \in \max (I)}R_\mu\bigg)/\bigg(\bigoplus_{\mu \in \max (I)}R_\mu\bigg)\bigcap \CI.$$

(iii) $\varinjlim{}^{K}_{i\in I} R_i =\bigg( \bigoplus_{i\in I}R_i\bigg)/\CI  \simeq \bigg( \bigoplus_{\mu \in \max (I)}R_\mu\bigg)/\bigg(\bigoplus_{\mu \in \max (I)}R_\mu\bigg)\bigcap \CI$: The statement (iii) follows from   the statements (i) and (ii).

(iv) {\em The ring  $*_{\mu\in \max (I)}R_\mu/ \gf$ satisfies the universal property of the ring $\varinjlim_{i\in I} R_i $}: Let $g_i': R_i\ra A$, $i\in I$, be  a set of ring homomorphisms such that $g_i'=g_j'f_{ji}$ for all $i\leq j$. By statement 1 and the statement (iii),  the homomorphisms $g_i'$ induce a unique ring homomorphism, $\xi :  *_{\mu\in \max (I)}R_\mu/ \gf\ra A$ such that $g_\mu'=\xi\th_\mu$ for all $\mu \in \max (I)$. 
 Now,   statement 2 follows from the statement (iv). 
 
3. Statement 3 follows from statement 2. $\Box$

Let $T\in \mrL (R)$. Then  the power set of $T$, i.e. the set of all subsets of  the set $T$,  $(\CP (T), \subseteq )$, is a poset. The system  $\{ R\langle S^{-1}\rangle, \phi_{SS'}\}$ is a direct system where $S, S'\in \CP (T)$. Corollary \ref{a7Dec22}.(1) describe the direct limit $\varinjlim_{\{ S\subseteq T\}} R\langle S^{-1}\rangle$.

 \begin{corollary}\label{a7Dec22}%\marginpar{a7Dec22}
Let $R$ be a ring and $T\in \mrLR$. Then 
\begin{enumerate}
\item $R\langle T^{-1}\rangle\simeq \varinjlim_{\{ S\subseteq T\}} R\langle S^{-1}\rangle$.
\item $\ass_R(T)=\bigcup_{S\subseteq T}\ass_R(S)$. 
\end{enumerate}
\end{corollary}

{\it Proof}. 1. Statement 1 follows from Proposition \ref{A7Dec22}.(2).

2.  Statement 2 follows from statement 1.  $\Box $

We collect two results on direct limits, Lemma \ref{a21Dec22} and Lemma \ref{a27Dec22},  that are used later in the paper and are interesting on their own.

\begin{lemma}\label{a21Dec22}%\marginpar{a21Dec22}
Let $(I, \leq )$ be a directed set  and    $(R_i,f_{ij})$ be  a  direct system of rings on $I$.  Suppose that all rings $R_i$, where $i\in I$,  are isomorphic to a ring $R$ and  all the maps $f_{ij}$ are isomorphisms. Then $\varinjlim_{i\in I}  R_i\simeq R$.
\end{lemma}

{\it Proof}. Fix an element $i\in I$ and an isomorphism $g_i:R_i\ra R$. For each element $j\in I$ such that $i\leq j$, let $g_j:= g_if_{ij}: R_j\ra R$ where $f_{ij}:=f_{ji}^{-1}$. For each element $k\in I$ such that $i\not\leq k$, fix an element $j\in I$ such that $i\leq j$ and $k\leq j$. Let $$g_k:= g_if_{ij}f_{jk}: R_k\ra R\;\; {\rm  where}\;\; f_{jk}:=f_{kj}.$$ Then it is easy to verify that the map $g_k$ does not depend on the choice of the element $j$ and $(R, g_i)_{i\in I}$ satisfies the universal property of the direct limit for the direct system $(R_i, f_{ij})$, and so $R=\varinjlim_{i\in I}R_i$. 
 $\Box $
 
 Lemma \ref{a27Dec22} gives an example of a direct system of rings with zero direct limit.

\begin{lemma}\label{a27Dec22}%\marginpar{a27Dec22}
Let $J=\{1\}\coprod I$ be a disjoint union of two sets where the set $I$  consists at least of two elements. The set $J$ is a poset where $1<i$ for all $i\in I$. Let $\{ R_i\, | \, i\in I\}$ be a set of rings and  $A_1=\prod_{i\in I}R_i$. For each $i\in I$, let the map $f_{1i}:R_1\ra R_i$ be  the  projection onto $R_i$, $f_{11}=\id _{R_1}$ and $f_{ii}=\id_{R_i}$ (the identity maps). Then $\varinjlim_{j\in J}  R_i\simeq \{0\}$.
\end{lemma}

{\it Proof}. For each $i\in I$, let $1_i$ be the identity element of the ring $R_i$, and $1_i'$ be an element of $R_1$ where on the $i$'th place is $1_i$ and zeros are elsewhere. Clearly, $\max (J)=I$. Now, by Theorem \ref{A14Dec22}.(2), $\varinjlim_{j\in J}  R_j\simeq *_{i\in I}R_i/\gf$. Since for all $i\in I$,
 $$1_i=1_i-0=(1_i-1_i')+ (1_i'-p_j(1_i'))\in \gf$$ where $j\neq i$, $\varinjlim_{j\in J}  R_i\simeq \{0\}$ (see Theorem \ref{A14Dec22}.(2) for the definition of the ideal $\gf$). $\Box$\\

{\bf The absolute quotient ring $Q_a(R)$ of a ring $R$.} The set   $(\mrLR, \subseteq )$   is a poset and the system $\{ R\langle S^{-1}\rangle, \phi_{ST}\}$ is a direct system of rings where $S, T\in \mrLR$.\\

{\it Definition.} The direct limit $Q_a(R):=\varinjlim_{S\in \mrLR}\RSm$ is called the {\em absolute quotient ring}  of the ring $R$. Let $\s_a: R\ra Q_a(R)$ be a unique $R$-homomorphism, $\ga_a(R):=\ker (\s_a)$, $R_a:=R/\ga_a(R)$,  $\pi_a:R\ra R_a$, $r\mapsto r+\ga_a$, and $\overline{\s_a}:R_a\ra Q_a (R)$ is a unique $R$-homomorphism. Clearly, $\s_a=\overline{\s_a}\pi_a$.

\begin{theorem}\label{XA30Dec22}%\marginpar{XA30Dec22}
$Q_a(R)=\begin{cases}
\RSm& \text{if }\max\mrLR =\{ S\},\\
\{ 0\}& \text{otherwise}, \\
\end{cases}$ and $\ga_a(R)=\begin{cases}
\ass_R(S)& \text{if }\max\mrLR =\{ S\},\\
R& \text{otherwise}.\\
\end{cases}$
\end{theorem}

{\it Proof}. Suppose that $\max\mrLR =\{ S\}$. Then, by  Theorem \ref{SC12Jan19}.(2), the element $S$ is the largest element of the poset $\mrLR$, and the statement  follows from Theorem \ref{A14Dec22}.(3).

Suppose that the set $\max\mrLR$ contains at least two  elements, say $S$ and $T$. The set $\CS :=\s_a^{-1}(Q_a(R)^\times)$ is a localizable set in $R$ since $\s_a(\CS)\subseteq Q_a(R)^\times$. 
The diagram below is commutative 
$$\begin{array}[c]{ccc}
R&\stackrel{\s_S}{\ra} &\RSm\\
&\stackrel{\s_a}{\searrow}&{\downarrow}^{\phi_S}\\
&  & Q_a(R)
\end{array}$$
where $\phi_S$ is a unique $R$-homomorphism. Then it follows from the inclusion $\phi_S (\RSm^\times)\subseteq  Q_a(R)^\times$ that  $$\CS =\s_a^{-1}(Q_a(R)^\times)=(\phi_S\s_S)^{-1}(Q_a(R)^\times)=\s_S^{-1}\phi_S^{-1}(Q_a(R)^\times)\supseteq \s_S^{-1}(\RSm^\times)=S,$$
by Corollary \ref{aSA19Aug21}.(1). Similarly, $\CS\supseteq T$. Therefore, $\CS\supseteq S\cup T$, and so the set $S\cup T$ is a localizable set of $R$ as a subset of the localizable set $\CS$, a contradiction (since $S,T\in \max \mrLR$ and $S\neq T$). The proof of the first equality of the theorem is complete. The second equality of the theorem follows from the first.  $\Box$

{\it Example.} If $R$ is a commutative domain then the set $\CC_R=R\backslash \{ 0\}$ is the largest element of $\mrLR$. By Theorem  \ref{XA30Dec22},  $Q_a(R)=Q_{cl}(R)$ is the classical quotient ring of $R$, i.e. the field of fractions of $R$, and $\ga_a(R)=0$. 

\begin{lemma}\label{xa30Dec22}%\marginpar{xa30Dec22}
If  $R$ is  a  simple Artinian rings then $\max\mrL (R)=\{ R^\times\}$,  $Q_a(R)=R$ and $\ga_a(R)=0$.
\end{lemma}

{\it Proof}.  Every non-unit of the simple Artinian $R$ ring is a zero divisor. Since the ring $R$ is simple every non-unit $r\in R$  is a non-localizable element since $\ass_R(r)=R$. Therefore, $\CL \mrLR = R^\times$ and $\CN\CL \mrLR = R\backslash R^\times$. Hence, $\max \mrLR = \{ R^\times\}$. Now, by Theorem \ref{XA30Dec22}, $Q_a(R)=R\langle (R^\times)^{-1}\rangle=R$ and $\ga_a(R)=0$. $\Box$

We need the lemma below in order to prove Proposition \ref{a30Dec22} which gives a practical sufficient condition for $Q_a(R)=\{0\}$ and $|\max \mrLR|\geq 2$.

\begin{lemma}\label{A30Dec22}%\marginpar{A30Dec22}

\begin{enumerate}
\item $Q_a(R)\simeq \Big(*_{S\in \max \mrLR}\RSm \Big)/\gf_a$ where $\gf_a$ is an ideal of the $K$-algebra $*_{S\in \max \mrLR}\RSm$ which is generated by the set $\th \Big( \bigoplus_{S\in \max \mrLR}\RSm \cap \CI\Big)$, see Theorem  \ref{A14Dec22}.(2).
\item $\ga_a(R)\supseteq \sum_{S\in \max \mrLR} \ass_R(S)$. 
\item $Q_a(R)=\{0 \}$ iff  $\ga_a(R)=R$. If $\sum_{S\in \max \mrLR} \ass_R(S)=R$ then  $Q_a(R)=\{0 \}$. 
\end{enumerate}
\end{lemma}

{\it Proof}. 1. Statement 1 is a particular case of  Theorem  \ref{A14Dec22}.(2).

2. Let $S\in \max \mrLR$. Then for each element $r\in \ass_R(S)$ and each $T\in \max \mrLR$, 
$$\s_T(r)=(\s_T(r)-r)-(0-r)=(\s_T(r)-r)-(\s_S(r)-r)\in \th \bigg( \bigoplus_{S\in \max \mrLR}\RSm \cap \CI\bigg),$$ and so $\ass_R(S)\subseteq \ga_a(R)$, and statement 2 follows.

3.  $Q_a(R)=\{0 \}$ iff $1=0$ in  $Q_a(R)$ iff  $\ga_a(R)=R$.  If $\sum_{S\in \max \mrLR} \ass_R(S)=R$ then, by statement 2,  $\ga_a(R)=R$, and so $Q_a(R)=\{0 \}$.  $\Box $

\begin{proposition}\label{a30Dec22}%\marginpar{a30Dec22}
 Given localizable sets $\{S_i \, | \, i\in I\}$ of a ring $R$ such that $\sum_{i\in I}\ass_R(S_i)=R$. Then $Q_a(R)=\{0\}$ and $|\max \mrLR|\geq 2$.
\end{proposition}

{\it Proof}.  By Theorem \ref{SC12Jan19}.(2), every localizable set $S_i$ is contained in a maximal localizable set $S_i'$. Since $\ass_R(S_I)\subseteq \ass_R(S_I')$ and $\sum_{i\in I}\ass_R(S_i)=R$, we have that $\sum_{i\in I}\ass_R(S_i')=R$, and so  $\sum_{S\in \max \mrLR} \ass_R(S)=R$. By Lemma \ref{A30Dec22}.(2), $\ga_a(R)=R$, and so $Q_a(R)=0$, by Lemma \ref{A30Dec22}.(3). By Theorem  \ref{XA30Dec22}, $|\max \mrLR|\geq 2$. $\Box $

\begin{corollary}\label{ya30Dec22}%\marginpar{ya30Dec22}
 Let $A=\prod_{i\in I}A_i$ be a product of  rings where $I$ is an arbitrary set that contains at least two elements. Then 
$Q_a(A)=0$ and $|\max \mrLR|\geq 2$. 
\end{corollary}

{\it Proof}.  By the assumption,  the set $I$ contains at least two elements. 
 For each $i\in I$, the set $S_i=A_i^\times \times \prod_{j\neq i}A_j$ is a  localizable set of $A$ with $\ass_R(S_i)=\{0\}\times \prod_{j\neq i}A_j$ since $A\langle S_i^{-1}\rangle\simeq A_i$. Let $i$ and $j$ be  distinct elements of the set $I$. Since $A=\ass_A(S_i)+\ass_A (S_j)$, we   have that $Q_a(A)=\{ 0\}$ and $|\max \mrLR|\geq 2$, by Proposition \ref{a30Dec22}.  $\Box $\\

{\bf Every Ore set is a localizable set.} Let $S$ be an Ore set of the ring $R$. Theorem \ref{S10Jan19}  states  that  every Ore set is localizable,  gives an explicit description of this ideal and the ring $R\langle S^{-1}\rangle$, and $\ga(S)=\ass_R (S)$. So, the localization an any Ore set is an example of the perfect localization. Furthermore,   Theorem \ref{S10Jan19} also  states that  the ring  $R\langle S^{-1}\rangle$ is $R$-isomorphic to the localization $\bS^{-1}\bR$ of the ring $\bR$ at the denominator set $\bS$ of $\bR$.

\begin{theorem}\label{S10Jan19}%\marginpar{S10Jan19}
Let $R$ be a ring and $S\in \Ore (R)$.
\begin{enumerate}
\item \cite[Theorem 4.15]{Bav-intdifline} Every Ore set is a localizable set. 
\item \cite[Theorem 1.6.(1)]{LocSets} $\ga :=\{ r\in R\, | \, srt=0$ for some elements $s,t\in S\}$ is an ideal of $R$ such that $\ga \neq R$. 
\item \cite[Theorem 1.6.(2)]{LocSets} Let $\pi : R\ra \bR :=R/\ga$, $r\mapsto \br =r+\ga$. Then $\bS :=
\pi (S) \in \Den (\bR , 0)$, $\ga = \ga (S)=\ass_R(S)$, $S\in \mrL (R, \ga )$,  and $S^{-1}R\simeq \bS^{-1}\bR$, an $R$-isomorphism. In particular, every Ore set is localizable. 
%\item \cite[Theorem 1.6.(3)]{LocSets} Let $\gb$ be an ideal of $R$ and $\pi^\dag :R\ra R^\dag :=R/\gb$, $r\mapsto r^\dag =r+\gb$. If $S^\dag := \pi^\dag (S)\in \Den (R^\dag ,0)$ then $\ga \subseteq \gb $ and the map $$\bS^{-1}\bR \ra {S^\dag}^{-1}R^\dag , \;\; \bs^{-1}\br \mapsto {s^\dag}^{-1}r^\dag$$ is a ring epimorphism. 
%\item \cite[Theorem 1.6.(4)]{LocSets} Let $f: R\ra Q$ be a ring homomorphism such that $f(S)\subseteq Q^\times$ and the ring $Q$ is generated by $f(R)$ and $\{ f(s)^{-1}\, | \, s\in S\}$. Then  
%\begin{enumerate}
%\item $\ga \subseteq \ker (f)$ and the map $$\bS^{-1}\bR \ra Q, \;\;  \bs^{-1}\br \mapsto f(s)^{-1}f(r)$$ is a ring epimorphism with kernel $\bS^{-1}(\ker (f) / \ga )$. 
%\item Let $\widetilde{R} = R/\ker (f)$ and $\tpi : R\ra \widetilde{R}$, $r\mapsto r+\ker (f)$. Then $\widetilde{S}:= \tpi (S) \in  \Den (\widetilde{R}, 0)$ and $\widetilde{S}^{-1}\widetilde{R}\simeq Q$, an $\widetilde{R}$-isomorphism. 
%\end{enumerate}
\end{enumerate} 
\end{theorem}

{\bf Criterion for a left Ore set to be a left localizable set.} 
 For a ring $R$ and  its ideal $\ga$, let 
\begin{eqnarray*}
 {}'\Den_l(R , \ga)&:=&\{ S\in \Den_l(R)\, | \, \ass_l(S) = \ga, S\subseteq \pCCR \}, \\
\Den_r'(R, \ga )&:=&\{ S\in \Den_r(R)\, | \, \ass_r(S) = \ga, S\subseteq \CC_R' \}.
\end{eqnarray*}
Theorem \ref {B16Jan19}.(1)  is  a criterion  for a left Ore set to be a left localizable set and  Theorem \ref {B16Jan19}.(2) describes the structure of the localization of a ring at a localizable left Ore set.

\begin{theorem}\label{B16Jan19}%\marginpar{B16Jan19}
(\cite[Theorem 1.5]{LocSets}) Let $R$ be a ring, $S\in \Ore_l(R)$, and  $\ga = \ass_R(S)$. Then
\begin{enumerate}
\item $S\in \mrLR$ iff ${}'\ga \neq R$ where the ideal ${}'\ga = {}'\ga (S)$ of $R$ is as in Proposition \ref{A12Jan19}.(2) and (\ref{ala1}). 
\item  Suppose that ${}'\ga \neq R$. Let ${}'\pi : R\ra {}'
R:=R/ {}'\ga$, $r\mapsto {}'r= r+{}'\ga$ and ${}'S={}'\pi (S)$. Then 
\begin{enumerate}
\item $'S\in {}'\Den_l({}'R)$.
\item $\ga = {}'\pi^{-1}(\ass_l({}'S))$.
\item $S^{-1}R\simeq {}'S^{-1}{}'R$, an $R$-isomorphism. 
\end{enumerate}
\end{enumerate}
\end{theorem}
\cite[Theorem 2.11]{LocSets} is a criterion for a right Ore set to be a localizable set.

\begin{definition} (\cite{LocSets}) A multiplicative set $S$ of a ring $R$ is called a {\em left localizable set} of $R$ if  
$$\RSm = \{ \bs^{-1} \br \, | \, \bs \in \bS, \br \in \bR\}\neq \{ 0\}$$
where $\bR = R/ \ga$,  $\ga = \ass_R(S)$ and $\bS = (S+\ga ) / \ga$, i.e., every element of the ring $\RSm$ is a left fraction $\bs^{-1} \br$ for some elements  $\bs \in \bS$ and $ \br \in \bR$. Similarly,  a multiplicative set $S$ of a ring $R$ is called a {\em right  localizable set}  of $R$  if 
$$\RSm = \{  \br \bs^{-1}\, | \, \bs \in \bS, \br \in \bR\}\neq \{ 0\}, $$
 i.e., every element of the ring $\RSm$ is a right  fraction $ \br\bs^{-1}$ for some elements  $\bs \in \bS$ and $ \br \in \bR$. A right and left localizable set of $R$ is called a {\em localizable set} of $R$.
\end{definition}

%%%%%%%%%%%%%%%%%%%%%  Section 3   %%%%%%%%%%%%%%%%%%%%%

\section{Maximal localizable sets, the localization radical, the sets of localizable and non-localizable elements}\label{MAXLOCSETS}%\marginpar{MAXLOCSETS}

 For a ring $R$, the following concepts are introduced: the maximal localizable sets, the localization  radical $\mrL\rad (R)$, the set of localizable elements $\CL\mrL (R)$, the set of non-localizable elements $\CN\CL\mrL (R)$, the set of completely  localizable elements $\CC\mrL (R)$. All these objects  are characteristic subsets of $R$, i.e. they are  invariant under the action of the automorphism group of $R$. 
 Proposition \ref{A29Dec22} shows that there are tight connections between the sets $\CL \mrLR$, $\CN \CL\mrLR$,  and $ \mrL\rad (R)$. For an arbitrary ring, 
Theorem \ref{SC12Jan19} states that  the set of maximal localizable set is a non-empty set and every localizable set  is a subset of a maximal localizable set. 
 
Proposition \ref{SA19Aug21}.(1), is an explicit description of    the largest element  $ \CS (R, \ga , A)$ of the poset $(\mrL(R, \ga,  A), \subseteq )$. Corollary \ref{a20Dec22} is a criterion for a set  $S\in \mrL (R, \ga, A )$ to be the largest element of the poset $\mrL (R, \ga, A )$.  
Theorem \ref{22Dec22} describes the sets  $\max \mrL (R)$ and $\max \mrLoc (R)$, presents a bijection between them, and shows that $\max \mrLoc (R)\neq \emptyset$.  Similarly, for each ideal $\ga\in \ass \mrLR$, Theorem \ref{A23Dec22} describes the sets  $\max \mrLRa$ and $\max \mrLoc (R, \ga)$, presents a bijection between them, and shows that they are nonempty sets.

 The $\ga$-absolute quotient ring $Q_a(R, \ga )$ is introduced for each ideal $\ga \in \ass \, \mrLR$. Theorem \ref{aXA30Dec22} gives an explicit description of the ring $Q_a(R,\ga)$, a criterion for $Q_a(R,\ga)\neq \{0\}$, and a criterion for $\ga_{a,\ga}=\ga$.
  
  Theorem \ref{20Dec22}.(1) is an $R$-isomorphism  criterion for  localizations of the ring $R$ and Theorem \ref{20Dec22}.(2) is a criterion for existing of an $R$-homomorphism between localizations of $R$. Proposition \ref{B29Dec22}  contains some properties of $\CC\mrLR$ and $\gc_R$. 
Theorem \ref{A17Jan23} is a criterion for the ring $\RSm$ to be an $R$-isomorphic to a localization of $R$ at a denominator set and  it also  describes all such denominator sets.

For a denominators set $T\in \Den_*(R, \ga)$, Proposition \ref{B17Jan23}  describes all localizable sets $S$ of the ring $R$ such that the ring $\RSm$ is an $R$-isomorphic to the localization of $R$ at $T$. For a localizable  multiplicative   set $S$ of the ring $R$, Proposition \ref{C17Jan23} is a criterion for the set $S$ to be a left denominator set of $R$.  \\
 
 {\bf The set $\max  \mrLR $ of maximal elements of $\mrLR$ is a non-empty set.}   For a ring $R$, the set $\maxDen_l(R)$ of maximal left denominator sets  (w.r.t. $\subseteq$) is a {\em non-empty}
 set, \cite[Lemma 3.7.(2)]{larglquot} and the sets of maximal localizable left or right or two-sided Ore sets are also {\em non-empty} sets \cite[Theorem 1.9]{LocSets}. Let $\max \, \mrLR$ be the set of maximal (w.r.t. $\subseteq $)  elements of the set $\mrLR$. An element of the set  $\max \, \mrLR$ is called a {\em maximal localizable set} of the ring $R$.

 \begin{theorem}\label{SC12Jan19}%\marginpar{SC12Jan19}
Let $R$ be a ring. Then 
\begin{enumerate}
\item $\max  \mrLR\neq \emptyset$.
\item Each set $S\in  \mrLR$ is contained in a set $\CS\in  \max  \mrLR$.
\end{enumerate}
\end{theorem}

{\it Proof}. 1. The theorem follows from Zorn's Lemma. 
 Let $\{S_i\}_{i\in \G}\subseteq \mrL (R)$  be an ordered by inclusion set of localizable sets in $R$ indexed by an ordinal $\G$, i.e. for each pair $i\leq j$ in $\G$, $S_i\subseteq S_j$ and if $i \in \G$ is a limit ordinal then $S_i=\bigcup_{j<i}S_j$. In view of Zorn's Lemma, it suffices  to show that  $\bigcup_{i\in \G}S_i\in \mrL (R)$. By (\ref{RSmRT}) and  (\ref{RSmRT1}), we have the direct system of ring homomorphisms, $(R\langle S_i^{-1}\rangle, \phi_{S_jS_i})$, where for each pair $i\leq j$, the map $\phi_{S_jS_i}: R\langle S_i^{-1}\rangle\ra R\langle S_j^{-1}\rangle$ is defined in (\ref{RSmRT}). Since the rings $R\langle S_i^{-1}\rangle$ are unital and $\phi_{S_jS_i}(1)=1$ for all $i\leq j$, the direct limit $\varinjlim R\langle S_i^{-1}\rangle$  is a nonzero ring and  and is isomorphic to $R\langle \Big(\bigcup_{i\in \G}S_i\Big)^{-1}\rangle$. Therefore, $\bigcup_{i\in \G}S_i\in \mrL (R)$. 
 
 2. Repeat the proof of statement 1 for the set of all localizable sets that contain the set $S$. $\Box$

  \begin{corollary}\label{b4Dec22}%\marginpar{b4Dec22}
Let $R$ be a ring and $T\in \max \, \mrLR$. Then 
\begin{enumerate}
\item $R\langle T^{-1}\rangle\simeq \varinjlim_{\{ S\subseteq T\}} R\langle S^{-1}\rangle$.
\item $\ass_R(T)=\bigcup_{S\subseteq T}\ass_R(S)$. 
\end{enumerate}
\end{corollary}

{\it Proof}. The corollary is a particular case of Corollary \ref{a7Dec22}. $\Box $

Let $\max \ass \, \mrLR$ be the set of maximal elements of the set $\ass \, \mrLR = \{ \ass_R(S)\, | \, S\in \mrLR\}$ and $\ass \max \mrLR := \{ \ass_R(\CS)\, | \, \CS \in \max \mrLR\} $. 

 \begin{corollary}\label{b22Dec22}%\marginpar{b22Dec22}
$\max \ass \, \mrLR\subseteq \ass \max \mrLR$. 
\end{corollary}

{\it Proof}. The corollary follows from Theorem  \ref{SC12Jan19}. $\Box $\\

A set of subsets is called an {\em antichain} if the sets are {\em incomparable}, that is none  of the sets is a subset of another one. 

%*******************************

%Give the results for Den sets and Ore sets

%***************************

{\em Question. Is $\max \ass \, \mrLR= \ass \max \mrLR$? I.e. the set $\ass \max \mrLR$ is an antichain.} \\

%**********************************************

{\bf The largest element $\CS(R, \ga , A)$ in $(\mrL(R, \ga , A) \subseteq )$ and its characterizations where $A\in \mrLoc (R, \ga )$.}

\begin{lemma}\label{a22Dec22}%\marginpar{a22Dec22}
Given  sets $S, \CS\in \mrL (R, \ga , A)$ such that   $S\subseteq \CS$. Then $T\in   \mrL (R, \ga , A)$ for all subsets $T$ of $R$ such  $S\subseteq T\subseteq \CS$. 
\end{lemma}

{\it Proof}. The inclusions $S\subseteq T\subseteq \CS$, imply that $\ga =\ass_R(S)\subseteq \ass_R(T)\subseteq \ass_R(\CS)=\ga$, and so $\ass_R(T)=\ga$. 

It suffices to show that the rings $R\langle T^{-1}\rangle$ and $A$ are $R$-isomorphic.   Since $S\subseteq \CS$ and  $S, \CS\in \mrL (R, \ga , A)$, the unique $R$-homomorphism $A\simeq \RSm \ra R\langle \CS^{-1}\rangle\simeq A$ is the identity map $\id_A$. In particular, $\s_S (T)\subseteq \s_S (\CS)\subseteq A^\times$ since 
$T\subseteq \CS$ where $\s_S:R\ra \RSm\simeq A$. Now, 
$$R\langle T^{-1}\rangle \simeq R\langle S^{-1}\rangle\langle \s_S(T)^{-1}\rangle \simeq A\langle \s_S(T)^{-1}\rangle=A$$ since $\s_S(T)\subseteq A^\times $, as required.  $\Box $

Proposition \ref{SA19Aug21}.(1), is an explicit description of    the largest element  $ \CS (R, \ga , A)$ of the poset $(\mrL(R, \ga,  A), \subseteq )$.
 
 \begin{proposition}\label{SA19Aug21}%\marginpar{SA19Aug21}
Let  $A\in \mrLoc (R, \ga)$, $A^\times$ be the group of units of the ring $A$,  
and  $\s  : R\ra A$, $ r\mapsto r+\ga$.  Then 
\begin{enumerate}
\item The set $$ \CS (R, \ga , A):=\s^{-1}(A^\times )=\bigcup_{S\in \mrL (R, \ga , A)}S$$ is the largest element of the poset $(\mrL(R, \ga,  A), \subseteq )$. In particular,  the poset $(\mrL(R, \ga,  A), \subseteq )$ is a directed set such that if $S\in \mrL(R, \ga,  A)$ then $T\in  \mrL(R, \ga,  A)$ for all subsets $T$ of $R$ such that $S\subseteq T\subseteq \CS (R, \ga,  A)$. 
\item The set $ \CS (R, \ga , A)$ is a multiplicative set in  $R$ such that $ \CS (R, \ga , A)+\ga = \CS (R, \ga , A)$. 
\item $A\simeq R\langle \CS (R, \ga , A)^{-1}\rangle\simeq  \varinjlim_{S\in \mrL (R, \ga , A)}\RSm$, $R$-isomorphisms. 
\end{enumerate}
\end{proposition} 

{\it Proof}. 1.  Let $\CS:=\s^{-1}(A^\times )$. 

(i) $\CS\in \mrL (R, \ga,  A)$:  Fix an element $S\in \mrL (R, \ga , A)$. Then the statement (i) follows from Lemma \ref{c17Dec22}.(3). 

%Since $\s (\CS)\subseteq A^\times$, $\CS\in \mrL (R)$, by Proposition \ref{a4Dec22}. Clearly, $\ass_R(\CS)\subseteq \ga$. Fix an element $S\in \mrL (R, \ga , A)$. Then $S\subseteq \CS$ since $\s_S (S)\subseteq A^\times$. Hence, $$\ga = \ass_R(S)\subseteq \ass_R(\CS)\subseteq \ga,$$ and so $\ga = \ass_R(\CS)$ and $\CS \in \mrL(R, \ga,  A)$. Therefore, 

(ii) $ \CS =\bigcup_{S\in \mrL (R, \ga , A)}S$: By the statement (i),   $\CS\in \mrL (R, \ga,  A)$. For all $S\in  \mrL (R, \ga , A)$, $\s (S)\subseteq A^\times$, hence $S\subseteq \CS$, and the statement (ii) follows.

By the statement (ii), the set $\CS$ is the largest element of the poset $(\mrL(R, \ga,  A), \subseteq )$, and so the poset $(\mrL(R, \ga,  A), \subseteq )$ is a directed set. By Lemma \ref{a22Dec22}, if  $S\in \mrL(R, \ga,  A)$ then $T \in \mrL(R, \ga,  A)$ for all subsets $T$ of $R$ such that $S\subseteq T\subseteq \CS (R, \ga,  A)$.

2.  By  statement 1, $\CS (R, \ga , A)= \s^{-1}(A^\times)$, and so  the set $ \CS (R, \ga , A)$ is a multiplicative set in  $R$ such that $ \CS (R, \ga , A)+\ga = \CS (R, \ga , A)$.

3. Since $\CS\in \mrL (R, \ga , A)$, there is an   $R$-isomorphism $A\simeq R\langle \CS (R, \ga , A)^{-1}\rangle$. By statement 1, the set $ \mrL (R, \ga , A)$ is a directed set w.r.t. inclusion such that for all sets $S,T\in \mrL (R, \ga , A)$ such that   $S\subseteq T$, $\RSm\simeq R\langle T^{-1}\rangle\simeq A$, $R$-isomorphisms. 
By Lemma \ref{a21Dec22},  $ \varinjlim_{S\in \mrL (R, \ga , A)}\RSm\simeq A$, an $R$-isomorphism.   $\Box$ 

Corollary \ref{a20Dec22} is a criterion for a set  $S\in \mrL (R, \ga, A )$ to be the largest element of the poset $\mrL (R, \ga, A )$.  

 \begin{corollary}\label{a20Dec22}%\marginpar{a20Dec22}
Let $S\in \mrL (R, \ga, A )$. Then $S=\CS (R, \ga, A )$ iff for every  element $r\in R\backslash S$  the rings $A$ and $R\langle (S\cup\{ r\})^{-1}\rangle$ are not  $R$-isomorphic.
\end{corollary} 

{\it Proof}. (i) {\em For all elements $x\in \CS \backslash S$, $S\cup \{ x\} \in \mrL (R, \ga, A)$}: Since $S\in \mrL (R, \ga, A )$, we must have that $S\subseteq \CS := \CS (R,\ga , A)$, by  the maximality of the element $\CS$. Let $x\in \CS \backslash S$ and  $T=S\cup \{ x\}$.  Then $S\subset T\subseteq \CS$, and so $\ga = \ass_R(S)\subseteq \ass_R(T)\subseteq \ass_R(\CS)=\ga$, i.e. $\ass_R(T)=\ga$ and there are unique $R$-homomorphisms $$A\simeq \RSm\ra R\langle T^{-1}\rangle\ra  R\langle \CS^{-1}\rangle\simeq A.$$  Hence, $T\in \mrL (R, \ga, A)$.

  (ii) {\em  For all elements $x\in R \backslash \CS$, $\CS\cup \{ x\} \not\in \mrL (R, \ga, A)$}: The statement (ii) follows from the maximality of the element $\CS$.

 Now, the corollary follows from the statements (i) and (ii). $\Box$ 

By Proposition  \ref{SA19Aug21}, for rings $A_1, A_2\in \mrLoc (R)$, $$A_1\ra A_2\;\; {\rm  iff}\;\;
\CS_l (R, \ga , A_1)\subseteq \CS_l (R, \ga , A_2).$$ 

  Corollary \ref{aSA19Aug21} is a description of  the largest  element  of the directed set $(\mrL (R, \ga , A), \subseteq )$. 
 \begin{corollary}\label{aSA19Aug21}%\marginpar{aSA19Aug21}
Let  $\CS\in \max \mrLR$, $\ga = \ass_R(\CS)$, $A=R\langle \CS^{-1}\rangle$,  and  $\s =\s_\CS : R\ra A$.  Then 
%Let  $A\in \max \mrLoc (R)$,  $A^\times$ be the group of units of the ring $A$, $\s  : R\ra A$, and $\ga =\ker (R\ra A)$.  Then 
\begin{enumerate}
\item The set $$ \CS (A):= \CS (R, \ga , A)=\s^{-1}(A^\times )=\bigcup_{S\in \mrL (R, \ga , A)}$$ is the largest element of the directed set $(\mrL(R, \ga,  A), \subseteq )$, and  if $S\in \mrL(R, \ga,  A)$ then $T\in  \mrL(R, \ga,  A)$ for all subsets $T$ of $R$ such that $S\subseteq T\subseteq \CS$.  
\item The set $ \CS (A)$ is a multiplicative set in $R$ such that $ \CS (A)+\ga = \CS ( A)$. 
\item $A\simeq R\langle \CS ( A)^{-1}\rangle\simeq  \varinjlim_{S\in \mrL (R, \ga , A)}\RSm$, $R$-isomorphisms. 
\end{enumerate}
\end{corollary} 

{\it Proof}. The corollary is a particular case of Proposition \ref{SA19Aug21}. $\Box$

 By (\ref{LocR}), for each ring $A\in \mrLoc (R)$ there is a unique ideal $\ga \in \ass \, \mrLR$ such that $A\in \mrLoc (R, \ga)$. By Proposition \ref{SA19Aug21}.(1), for each ring $A\in \mrLoc (R, \ga)$, there is the largest element $\CS (R, \ga, A)$ of the set $\mrL (R, \ga , A )$. By Proposition \ref{SA19Aug21}.(1,2), the set $$\CS (R, \ga, A))=\s^{-1}(A^\times)$$ is a multiplicative submonoid of $R$ where $\s: R\ra A$, $r\mapsto r+\ga$. Hence, the map
%\marginpar{Sm1RD2}
\begin{equation}\label{Sm1RD2}
 \mrLoc (R)= \coprod_{\ga \in \ass\,  \mrLR } \mrLoc (R,\ga ) \ra \mrLR, \;\;  \mrLoc (R, \ga )\ni A \mapsto \CS (R, \ga, A),
\end{equation}
is a section of the  epimorphism of the posets (\ref{Sm1RD}). The set 
%\marginpar{Sm1RD3}
\begin{equation}\label{Sm1RD3}
\CM (R)=\{ \CS (R, \ga, A)\, |\, \ga \in \ass\,  \mrLR, A\in \mrLoc (R, \ga )\}
\end{equation}
 is the image of the section above. In particular, the map 
%\marginpar{Sm1RD4}
\begin{equation}\label{Sm1RD4}
 \mrLoc (R)= \coprod_{\ga \in \ass\,  \mrLR } \mrLoc (R,\ga ) \ra \CM (R), \;\;  \mrLoc (R, \ga )\ni A \mapsto \CS (R, \ga, A),
\end{equation}
is an {\em isomorphism} of the posets $(\mrLoc (R), \ra )$ and $(\CM (R), \subseteq )$ with inverse $\CS\mapsto R\langle \CS^{-1}\rangle$. 

Theorem \ref{22Dec22} describes the sets  $\max \mrL (R)$ and $\max \mrLoc (R)$, presents a bijection between them, and shows that $\max \mrLoc (R)\neq \emptyset$.

\begin{theorem}\label{22Dec22}%%\marginpar{22Dec22}

\begin{enumerate}
\item The map $\max \mrLR\ra \max \mrLoc (R)$, $\CS \mapsto R\langle \CS^{-1}\rangle$ is a bijection with inverse $A\mapsto \s_A^{-1}(A^\times)$ where $\s_A: R\ra A$. 
\item  $\max \mrL (R)=\{  \s_A^{-1}(A^\times)\, | \, A\in \max \mrLoc (R)\}$.
\item $\max \mrLoc (R)=\{  R\langle \CS^{-1}\rangle\, | \, \CS\in \max \mrL (R) \}\neq \emptyset$.
\end{enumerate}
\end{theorem}

{\it Proof}. 1. Statement 1 follows from the isomorphisms of posets  given in (\ref{Sm1RD4}) and Corollary \ref{aSA19Aug21}.

2 and 3. Statements 2 and 3 follow from statement 1 and  Corollary \ref{aSA19Aug21} ($\max \mrLoc (R)\neq \emptyset$ since  $\max \mrLR\neq \emptyset$, by Theorem \ref{SC12Jan19}).(1)). $\Box $

Similarly, for each ideal $\ga\in \ass \mrLR$, Theorem \ref{A23Dec22} describes the sets  $\max \mrLRa$ and $\max \mrLoc (R, \ga)$, presents a bijection between them, and shows that they are nonempty sets.

\begin{theorem}\label{A23Dec22}%\marginpar{A23Dec22}
Let $\ga \in \ass \, \mrLR$. 
\begin{enumerate}
\item The map $\max \mrLRa\ra \max \mrLoc (R, \ga)$, $\CS \mapsto R\langle \CS^{-1}\rangle$ is a bijection with inverse $A\mapsto \s_A^{-1}(A^\times)$ where $\s_A: R\ra A$. 
\item  $\max \mrLRa=\{  \s_A^{-1}(A^\times)\, | \, A\in \max \mrLoc (R, \ga)\}\neq \emptyset$.
\item $\max \mrLoc (R, \ga)=\{  R\langle \CS^{-1}\rangle\, | \, \CS\in \max \mrLRa \}\neq \emptyset$.
\end{enumerate}
\end{theorem}

{\it Proof}. 1. Statement 1 follows from the isomorphisms of posets  given in (\ref{Sm1RD4}) and Corollary \ref{aSA19Aug21}.

2. (i) $\max \mrLRa\neq \emptyset$:  The statement (i) follows from Zorn's Lemma. In more detail,  let $\{S_i\}_{i\in \G}\subseteq \mrLRa$  be an ordered by inclusion set of localizable sets in $R$ indexed by an ordinal $\G$, i.e. for each pair $i\leq j$ in $\G$, $S_i\subseteq S_j$ and if $i \in \G$ is a limit ordinal then $S_i=\bigcup_{j<i}S_j$. In view of Zorn's Lemma, it suffices  to show that  $\bigcup_{i\in \G}S_i\in \mrLRa$. By (\ref{RSmRT}) and  (\ref{RSmRT1}), we have the direct system of ring homomorphisms, $(R\langle S_i^{-1}\rangle, \phi_{S_jS_i})$, where for each pair $i\leq j$, the map $\phi_{S_jS_i}: R\langle S_i^{-1}\rangle\ra R\langle S_j^{-1}\rangle$ is defined in (\ref{RSmRT}). Since the rings $R\langle S_i^{-1}\rangle$ are unital and $\phi_{S_jS_i}(1)=1$ for all $i\leq j$, the direct limit $\varinjlim R\langle S_i^{-1}\rangle$  is a nonzero ring and  and is isomorphic to $R\langle \Big(\bigcup_{i\in \G}S_i\Big)^{-1}\rangle$. Therefore, $\bigcup_{i\in \G}S_i\in \mrL (R)$. Finally, $$\ann_R(\bigcup_{i\in \G}S_i)=\bigcup_{i\in \G}\ann_R(S_i)=\bigcup_{i\in \G}\ga  = \ga,$$ and the statement (ii) follows. 

(ii)  $\max \mrLRa=\{  \s_A^{-1}(A^\times)\, | \, A\in \max \mrLoc (R, \ga)\}$: The statement (ii) follows from statement 1 and  Corollary \ref{aSA19Aug21}.

3.  Statement 3 follows from statement 1 and  Corollary \ref{aSA19Aug21} (by statements 1 and 2, $\max \mrLoc (R, \ga)\neq \emptyset$). $\Box $\\

{\bf The $\mrL$-radical $\mrL\rad (R)$ of the ring $R$, the sets of localizable, non-localizable and completely localizable elements of $R$.} 

{\it Definition.} An element $S$ of the set	$\max  \mrLR$  is called a {\em maximal  localizable set}, and the rings $\RSm$ is  called the {\em maximal  localization} of $R$. The set of maximal localizations of the ring $R$ is denoted my $\max\mrLoc (R)$. The intersection and the sum, 
%\marginpar{SLLradR}
\begin{equation}\label{SLLradR}
\gll (R):=\mrL\rad (R) := \bigcap_{S\in \max  \mrLR} \ass _R(S)\;\; {\rm and}\;\; \gc\gll (R):= {\rm co-}
\mrL\rad (R) := \sum_{S\in \max  \mrLR} \ass _R(S), 
\end{equation}
are called the {\em localization radical} and the {\em localization coradical}, respectively, or the $\mrL$-{\em radical} and  the $\mrL$-{\em coradical} of $R$, for short. 
By Theorem \ref{SC12Jan19}.(2), $$\gc\gll (R)= {\rm co-}
\mrL\rad (R) = \sum_{S\in  \mrLR} \ass _R(S).$$
For the ring $R$, there is the canonical exact sequence,  
%\marginpar{SLLradR1}
\begin{equation}\label{SLLradR1}
0\ra \mrL\rad (R)\ra R\stackrel{\s }{\ra} \prod_{S\in \max \mrLR}\RSm , \;\; \s := \prod_{S\in\max \mrLR}\, \s_S,
\end{equation}
where $\s_S:R\ra \RSm$. \\

{\it Definition.} The sets $\CL \mrLR:=\bigcup_{S\in \mrLR}S=\bigcup_{S\in \max \mrLR}S$ and $\CN \CL\mrLR := R\backslash \CL \mrLR$ are called the set of {\em localizable} and {\em non-localizable} elements of $R$, resp., and the intersection
$$ \CC\mrLR :=\bigcap_{S\in \max \mrLR}S$$ is called the {\em set of completely  localizable elements} of the ring $R$. The ring $Q_c(R):= R\langle \CC\mrLR^{-1}\rangle$ is called the {\em complete localization} or the $\CC$-{\em localization}   of $R$. Let $\gc_R:= \ass_R(\CC\mrLR)$, $R_c:=R/\gc_R$,   $\pi_{\gc_R}: R\ra R_c$, $r\mapsto r+\gc_R$, and $\CC_c:=\CC\mrLR_c :=\pi_{\gc_R}(\CC\mrLR)$. \\

By Corollary \ref{a17Dec22}.(2), $R_c\subseteq Q_c$, $\CC_c\in \mrL (R_c)$, and 
%\marginpar{QcR}
\begin{equation}\label{QcR}
Q_c(R)\simeq R_c\langle \CC_c^{-1}\rangle. 
\end{equation}

By the  definition,  the sets $\CL \mrLR$, $\CN \CL\mrLR$,  and $ \CC\mrLR$ are invariant under the action of the automorphism  group of the ring $R$, i.e.,  they are {\em characteristic sets}. 

Proposition \ref{A29Dec22} shows that there are tight connections between the sets $\CL \mrLR$, $\CN \CL\mrLR$,  and $ \mrL\rad (R)$.

\begin{proposition}\label{A29Dec22}%\marginpar{A29Dec22}

\begin{enumerate}
\item $\CL\mrLR \cap \mrL\rad (R)=\emptyset$.
\item $\mrL\rad (R)\subseteq \CN\CL\mrLR $.
\item $\CL\mrLR + \mrL\rad (R)=\CL\mrLR$.
\item $\CN\CL\mrLR + \mrL\rad (R)=\CN\CL\mrLR$.
\end{enumerate}
\end{proposition}

{\it Proof}. 1. For each element $s\in \CL\mrLR$, there is a maximal localizable set, say $S\in \max \mrLR$, such that $s\in S$. Then $s\not\in \ass_R(S)$ since the element $s$ is a unit in the ring $\RSm$. Hence, $s\not\in \mrL\rad (R)$, and so $\CL\mrLR \cap \mrL\rad (R)=\emptyset$.

2. Statement 2 follows from statement 1:  $\mrL\rad (R)\subseteq R\backslash \CL\mrLR = \CN\CL\mrLR $. 

3. Since $0\in \mrL\rad (R)$, to prove statement 3 it suffices to show that $\CL\mrLR + \mrL\rad (R)\subseteq \CL\mrLR$. For each element $s\in  \CL\mrLR$, there is a maximal localizable set, say $S\in \max \mrLR$, that contains the element $s$. Then the image of the elements $s+\mrL\rad (R)$ in the ring $\RSm$ are units (since $\mrL\rad (R)\subseteq\ass_R(S)$). Therefore, $s+\mrL\rad (R)\subseteq \CL\mrLR$, and statement 3 follows.

4. Recall that  $R=\CL\mrLR\coprod \CN\CL\mrLR$, a disjoint union. Suppose that statement 4 does not hold, i.e. there are elements $n\in \CN\CL\mrLR$ and $r\in \mrL\rad (R)$ such $n+r\not\in \CN\CL\mrLR$. Then $n+r\in \CL\mrLR$, and so 
$$\CL\mrLR \not\ni n\in -r+\CL\mrLR\subseteq \mrL\rad (R)+\CL\mrLR\stackrel{{\rm st.3}}{=}\CL\mrLR,$$ a contradiction. $\Box $

Corollary \ref{a28Dec22} provides a useful sufficient condition for the localization radical to be zero. Corollary  \ref{a28Dec22} is used in the description of the localization radical of an arbitrary direct product of division rings (Theorem \ref{A28Dec22}.(3)).

\begin{corollary}\label{a28Dec22}%\marginpar{a28Dec22}
If $\CL\mrLR =R\backslash \{ 0\}$  then $\mrL\rad (R)=0$.
\end{corollary}

{\it Proof}. By Proposition \ref{A29Dec22}.(2), $\mrL\rad (R)\subseteq \CN\CL\mrLR =R\backslash \CL\mrLR =\{0\}$ since $\CL\mrLR =R\backslash \{ 0\}$. Hence $\mrL\rad (R)=0$.  $\Box$

Lemma \ref{a19Dec22} is a criterion for an element of $R$ to be a localizable element. 
\begin{lemma}\label{a19Dec22}%\marginpar{a19Dec22}
Let $r\in R$. The following statements are equivalent:
\begin{enumerate}
\item $r\in \CL\mrLR$.
\item $\{r\}\in \mrLR$. 
\item $\ass_R(r)\neq R$.
\item There is an $R$-homomorphism of rings $f: R\ra A$ such that $f(r)\in A^\times$. 
\end{enumerate}
\end{lemma}

{\it Proof}. $(1\Rightarrow 2 )$ If $r\in S$ for some $S\in \mrLR$, then $\{r\}\in \mrLR$, by Proposition \ref{A7Dec22}.(2).

$(2\Rightarrow 1 )$ The implication follows from the inclusions $r\in \{ r\}$ and $ \{ r\} \in \mrLR$.

$(2\Leftrightarrow 3 \Leftrightarrow 4 )$ Proposition \ref{a4Dec22}. $\Box $

\begin{corollary}\label{a25Dec22}%\marginpar{a25Dec22}
\begin{enumerate}
\item Every localizable set consists of localizable elements.
\item Every subset of a ring that consists of localizable elements is a localizable set iff there is  the largest localizable set.
\end{enumerate}
\end{corollary}

{\it Proof}. 1. Statement 1 follows from  Proposition \ref{A7Dec22}.(2) or  Lemma \ref{a19Dec22}.(1,2).

2. If $\CS$ is the largest localizable set of a ring $R$ then  $\CL\mrLR =\CS$,  and so  every subset of $R$  that consists of localizable elements is a subset of $\CS$, hence localizable, by Proposition \ref{A7Dec22}.(2).

Suppose that  a ring $R$ has two maximal localizable sets, $\CS$ and $\CT$. Then their union consists of localizable elements but it is not a localizable set otherwise, by Theorem \ref{SC12Jan19}.(2),  it would have contained in a maximal localizable set, which is not possible (since $\CS$ and $\CT$ are  maximal localizable sets). $\Box $\\

{\bf The set $\CC\mrLR$ of completely localizable elements of $R$.} Let $\gc_R :=\ass_R(\CC\mrLR)$. Consider the maps
%\marginpar{gcLradR}
\begin{equation}\label{gcLradR}
0\ra \mrL\rad (R)\ra R \stackrel{\pi_{\gc_R}}{\ra} R/\gc_R\stackrel{\overline{\kappa} }{\ra} Q_c(R)\ra  \prod_{S\in \max \mrLR}\RSm \;\; {\rm and} \;\; \kappa :=\overline{\kappa} \pi_{\gc_R}:R\ra  Q_c(R)
\end{equation}
is the unique $R$-homomorphism from $R$ to its localization $Q_c(R)$, $\pi_{\gc_R}(r):=r+\gc_R$,  and $\overline{\kappa}$ is the unique $R$-homomorphism from $R/\gc_R$ to $Q_c(R)$ (induced by $\kappa$).

Proposition \ref{B29Dec22} contains some properties of $\CC\mrLR$ and $\gc_R$.

\begin{proposition}\label{B29Dec22}%\marginpar{B29Dec22}
Let $\CC :=\CC\mrLR$. 
\begin{enumerate}
\item $\CC=\s^{-1}\Big(  \prod_{S\in \max \mrLR}\RSm^\times \Big)$ where $\s :R\ra  \prod_{S\in \max \mrLR}\RSm$.
\item $\CC\mrLR + \mrL\rad (R)=\CC\mrLR$.
\item $\gc_R\subseteq \mrL\rad (R)$.
\item If the map $\overline{\kappa}:R/\gc_R\ra Q_c(R) $  is an injection then  $\gc_R = \mrL\rad (R)$.
\item $\CC +\gc_R=\CC$.
\end{enumerate}
\end{proposition}

{\it Proof}. 1. Statement 1 follows from the definition of the set $\CC$ and 
Corollary \ref{aSA19Aug21}.(1): 
 $$\CC=\bigcap_{S\in \max \mrLR}S=\s^{-1}\bigg(  \prod_{S\in \max \mrLR}\RSm^\times \bigg).$$

2. By Corollary \ref{aSA19Aug21}.(3), $S+\ass_R(S)=S$ for all $S\in \max \mrLR$. Now,  $S+\mrL\rad (R)=S$ for all $S\in \max \mrLR$  since $\mrL\rad (R)=\bigcap_{S\in \max \mrLR}\ass_R(S)$. It follows that  $\CC\mrLR + \mrL\rad (R)=\CC\mrLR$ since $\CC\mrLR = \bigcap_{S\in \max \mrLR} S$. 

3. The inclusion $\gc_R\subseteq \mrL\rad (R)$ follows from (\ref{gcLradR}).

4. If the map $\overline{\kappa}$  is an injection then  $\gc_R\supseteq \mrL\rad (R)$ (by (\ref{gcLradR})), and statement 4 follows from statement 3.

5. Statement 5 follows from statements 2 and 3. $\Box $\\

{\bf The $\ga$-absolute quotient ring $Q_a(R, \ga))$ of a ring $R$ where $\ga \in \ass\, \mrLR$.} Let $\ga \in \ass \, \mrLR$ and 
%\marginpar{URa}
\begin{equation}\label{URa}
U(R,\ga):=\bigcup_{S\in \max \mrLRa}S.
\end{equation}
The  set  $(\mrLRa, \subseteq )$   is a poset and the system $\{ R\langle S^{-1}\rangle, \phi_{ST}\}$ is a direct system of rings where $S, T\in \mrLRa$.

{\it Definition.} The direct limit $Q_a(R,\ga):=\varinjlim_{S\in  \mrLRa}\RSm$ is called the $\ga$-{\em absolute quotient ring}  of the ring $R$. Let $\s_{a,\ga}: R\ra Q_a(R,\ga)$ be a unique $R$-homomorphism and  $\ga_{a,\ga}(R):=\ker (\s_{a,\ga})$, $R_{a,\ga}:=R/\ga_{a,\ga}(R)$,  $\pi_{a,\ga}:R\ra R_{a,\ga}$, $r\mapsto r+\ga_{a,\ga}$, and $\overline{\s_{a,\ga}}:R_{a,\ga}\ra Q_a (R,\ga)$ is a unique $R$-homomorphism. Clearly, $\s_{a,\ga}=\overline{\s_{a,\ga}}\pi_{a,\ga}$ and $\ga_{a,\ga}(R)\supseteq \ga$.

Theorem \ref{aXA30Dec22} gives an explicit description of the ring $Q_a(R,\ga)$, a criterion for $Q_a(R,\ga)\neq \{0\}$, and a criterion for $\ga_{a,\ga}=\ga$.

\begin{theorem}\label{aXA30Dec22}%\marginpar{aXA30Dec22}
Let $\ga \in \ass \,\mrLR$ and $U=U(R,\ga )$.
\begin{enumerate}
\item $Q_a(R,\ga)\neq \{0\}$ iff $U\in \mrLR$.
\item Suppose that $U\in \mrL (R, \gb)$. Then 
\begin{enumerate}
\item The rings $Q_a(R,\ga)$ and $ R\langle U^{-1}\rangle$ are $R$-isomorphic and $\ga_{a,\ga}=\gb$.
\item $\ga_{a,\ga}=\ga$ iff $\max \mrLRa =\{ S\}$ iff the rings $Q_a(R,\ga)$ and $ R\langle S^{-1}\rangle$ are $R$-isomorphic for some $S\in \max\mrLRa$.
\end{enumerate}
\end{enumerate}
\end{theorem}

{\it Proof}.  1.  $(\Rightarrow)$ Suppose that $Q_a(R,\ga )\neq \{ 0\}$. Then 
 the set $\CS :=\s_{a,\ga}^{-1}(Q_a(R,\ga )^\times)$ is a localizable set in $R$ since $\s_{a,\ga}(\CS)\subseteq Q_a(R,\ga )^\times$.  Then for each $S\in \max\mrLRa$, 
 there is a commutative  diagram
$$\begin{array}[c]{ccc}
R&\stackrel{\s_S}{\ra} &\RSm\\
&\stackrel{\s_{a,\ga}}{\searrow}&{\downarrow}^{\phi_S}\\
&  & Q_a(R,\ga)
\end{array}$$
where $\phi_S$ is a unique $R$-homomorphism. Then it follows from the inclusion $\phi_S (\RSm^\times)\subseteq  Q_a(R,\ga )^\times$ that  
$$ \CS =\s_{a,\ga}^{-1}(Q_a(R,\ga )^\times)=(\phi_S\s_S)^{-1}(Q_a(R,\ga )^\times)=\s_S^{-1}\phi_S^{-1}(Q_a(R,\ga )^\times)\supseteq \s_S^{-1}(\RSm^\times)=S,$$
by Corollary \ref{aSA19Aug21}.(1).
Therefore, $\CS \supseteq \bigcup_{S\in \max \mrLRa}S=U$, and so $U\in\mrLR $ since $\CS \in \mrLR$.

$(\Leftarrow)$  Suppose that $U\in \mrLR$. Let us show that the 
rings $Q_a(R,\ga)$ and $R\langle U^{-1}\rangle$ are $R$-isomorphic. By the definition,  the ring $Q_a(R,\ga)$ is the direct limit of the direct system $\CD:= \{ R\langle S^{-1}\rangle, \phi_{ST}\}$ where $S, T\in \mrLRa$. Since for all $S\in \mrLRa$, $S\subseteq U$, we can extend this direct system to a direct system $\CD'$ by adding the $R$-homomorphisms 
$\RSm\ra R\langle U^{-1}\rangle$ for all $S\in \mrLRa$.
The element $U$ is the largest element of the poset $\mrL':=\mrLRa\cup \{ U\}$. By Theorem \ref{A14Dec22}.(3), 
$$\varinjlim_{S\in \mrL'}\RSm\simeq R\langle U^{-1}\rangle,$$ an $R$-isomorphism. By the definitions of the two direct systems, $\CD$ and $\CD'$, there is necessarily a unique $R$-homomorphism $\alpha : Q_a(R,\ga)\ra R\langle U^{-1}\rangle$. 
Since 
$$\s_{a,\ga} (U)\subseteq Q_a(R, \ga)^\times,$$ 
where $\s_{a,\ga}:R\ra Q_a(R, \ga)$,  there is necessarily a unique $R$-homomorphism 
$$\beta :  R\langle U^{-1}\rangle\ra Q_a(R,\ga).$$  
Hence, the rings $Q_a(R,\ga)$ and $R\langle U^{-1}\rangle$ are $R$-isomorphic, by the universal properties of localization and the direct limit. In particular, $\ga_{a,\ga}=\gb$. This finishes the proof of statements 1 and 2(a).

2(b). (i) $\ga_{a,\ga}=\ga \Rightarrow\max \mrLRa =\{ S\}$: Suppose that the set $\max\mrLRa$ contains at least two  elements, say $S$ and $T$. Then $$S\subset S\cup T\subseteq U\;\; {\rm  and}\;\;S,U\in \mrLRa.$$ This contradicts to the fact that $S\in 
\max \mrLRa$, and the statement (i) follows.

(ii) $\max \mrLRa =\{ S\}\Rightarrow $ {\em the rings $Q_a(R,\ga)$  and $ R\langle S^{-1}\rangle$ are $R$-isomorphic}: The implication follows from Theorem  \ref{A14Dec22}.(3).

(iii) {\em The rings $Q_a(R,\ga)$ and $ R\langle S^{-1}\rangle$ are $R$-isomorphic for some $S\in \max\mrLRa$ then }
 $\ga_{a,\ga}=\ga$: The implication is obvious. $\Box$

{\it Example.} If $R$ is a commutative domain then the set $\CC_R=R\backslash \{ 0\}$ is the largest element of $\mrL (R,0)$. By Theorem  \ref{aXA30Dec22},  $Q_a(R, 0)=Q_{cl}(R)$ is the classical quotient ring of $R$, i.e. the field of fractions of $R$, and $\ga_{a,0}(R)=0$. 

\begin{lemma}\label{axa30Dec22}%\marginpar{axa30Dec22}
If  $R$ is  a  simple Artinian rings then $\max\mrL (R,0)=\{ R^\times\}$,  $Q_a(R,0 )=R$ and $\ga_{a,0}(R)=0$.
\end{lemma}

{\it Proof}. Lemma follows from Lemma \ref{xa30Dec22}. $\Box$\\

{\bf Isomorphism criterion for localizations of a ring.} Theorem \ref{20Dec22}.(1) is an $R$-isomorphism  criterion for  localizations of the ring $R$ and Theorem \ref{20Dec22}.(2) is a criterion for existing of an $R$-homomorphism between localizations of $R$.

\begin{theorem}\label{20Dec22}%\marginpar{20Dec22}
Let $S, T\in \mrL (R)$. 
\begin{enumerate}
\item The rings $\RSm$ and $R\langle T^{-1}\rangle$ are $R$-isomorphic iff $$\s_S^{-1}(\RSm^\times)=\s_T^{-1}(R\langle T^{-1}\rangle^\times)$$ 
where $\s_S: R\ra \RSm$ and $\s_T: R\ra R\langle T^{-1}\rangle$.
\item There is an $R$-homomorphism $\RSm\ra R\langle T^{-1}\rangle$ iff $\s_S^{-1}(\RSm^\times)\subseteq \s_T^{-1}(R\langle T^{-1}\rangle^\times).$
\end{enumerate}

\end{theorem}

{\it Proof}. 1.  $(\Rightarrow)$ The implication is obvious.

$(\Leftarrow )$ The implication follows from the fact that the rings $\RSm$ and $R\langle \s_S^{-1}(\RSm^\times)\rangle$ are isomorphic for all $S\in \mrL (R)$, by Proposition \ref{SA19Aug21}.(1). 

2. $(\Rightarrow)$ The implication is obvious.

$(\Leftarrow )$ The implication follows from the fact that the rings $\RSm$ and $R\langle \s_S^{-1}(\RSm^\times)\rangle$ are isomorphic for all $S\in \mrL (R)$ (Proposition \ref{SA19Aug21}.(1)) and the fact that if $S_1,S_2\in \mrL (R)$ and $S_1\subseteq S_2$ then there is a unique $R$-homomorphism $R\langle S_1^{-1}\rangle \ra R\langle S_2^{-1}\rangle$.  $\Box $

 Proposition \ref{SA20Aug21} is a practical criterion for  a subset $S$  of a ring $R$ to be an element of the set $\mrL(R, \ga )$.

\begin{proposition}\label{SA20Aug21}%\marginpar{SA20Aug21}
Let $S$ be a subset of a ring $R$.
\begin{enumerate}
%\item Suppose that $S\in \mrL (R, \ga )$. Then the ideal $\ga$ is the least ideal of $R$ such that $S+\ga \in \Den(R/\ga , 0)$ where $*\in \{ l,r,\emptyset\}$.
\item Suppose that there exists an ideal $\ga$ of $R$ such that $\ga \subseteq \ass_R(S)$ and $\bS :=\pi (S) \in  \mrL(\bR , 0)$ where  $\pi : R\ra \bR:=R/\ga$, $r\mapsto r+\ga$. Then  $S\in \mrL(R, \ga )$.

\item $S\in \mrL(R, \gb)$ iff there is an ideal $\ga$ of $R$ such that $\ga \subseteq \gb$ and  $\bS :=\pi (S)\in  \mrL(\bR , 0)$  where  $\pi : R\ra \bR:=R/\ga$, $r\mapsto r+\ga$.
\end{enumerate}

\end{proposition} 

{\it Proof}. % 1. By Theorem \ref{L11Jan19}.(1),  the ideal $\ass_R(S)$ satisfies the condition of statement 4. Conversely, suppose that an ideal $\ga'$ of the ring $R$ be such that $\bS':=S+\ga'\in \Den(\bR', 0)$ where $\bR':=R/\ga'$. Since the elements of the set $\bS'$ are invertible in the localization $\CR'$ of the ring $\bR'$ at $\bS'$, there is an $R$-epimorphism from $R\langle S^{-1}\rangle$ to $\CR'$. In particular, $\ass_R(S)\subseteq \ga'$, and statement 4 follows. 
1. Since the elements of the set $\bS$ are invertible in the ring $\bR\langle \bS^{-1}\rangle$, $\ass_R(S)\subseteq \ga$, by Proposition \ref{A18Dec22}. By the assumption,  $\ass_R(S)\supseteq \ga$.  Therefore, $\ass_R(S)=\ga$ and  $S\in \mrL(R, \ga )$. 

2. $(\Rightarrow)$ If $S\in \mrL(R, \gb )$ then it suffices to take $\ga=\gb$, by Corollary \ref{a17Dec22}.(2).

$(\Leftarrow)$ This implication follows from statement 1. $\Box$

 \begin{proposition}\label{TA19Aug21}%\marginpar{TA19Aug21}
 Suppose that $S\in \mrLRa$,  $T\in \mrL (\RSm, 0)$ and   $\bS \subseteq T$ where $\bS =\s_S (S)$ and $\s_S: R\ra \RSm$.   Then 
 \begin{enumerate}
\item $\bT :=T\cap \bR\in  \mrL (\bR , 0)$, $\bT \in  \mrL (\RSm , 0)$, $$\bR\langle \bS^{-1}\rangle \subseteq \bR \langle \bT^{-1}\rangle =\RSm \langle \bT^{-1}\rangle\subseteq \RSm \langle T^{-1} \rangle,$$ and $T':= \s_S^{-1} (T) \in \mrL(R, \ga ,\bR \langle \bT^{-1}\rangle )$.
\item Let $\s: R\ra A:=\RSm \langle T^{-1} \rangle$ and $\CT := \s^{-1}(A^\times)$. Then $\CT \in   \mrL(R, \ga ,A )$.

\end{enumerate}
\end{proposition} 

{\it Proof}. 1.  Recall that the rings $\RSm$ and $\bR\langle \bS^{-1}\rangle$ are $R$-isomorphic. 

(i) $\bT \in \mrL (\bR , 0)$ {\em and} $\bT \in \mrL (\RSm , 0)$: The  chain of inclusions below follows from the inclusions $\bR\subseteq \bR\langle \bS^{-1}\rangle$ and $\bT\subseteq T$, 
$$ \ass_{\bR}(\bT )\subseteq \ass_{\bR\langle \bS^{-1}\rangle}(\bT)=\ass_{\RSm}(\bT) \subseteq \ass_{\RSm}(T)=0,$$
and the statement (i) follows.

(ii) $\bS \subseteq \bT$ {\em and} $\bR\langle \bS^{-1}\rangle \subseteq \bR \langle \bT^{-1}\rangle =\RSm \langle \bT^{-1}\rangle\subseteq \RSm \langle T^{-1} \rangle$: Since $\bS\subseteq T$ and $\bS\subseteq \bR$, we have that $\bS \subseteq T\cap \bR =\bT$. Hence $$\bR\langle \bS^{-1}\rangle \subseteq \bR \langle \bT^{-1}\rangle $$ since  $\bT\in \mrL (\RSm, 0)=\mrL (\bR\langle \bS^{-1}\rangle , 0)$ (the statement (i)).  Now,
$$  \bR \langle \bT^{-1}\rangle = \bR \langle (\bS\cup \bT )^{-1}\rangle= \bR \langle \bS^{-1}\rangle  \langle \bT^{-1}\rangle =\RSm \langle \bT^{-1}\rangle\subseteq  \RSm \langle T^{-1}\rangle $$
since $\bS\subseteq \bT\subseteq T$  and $\bT , T\in \mrL (\RSm, 0)$. 

Recall that $\pi_S: R\ra \bR$, $r\mapsto r+\ga$. 

(iii) $T'=\pi_S^{-1}(\bT)\supseteq S$ {\em and} $\pi_S(T')=\bT$: The map $\s_S: R\ra \RSm$ is the composition of two maps $$\s_S: R\stackrel{\pi_S}{\ra}\bR\ra \RSm$$ where the second map is a natural inclusion. Hence,  $T'=\pi_S^{-1}(\bT)$. Clearly,  $\pi_S^{-1}(\bT)\supseteq S$ since $\bS\subseteq T$ and $\pi_S(S)=\bS$.
 The map $\pi_S: R\ra \bR$ is an epimorphism, and so $\pi_S(T')=\pi_S(\pi_S^{-1}(\bT))=\bT$.

(iv) $T'\in \mrL(R, \ga ,\bR \langle \bT^{-1}\rangle )$:
 It follows from the inclusion $S\subseteq T'$ that $$\ga = \ass_R(S)\subseteq \ass_R(T').$$ By the statements (i) and (iii), $\pi_S(T')=\bT\in \mrL (\bR , 0)$. Therefore, $\ass_R(T')=\ga$, by Proposition \ref{SA20Aug21}.(1). Now, 
$$ R\langle T'{}^{-1}\rangle =  R\langle (S\cup T)'{}^{-1}\rangle = \RSm\langle \pi_S(T')^{-1}\rangle=  \RSm\langle \bT^{-1}\rangle =\bR\langle \bT^{-1}\rangle,$$
by the statement (ii), and the statement (iv) follows.

2. Since $T\in \mrL (\RSm, 0)$, $\ker (\s) = \ga$. Now, statement 2 follows from Proposition \ref{SA19Aug21}.(1). $\Box$\\

{\bf Localizable sets and denominator sets.} Theorem \ref{A17Jan23} is a criterion for the ring $\RSm$ to be an $R$-isomorphic to a localization of $R$ at a denominator set and  it also  describes all such denominator sets.

\begin{theorem}\label{A17Jan23}%\marginpar{A17Jan23}
Let $S\in \mrLRa$,   $\s_S:R\ra \RSm$, and $\CS := \s_S^{-1}(\RSm^\times)$. Then 
\begin{enumerate}
\item The ring $\RSm$ is  $R$-isomorphic to a localization of the ring $R$ at a *-denominator set where $*\in \{ l,r, \emptyset\}$ iff $\CS\in \Den_*(R,\ga)$.
\item The set of all *-denominator sets of $R$  with localization $R$-isomorphic to the ring $\RSm$ is equal to $\{ T\in \Den_*(R,\ga )\, | \, T\subseteq \CS, \CS+\ga \subseteq (T^{-1}R)^\times\}$.

\end{enumerate}
\end{theorem}

{\it Proof}. 1. Let us prove statement 1 for left denominator sets. The other two case can be proven in a similar way (or  reduced to this case using an opposite ring construction). 

 $(\Rightarrow )$ Suppose that the ring $\RSm$ is  $R$-isomorphic to a localization of the ring $R$ at a left denominator set $T$ of $R$. Then $T\in \Den _l(R, \ga )$. By \cite[Proposition 3.1.(1)]{larglquot}, $\CS \in \Den _l(R, \ga )$ and the rings  $\CS^{-1}R$ and $T^{-1}R$ are $R$-isomorphic.

$(\Leftarrow )$ Suppose that $\CS \in \Den _l(R, \ga )$. Then the rings    $\CS^{-1}R$ and $R\langle \CS^{-1}\rangle$ are $R$-isomorphic.  By Proposition \ref{SA19Aug21}.(1), the rings    $\RSm$ and $R\langle \CS^{-1}\rangle$ are $R$-isomorphic, and the implication follows. 

2. Statement follows from Corollary \ref{a4Jan23}.  $\Box $

For a denominators set $T\in \Den_*(R, \ga)$, Proposition \ref{B17Jan23}  describes all localizable sets $S$ of the ring $R$ such that the ring $\RSm$ is an $R$-isomorphic to the localization of $R$ at $T$.

\begin{proposition}\label{B17Jan23}%\marginpar{B17Jan23}
Let $T\in \Den_*(R,\ga )$ where $*\in \{ l,r, \emptyset\}$, $A$ be the localization of $R$ at $T$, and  $\pi_T:R\ra \bR = R/\ga$, $r\mapsto r+\ga$. Then the set of all localizable sets $S$ of $R$ such that the rings $\RSm$ and $A$ are $R$-isomorphic is equal to $\{ S\in \mrLRa\, | \, \pi_S(S)\subseteq A^\times,  \pi_S(T)\subseteq \RSm^\times\}$.
\end{proposition}

{\it Proof}. The statement follows from Corollary \ref{a4Jan23}.  $\Box $

For a localizable  multiplicative   set $S$ of the ring $R$, Proposition \ref{C17Jan23} is a criterion for the set $S$ to be a left denominator set of $R$.

\begin{proposition}\label{C17Jan23}%\marginpar{C17Jan23}
A localizable multiplicative  set  $S\in \mrLRa$ is a left denominator set of $R$ iff
\begin{enumerate}
\item $\ga = \sum_{s\in S}\ker_R (s\cdot )\supseteq \sum_{s\in S}\ker_R (\cdot s)$, and 
\item the left $R$-module $\RSm/\bR$ is $S$-torsion where $\bR = R\ga$.
\end{enumerate} 
\end{proposition}

{\it Proof}.  $(\Rightarrow )$ The implication is obvious.

$(\Leftarrow )$ Recall that $\ass_l(S)=\sum_{s\in S}\ker_R (s\cdot )$. By statement 1, 
$$\ga = \ass_l(S).$$ 
Let $a\in \RSm$. By statement 2,  $sa=\br\in \bR$ for some elements $s\in S$ and $r \in \bR$, and so the element $a=s^{-1}\br$ is a left fraction.  Therefore, the set $S$ is a left Ore set since $\ga = \ass_l(S)$. In fact, $S\in \Den_l(R,\ga)$ since $\ass_l(S)\supseteq \sum_{s\in S}\ker_R (\cdot s)$, see statement 1.  $\Box $

Corollary \ref{aC17Jan23} (resp.,  Corollary \ref{bC17Jan23}) is a criterion for the set $S$ to be a right  (resp., left and right) denominator set of $R$.

\begin{corollary}\label{aC17Jan23}%\marginpar{aC17Jan23}
A localizable multiplicative  set  $S\in \mrLRa$ is a right denominator set of $R$ iff
\begin{enumerate}
\item $\ga = \sum_{s\in S}\ker_R (\cdot s)\supseteq \sum_{s\in S}\ker_R (s\cdot )$, and 
\item the right $R$-module $\RSm/\bR$ is $S$-torsion where $\bR = R\ga$.
\end{enumerate} 
\end{corollary}

{\it Proof}. Repeat the proof of Proposition \ref{C17Jan23} replacing `left' to `right' everywhere.  $\Box $

\begin{corollary}\label{bC17Jan23}%\marginpar{bC17Jan23}
A localizable multiplicative  set  $S\in \mrLRa$ is a  denominator set of $R$ iff
\begin{enumerate}
\item $\ga = \sum_{s\in S}\ker_R (\cdot s)= \sum_{s\in S}\ker_R (s\cdot )$, and 
\item  $\RSm/\bR$ is an $S$-torsion left and  right $R$-module where $\bR = R\ga$.
\end{enumerate} 
\end{corollary}

{\it Proof}. The corollary follows from Proposition \ref{B17Jan23} and Corollary \ref{aC17Jan23}.  $\Box $

%%%%%%%%%%%%%%%%%   Section  4   %%%%%%%%%%%%%%%%%%%

\section{Localizations of commutative  rings}\label{COMRINGS}%\marginpar{COMRINGS}

In this section, let $R$ be a {\em commutative} ring and $\gn=\gn_R$ be its prime radical.

For an arbitrary commutative ring $R$, descriptions of the following sets are obtained: the sets of localizable and non-localizable elements (Lemma \ref{a29Dec22}), the maximal localizable sets, maximal localization rings, the set of completely localizable elements, and the ideal $\gc_R$  (Theorem \ref{C29Dec22}). Theorem \ref{C29Dec22}  also describes relations between the  ideals $\mrL\rad (R)$, $\gn_R$, and $\gc_R$.

For a commutative ring $R$ such that  $|\min (R)|<\infty$ and  $\gn_R$ is a nilpotent ideal (eg, $R$ is a commutative Noetherian ring), Proposition \ref{XC29Dec22} describes  $\Spec (Q_c(R))$ and  the rings  $Q_c (R)$ and $Q_c(R)/\gn_{Q_c(R)}$.

For a commutative ring $R$ such that  $|\min (R)|<\infty$ and  $\gn_R=0$  (eg, $R$ is a commutative semiprime Noetherian ring), Proposition \ref{XC29Dec22} shows that  $\CC\mrLR = \CC_R$, $Q_c(R)=Q_{cl}(R)$, and $\mrL\rad (R)=\gc_R=\gn_{Q_c(R)}=0$.   A characterization of commutative semiprime Goldie rings is given (Corollary \ref{a1Jan23}).

 Notice that the prime radical $\gn_R$ of a commutative ring $R$ contains precisely all the nilpotent elements of the ring $R$.  
 Lemma \ref{a29Dec22} describes the sets  of localizable and non-localizable elements of an arbitrary commutative ring.

\begin{lemma}\label{a29Dec22}%\marginpar{a29Dec22}
Let $R$ be a commutative ring.
\begin{enumerate}
\item $\CL\mrLR =R\backslash \gn_R$.
\item $\CN\CL\mrLR =\gn_R$.
\item Let $S\subseteq R$. Then $S\in \mrLR$ iff $S_{mon}$ is a multiplicative subset of $R$. If $S\in \mrLR$ then $\RSm =S_{mon}^{-1}R$ and $\ass_R(S)=\ass_R (S_{mon})=\{ r\in R\, | \, tr=0$ for some $t\in S_{mon}\}$. 
\end{enumerate}
\end{lemma}

{\it Proof}. 1 and 2.  The prime radical $\gn_R$ contains precisely all the nilpotent elements of the ring $R$. Hence,  $\gn_R\subseteq \CN\CL\mrLR$. The ring $R$ is a commutative. So, for each element $t\in R\backslash \gn_R$, the set $S_t:=\{ t^i\, | \, i\in \N\}$ is a multiplicative set of $R$, and so $s\in \CL\mrLR$, i.e. $R\backslash \gn_R\subseteq \CL\mrLR $. Since, $$R= \CL\mrLR \coprod  \CN\CL\mrLR =R\backslash \gn_R\coprod \gn_R,$$ we must have that 
 $\CL\mrLR =R\backslash \gn_R$ and  $\CN\CL\mrLR =\gn_R$. 
 
 3. Statement 3 follows from statement 1. $\Box $

For an arbitrary commutative ring $R$, Theorem \ref{C29Dec22} 
describes its maximal localizable sets, maximal localization rings, the set of completely localizable elements, and the ideal $\gc_R$. It also describes relations between the  ideals  $\mrL\rad (R)$, $\gn_R$, and $\gc_R$.

\begin{theorem}\label{C29Dec22}%\marginpar{C29Dec22}
Let  $R$ be a commutative ring and  $\CC := \CC \mrLR$. Then 
\begin{enumerate}
\item $\max \mrL (R)=\{S_\gp \, | \,\gp\in \min (R)\}$ where 
$S_\gp =R\backslash \gp$ is a multiplicative subset of $R$.
\item $R\langle S_\gp^{-1}\rangle =S_\gp^{-1}R:=R_\gp$ is the localization of the ring $R$ at the prime ideal $\gp$.
\item $\ass_R(S_\gp )=\{ r\in R\, | \, sr=0$ for some $s\in S_\gp\}\subseteq \gp$. 
\item $\mrL\rad (R)\subseteq \gn_R$.  
\item $\CC = R\backslash \Big( \bigcup_{\gp \in \min (R)}\gp\Big)$ is a multiplicative subset of $R$ and $Q_c(R)=\bigg( R\backslash \Big( \bigcup_{\gp \in \min (R)}\gp\Big)\bigg)^{-1}R$.
\item $\gc_R=\{ r\in R\, | \, cr=0$ for some $c\in \CC\}\subseteq \gn_R$. 
\item  $\gn_{Q_c(R)}=\CC^{-1}\gn_R$.
\end{enumerate}
\end{theorem}

{\it Proof}. 1. Let $S\in \mrLR$. By Lemma \ref{a29Dec22}.(3), the monoid $S_{mon}$ is a multiplicative subset of $R$ and  $\RSm = S_{mon}^{-1}R$. Choose  a prime ideal $P\in \Spec (S_{mon}^{-1}R)$. Then $\gp' := \s_S^{-1}(P)\in \Spec (R)$ where $\s_S: R\ra S_{mon}^{-1}R$, $r\mapsto \frac{r}{1}$, and so $$S\subseteq S_{\gp'}=R\backslash \gp'\subseteq S_\gp$$ for each $\gp\in \min (R)$ such that $\gp\subseteq \gp'$, and statement 1 follows.

2. Statement 2 follows from statement 1.

3. Clearly, $\ass_R(S_\gp )=\{ r\in R\, | \, sr=0$ for some $s\in S_\gp\}$. If $r\in  \ass_R(S_\gp )$ then $sr=0$ for some element $s\in S_\gp = R\backslash \gp$, and so   $r\in \gp$ (since $s\not\in \gp$ and the   ideal $\gp$ is a prime ideal). Hence,
$\ass_R(S_\gp )\subseteq \gp$. 

4. By statements 1 and  3, $\mrL\rad (R)=\bigcap_{\gp \in \min (R)}\ass_R(S_\gp)\subseteq \bigcap_{\gp \in \min (R)}\gp=\gn_R$.

5.  By statement 1 and Lemma \ref{a29Dec22}.(1), the set $$\CC =\bigcap_{\gp\in\min (R)} S_\gp= R\backslash \Big( \bigcup_{\gp \in \min (R)}\gp\Big)\subseteq  R\backslash \Big( \bigcap_{\gp \in \min (R)}\gp\Big)=R\backslash \gn_R=\CL\mrLR$$ is a localizable set. Hence,  the set $\CC$  is a multiplicative subset of $R$.

6. By statement 5, $\gc_R=\{ r\in R\, | \, cr=0$ for some $c\in \CC\}$. Since $c\not\in \gp$ for all $\gp \in \min (R)$, we must have $r\in \bigcap_{\gp \in \min (R)}\gp = \gn_R$, and so $\gc_R\subseteq \gn_R$.

7. An element $c^{-1}r\in Q_c(R)$, where $c\in \CC$ and $r\in R$,  belongs to the prime radical $\gn_{ Q_c(R)}$ iff it is a nilpotent element iff $\frac{r}{1}$ is a nilpotent element of $Q_c(R)$ iff $c'r^n=0$ for some element $c'\in \CC$ and a natural number $n\geq 1$ iff $(c'r)^n=0$ iff $c'r\in \gn_R$ iff $c^{-1}r\in \gn_{Q_c(R)}$.  $\Box$

For a commutative ring $R$ such that  $|\min (R)|<\infty$ and  $\gn_R$ is a nilpotent ideal (eg, $R$ is a commutative Noetherian ring), Proposition \ref{XC29Dec22} describes  $\Spec (Q_c(R))$ and  the rings  $Q_c (R)$ and $Q_c(R)/\gn_{Q_c(R)}$.

\begin{proposition}\label{XC29Dec22}%\marginpar{XC29Dec22}

Let  $R$ be a commutative ring  and  $\CC := \CC \mrLR$. Then 
\begin{enumerate}
\item  Suppose that  $|\min (R)|<\infty$. Then 
 $K\dim (Q_c(R))=0$ and $\Spec (Q_c(R))=\{ \CC^{-1}\gp \, | \, \gp \in \min (R)\}=\min  (Q_c(R))$.
\item  Suppose that  $|\min (R)|<\infty$ and  $\gn_R$ is a nilpotent ideal, eg, $R$ is a commutative Noetherian ring. Then 
\begin{enumerate}
\item $Q_c (R)\simeq \prod_{\gp \in \min (R)} R_{\gp}$ is a direct product of  local rings  $(R_{\gp}, \gp_\gp )$ and the maximal ideal $ \gp_\gp$ of $R_\gp$ is a nilpotent ideal.
\item The rings $\CC^{-1}(R/\gp)\simeq (R/\gp)_\gp$ are fields for all $\gp \in \min (R)$.
\item $Q_c(R)/\gn_{Q_c(R)}\simeq \prod_{\gp \min (R)}(R/\gp)_\gp$ is a direct product of fields.
\end{enumerate}

\end{enumerate}
\end{proposition}

{\it Proof}. 1. By Theorem \ref{C29Dec22}.(5), $\CC = R\backslash \Big( \bigcup_{\gp \in \min (R)}\gp\Big)$ is a multiplicative subset of $R$ and $$Q_c(R)=\CC^{-1}R=\bigg( R\backslash \Big( \bigcup_{\gp \in \min (R)}\gp\Big)\bigg)^{-1}R.$$ Since $|\min (R)|<\infty$,  $K\dim (Q_c(R))=0$ and $\Spec (Q_c(R))=\{ \CC^{-1}\gp \, | \, \gp \in \min (R)\}=\min  (Q_c(R))$.

2. By statement 1, the set $\Spec (Q_c(R))=\{ \CC^{-1}\gp \, | \, \gp \in \min (R)\}$ consists precisely  of maximal ideals of the ring $Q_c(R)$. Hence, they are co-prime ideals, that is $\CC^{-1}\gp+\CC^{-1}\gp'=Q_c(R)$ for all $\gp \neq \gp'$. Hence,
$$Q_c(R)/\gn_{Q_c(R)}=Q_c(R)/\bigcap_{\gp \min (R)}\CC^{-1}\gp\simeq \prod_{\gp \min (R)}Q_c(R)/\CC^{-1}\gp =\prod_{\gp \min (R)}\CC^{-1}R/\CC^{-1}\gp=\prod_{\gp \min (R)}\CC^{-1}(R/\gp)$$
is a direct product of fields. Therefore, $\CC^{-1}(R/\gp)\simeq (R/\gp)_\gp$ for all $\gp \in \min (R)$. Since $\gn_R$ is a nilpotent ideal, the orthogonal primitive idempotents that correspond to the direct product decomposition of the factor ring $Q_c(R)/\gn_{Q_c(R)}$ above can be lifted, i.e. the ring $Q_c(R) $ 
is a product of local rings. Hence,  
$$Q_c(R)\simeq \prod_{\gp \in \min (R)} R_\gp.\;\; \Box$$

For a commutative ring $R$ such that  $|\min (R)|<\infty$ and  $\gn_R=0$  (eg, $R$ is a commutative semiprime Noetherian ring), Proposition \ref{XC29Dec22} shows that  $\CC\mrLR = \CC_R$ and $Q_c(R)=Q_{cl}(R)$.

\begin{theorem}\label{YC29Dec22}%\marginpar{YC29Dec22}

Let  $R$ be a commutative ring  with  $|\min (R)|<\infty$ and  $\gn_R=0$, eg, $R$ is a commutative semiprime Noetherian ring, and  $\CC := \CC \mrLR$. Then 
\begin{enumerate}
\item The ring $R$ is a semiprime Goldie ring.
\item $Q_{cl}(R)=\prod_{\gp \in \min (R)}R_\gp$ is a product of fields.
\item  For all $\gp \in \min (R)$, $\ass_R(S_\gp)=\gp$.
\item $\CC = \CC_R$ and $Q_c(R)=Q_{cl}(R)$.
\item $\gn_{Q_c(R)}=0$.
\item $\mrL\rad (R)=0$.
\item $\gc_R=0$.
\end{enumerate}
\end{theorem}

{\it Proof}. 1. By the assumption, the ring $R$ is a commutative ring  with  $|\min (R)|<\infty$ and  $\gn_R=0$. By \cite[Theorem 5.1]{Bav-Crit-S-Simp-lQuot}, the ring $R$ is a semiprime Goldie ring since 
for each minimal prime ideal $\gp \in \min (R)$, the field $R/\gp$ is  a Goldie ring (any field is a Goldie ring). 

2 and 3. Statements 2 and follow from statement 1 and \cite[Theorem 4.1]{Bav-Crit-S-Simp-lQuot}.

4. By statements 2 and 3, the union $\bigcup_{\gp \min (R)}\gp$ is the set of all zero divisors of the ring $R$, and so $$\CC_R=R\backslash \bigcup_{\gp \min (R)}\gp =\CC,$$
by Theorem \ref{C29Dec22}.(5). Hence, $Q_c(R)=Q_{cl}(R)$.

5.  By Theorem \ref{C29Dec22}.(7) (or by statements 2 and 4), $\gn_{Q_c(R)}=\CC^{-1}\gn_R=0$.

6. By Theorem \ref{C29Dec22}.(4), $\mrL\rad (R)=0$.

7. Since $\CC=\CC_R$ (statement 4), $\gc_R=0$. $\Box$

Corollary \ref{a1Jan23} gives a characterization of commutative semiprime Goldie rings.

\begin{corollary}\label{a1Jan23}%\marginpar{a1Jan23}
Let $R$ be a  commutative semiprime ring. Then the ring $R$ is a Goldie ring iff $|\min (R)|<\infty$.
\end{corollary}

{\it Proof}. Suppose that the ring $R$ is a Goldie ring. Then the ring $Q_{cl}(R)$ is a commutative semisimple   Artinian   ring, i.e. a finite product of fields. Hence, $|\min (R)|<\infty$.

Suppose that the ring $R$ has only finitely many minimal primes. Then, by Theorem \ref{YC29Dec22}.(1),  the ring $R$ is a Goldie ring. $\Box $

% *** another proof*** Then, by Proposition \ref{XC29Dec22}.(3a,c), $Q_c(R)=Q_{cl}(R)$ is a  semisimple   Artinian  commutative ring. By Goldie's Theorem, the ring $R$ is a Goldie ring. $\Box $

%%%%%%%%%%%%%%%%%   Section  5   %%%%%%%%%%%%%%%%%%%

\section{Localizations of semiprime left Goldie rings and direct products of division rings}\label{LEFTGOLDIE}%\marginpar{LEFTGOLDIE}
 Proposition \ref{A26Dec22}.(4) describes the maximal localizable sets of finite direct product of rings and  their localizations. Proposition \ref{A26Dec22}.(1) describes localizations of the direct product of rings via localizations of its components. For a direct product of simple rings $A=\prod_{i=1}^sA_i$ such that $\CL\mrL (A_i)=A_i^\times$ for $i=1, \ldots , s$,  Theorem  \ref{b26Dec22} describes the sets $\max \mrL (A)$, $\ass_A(S)$ for all $S\in \max \mrL (A)$, the sets of localizable, non-localizable, and completely localizable elements of $A$. It also shows that $\mrL\rad (A)=0$. Every semisimple Artinian ring satisfies the assumptions of Theorem  \ref{b26Dec22}, see Corollary \ref{26Dec22} for detail. For a semiprime left Goldie ring, Lemma \ref{a27Dec22} describes all the  maximal localizable sets that contain the set of regular elements of the ring.
 
  Let $D=\prod_{i\in I}D_i$ be  a direct product of  division rings $D_i$ where $I$ is an arbitrary set. Proposition \ref{B28Dec22} describes all the localizable sets in the ring $D$ and all the localizations of $D$. Theorem \ref{A28Dec22} describes  the following sets: $\max \mrL (D)$, $\max \mrLoc (D)$, $\mrL\rad (D)$, $\CC\mrL (D)$, $Q_c(D)$, $\CL\mrL (D)$, and 
 $\CN\CL\mrL (D)$. Theorem \ref{A28Dec22}.(2) shows that for every $\CS\in \max \mrL (D)$, the ring $D\langle \CS^{-1}\rangle$ is a division ring.  \\

{\bf The maximal localizable sets of finite direct product of rings.} Let $A=\prod_{i=1}^s A_i$ be a finite direct product of rings $A_i$. For each $i=1, \ldots , n$, let $p_i :A\ra A_i$, $(a_1, \ldots , a_n)\mapsto a_i$ be the canonical projection and $\nu_i :A_i\ra A$, $a_i\mapsto (0, \ldots , 0, a_i, 0, \ldots , 0)$ is the canonical injection.

Proposition \ref{A26Dec22}.(4) describes the maximal localizable sets of finite direct product of rings and  their localizations.

\begin{proposition}\label{A26Dec22}%\marginpar{A26Dec22}
Let  $A=\prod_{i=1}^s A_i$ be a  direct product of rings, $S\subseteq A$, and $S_i:=p_i(S)$ for $i=1, \ldots , n$. Then 
\begin{enumerate}
\item $A\langle S^{-1}\rangle \simeq \prod_{i=1}^s A_i\langle S_i^{-1}\rangle$,   an $A$-isomorphism. 
\item $\ass_A(S)=\prod_{i=1}^s \ass_{A_i}(S_i)$. 
\item $S\in \mrL (A)$ iff $S_i \in \mrL (A_i)$ for some $i$.
\item $\max \mrL (A)=\coprod_{i=1}^s\max \mrL (A_i)$, where the set $\max \mrL (A_i)$ is identified with its image in $\max \mrL (A)$ via the  injection $\phi_i: \max \mrL (A_i)\ra \max \mrL (A)$, $S_i\mapsto A_1\times \cdots \times A_{i-1}\times  S_i\times A_{i+1}\times \cdots \times A_s$, and $A\langle S_i^{-1}\rangle\simeq A_i\langle S_i^{-1}\rangle$.

\end{enumerate}
\end{proposition}

{\it Proof}. 1. The ring $\prod_{i=1}^s A_i\langle S_i^{-1}\rangle$ satisfies the universal property of localization for the set $S$, and statement 1 follows.

2 and 3. Statements 2 and 3 follow from statement 1. 

4. By statement 1, $A\langle S_i^{-1}\rangle\simeq A_i\langle S_i^{-1}\rangle$ for all $S_i\in \max \mrL (A_i)$. Therefore, the map $\phi_i$ is an injection for every $i=1, \ldots , s$. By the definition of the maps $\phi_i$, 
$$\max \mrL (A)\supseteq \coprod_{i=1}^s\max \mrL (A_i).$$
Let $S\in \max \mrL (A)$. By statement 1, $A\langle S^{-1}\rangle \simeq \prod_{i=1}^s A_i\langle S_i^{-1}\rangle$,   where $S_i:=p_i(S)$ for $i=1, \ldots , s$. There is an index $i$ such that $A_i\langle S_i^{-1}\rangle\neq \{ 0\}$. Therefore, $$S\subseteq T_i:=A_1\times \cdots \times A_{i-1}\times  S_i\times A_{i+1}\times \cdots \times A_s,$$ and  
$A\langle T_i^{-1}\rangle \simeq A_i\langle S_i^{-1}\rangle$. Hence, $$\max \mrL (A)\subseteq \coprod_{i=1}^s\max \mrL (A_i),$$
and statement 4 follows. $\Box $

\begin{lemma}\label{a26Dec22}%\marginpar{a26Dec22}
Let  $A$ be a ring such that $\CL\mrL (A)=A^\times$. Then 
\begin{enumerate}
\item $\max \mrL (A)=\{ A^\times\}$ and $A\langle (A^\times)^{-1} \rangle = A$.
\item $\mrL\rad (A)=\ass_A(A^\times)=0$.
\item $\CC\mrL (A)=\CL\mrL (A)=A^\times$, and $Q_c(A)=A$, $\gc_A=0$.
\end{enumerate}
\end{lemma}

{\it Proof}. 1. By Proposition \ref{A7Dec22}.(1), $A^\times \in \mrL (A)$  and  $A\langle {A^\times}^{-1}\rangle =A$, and so  $\ass_A(A^\times)=0$. Hence, $A^\times \in \mrL (A)$. If $S\in \mrL (A)$ then $S\subseteq  \CL\mrL (A)=A^\times$, by Corollary \ref{a25Dec22}.(1). Therefore, $\max \mrL (A)=\{ A^\times\}$ and statement 1 follows.

2. By statement 1, $\max \mrL (A)=\{ A^\times\}$,  and so $\mrL\rad (A)=\ass_A(A^\times)=0$.  

3. Statement 3 follows from statement 1. $\Box$

{\it Example.} Let $(A,\gm)$ be a local ring where $\gm$ is a   maximal ideal  of $A$ which is assumed to be a {\em nil ideal} (that is every element of $\gm$ is a nilpotent element). Then $$\CL\mrL (A)=A^\times\;\; {\rm  and }\;\;\CN\CL\mrL (A)=\gm.$$ In particular, the ring $A$ satisfies the conditions of Lemma \ref{a26Dec22}, and so $\max \mrL (A)=\{ A^\times\}$, $A\langle (A^\times)^{-1} \rangle = A$, and  $\mrL\rad (A)=\ass_A(A^\times)=0$. 

Let  $A=\prod_{i=1}^s A_i$ be a  direct product of simple rings $A_i$. Then $\min (A)=\{ \gp_i\, | \, i=1, \ldots , s\}$ where $\gp_i:= A_1\times \cdots \times A_{i-1}\times  \{0\}\times A_{i+1}\times \cdots \times A_s$, and $A/\gp_i\simeq A_i$.

 For the  ring $A=\prod_{i=1}^s A_i$ such that  $\CL\mrL (A_i)=A_i^\times$  for $i=1, \ldots , n$,  Theorem \ref{b26Dec22} describes in an explicit way the following sets: $\max \mrL (A)$, $\mrL\rad (A)$, $\CC\mrL (A)$, $Q_c(A)$, $\CL\mrL (A)$, and 
 $\CN\CL\mrL (A)$.

\begin{theorem}\label{b26Dec22}%\marginpar{b26Dec22}
Let  $A=\prod_{i=1}^s A_i$ be a  direct product of  rings $A_i$ such that $\CL\mrL (A_i)=A_i^\times$  for $i=1, \ldots , n$. Then 
\begin{enumerate}
\item $\max \mrL (A)=\{S_i \, | \, i=1, \ldots , s\}$ where 
$S_i :=p_i^{-1}(A_i^\times)=\{ (a_1,\ldots , a_s) \in A\, | \, a_i\in A_i^\times, a_j\in A_j$ for all $j\neq i\}$.
\item $\ass_A (S_i )=\gp_i$ and $A\langle S_i^{-1} \rangle\simeq A/\gp_i = A_i$.
\item $\mrL\rad (A)=0$.
\item $\CC\mrL (A)=A^\times$ and $Q_c(A)=A$.
\item  $\CL\mrL (A)=\{ (a_1, \ldots , a_s)) \in A \, | \,  a_i \in A_i^\times$ for some $i\}$ and $\CN\CL\mrL (A)=\prod_{i=1}^sA_i\backslash A_i^\times$.
\end{enumerate}
\end{theorem}

{\it Proof}. 1. Statement 1 follows from Proposition \ref{A26Dec22}.(4) and Lemma \ref{a26Dec22}.(1).

2. Statement 2 follows from statement 1. 

3. By statement 2, $\mrL\rad (A)=\bigcap_{i=1}^s\gp_i=0$.

4. By statement 1,  $\CC\mrL (A)=\bigcap_{i=1}^sS_i=A^\times$, and so  $Q_c(A)=A\langle \CC\mrL (A)^{-1}\rangle=A$. 

5.  Statement 5 follows from statement 1. $\Box$

{\it Example.} $A=\prod_{i=1}^s A_i$ be a  direct product of local  rings $(A_i, \gm_i)$ where the maximal ideals $\gm_i$ are nil ideals. Then the ring $A$ satisfies the assumption of 
 Theorem \ref{b26Dec22} (see the Example above).

 Let $A$ be a semisimple Artinian ring. Then the set of minimal primes $\min (A)$ of  $A$ is a finite set, $A(\gp) :=A/\gp$ is a simple Artinian ring for every $\gp \in \min (A)$, and
$$A\simeq \prod_{\gp \in \min (A)} A(\gp),$$
 and vice versa. We identify the rings $A$  and $\prod_{\gp \in \min (A)} A(\gp)$. Then each element $a\in A$ is identified with $(a_\gp )_{\gp \in \min (A)}\in \prod_{\gp \in \min (A)} A(\gp)$, where $a_\gp \in A(\gp )$,  and 
  for each $\gp \in \min (A)$, $$\gp = \prod_{\gq \in \min (A)\backslash \{ \gp\} } A(\gq),\;\; A=A (\gp )\oplus \gp, \;\;  \pi_\gp : A\ra A(\gp), \;\; a\mapsto a_\gp := a+\gp.$$
  For the  semisimple Artinian ring $A$, Corollary \ref{26Dec22} describes  the following sets: $\max \mrL (A)$, $\mrL\rad (A)$, $\CC\mrL (A)$, $Q_c(A)$, $\CL\mrL (A)$, and 
 $\CN\CL\mrL (A)$. 
  
 \begin{corollary}\label{26Dec22}%\marginpar{26Dec22}
Let $A$ be a semisimple Artinian ring. We keep the notation as above. Then
\begin{enumerate}
\item $\max \mrL (A)=\{\CS_\gp \, | \, \gp \in \min (A)\}$ where 
$\CS_\gp :=\pi_\gp^{-1}(A(\gp)^\times)=\{ (a_\gq)_{\gq \in \min (A)}\in A\, | \, a(\gp)\in A(\gp )^\times\}$.
\item $\ass_A (\CS_\gp )=\gp$ and $A\langle \CS_\gp^{-1} \rangle\simeq A/\gp = A(\gp)$.
\item $\mrL\rad (A)=0$.
\item $\CC\mrL (A)=A^\times=\CC_A$ and $Q_c(A)=A$.
\item  $\CL\mrL (A)=\{ (a_\gp )_{\gp \in \min (A)}\in \prod_{\gp \in \min (A)} A(\gp)\, | \,  a_\gp \in A(\gp )^\times$ for some $\gp \in \min (A)\}$ and $\CN\CL\mrL (A)=\prod_{\gp \in \min (A)}A(\gp )\backslash A(\gp )^\times$.
\end{enumerate}
\end{corollary}

{\it Proof}. For each simple Artininian ring $\CA$, $\CC_\CA= {}'\CC_\CA =\CC_\CA'=\CA^\times$. Hence, $\ass_\CA (a)=\CA$ for all $a\in \CA\backslash \CA^\times$. Therefore, $\CL\mrL (\CA) = \CA^\times$, and so the condition of Theorem  \ref{b26Dec22} hold, and the corollary follows from  Theorem  \ref{b26Dec22}. $\Box$\\

{\bf Classification of maximal left localizable sets of a semiprime left Goldie ring that contain the set of regular elements of the ring.}  It was proved that the set of maximal left denominators sets of $R$, $\maxDen_l(R)$,  is a {\em finite set}  if the classical left quotient ring $Q_{ l, cl}(R):=\CC_R^{-1}R$ of $R$ is
 a semisimple  Artinian ring, \cite{Bav-Crit-S-Simp-lQuot},  or a left Artinian ring,  \cite{Bav-LocArtRing}, or a left Noetherian ring, \cite{Crit-lNoeth-lQuot}. In each of the three cases an explicit description of the set $\maxDen_l(R)$ is given. 
 
 % ***  For a ring $R$, the rings $$Q_{l, cl}(R):=\CC_R^{-1}R\;\; {\rm  and}\;\; Q_{r, cl}(R):=R\CC_R^{-1}$$ are called the {\em classical left and right quotient rings} provided they exist, respectively. If both rings exist then they are isomorphic and the ring
% $$Q_{cl}(R):= Q_{l, cl}(R)\simeq  Q_{r, cl}(R)$$
 %is called the {\em classical quotient ring} of $R$.  Goldie's Theorem states that the ring $Q_{l, cl}(R)$ is  a semisimple Artinian ring iff the ring $R$ is a semiprime left Goldie ring (the ring $R$ is called a {\em left Goldie ring} if   $\udim (R)<\infty$ and the ring $R$ satisfies the a.c.c. on left annihilators where $\udim$ stands for the {\em left uniform dimension}). In a similar way a {\em right Goldie ring} is defined. A left and right Goldie ring is called a {\em Goldie ring}.    ***

 Theorem \ref{S18Jan19} is a classification  of maximal left localizable sets of a semiprime left Goldie ring that contain the set of regular elements of the ring.

\begin{theorem}\label{S18Jan19}%\marginpar{S18Jan19}
Let $R$ be a semiprime left Goldie ring.   Then $\{ \CS\in \max  \mrLR \, | \, \CC_R\subseteq \CS\} = \max  \Den_l(R)=\{ \CC (\gp ) \, | \, \gp \in \min (R)\}$ where $\CC (\gp ):= \{ c\in R \, | \, c+\gp \in \CC_{R/\gp} \}$. 
\end{theorem}

{\it Proof}. (i) $\max  \Den_l(R)=\{ \CC (\gp ) \, | \, \gp \in \min (R)\}$: See \cite[Theorem 3.1.(1-4)]{Bav-Crit-S-Simp-lQuot}.

(ii)  $\{ \CC (\gp ) \, | \, \gp \in \min (R)\}\subseteq \max \mrL (R)$: The inclusion follows from Corollary  \ref{26Dec22}.(1). 

Suppose that $S\in \mrLR$  and $\CC (\gp) \subseteq S$. We have to show that $S\subseteq \CC (\gp)$. By \cite[Theorem 3.1.(3)]{Bav-Crit-S-Simp-lQuot}, $\CA :=\CC (\gp)^{-1}R$ is a simple Artinian ring. In particular, $\CL\mrL(\CA)=\CA^\times$. Let $\tau: R\mapsto \CA$, $r\mapsto \frac{r}{1}$. By Proposition \ref{A7Dec22}.(2),
$$\RSm \simeq R\langle \CC (\gp )^{-1}\rangle \langle \tau (S)^{-1}\rangle=\Big( \CC (\gp)^{-1}R\Big)\langle \tau (S)^{-1}\rangle=\CA\langle \tau (S)^{-1}\rangle=\CA ,$$
since $\tau (S)\subseteq \CL\mrL(\CA)=\CA^\times$. It follows from the definition of the set $\CC (\gp)$ that 
$$ S\subseteq \tau^{-1} (\tau (S))\subseteq \tau^{-1} (\CA^\times) =\CC (\gp).$$

(iii)  $\{ \CS\in \max  \mrLR \, | \, \CC_R\subseteq \CS\} \subseteq \{ \CC (\gp ) \, | \, \gp \in \min (R)\}$: Let $\CS\in \max  \mrLR $ be such that $\CC_R\subseteq \CS$. By Proposition \ref{A7Dec22}.(2),
$$R\langle \CS^{-1}\rangle \simeq R\langle \CC_R^{-1}\rangle \langle \s (\CS)^{-1}\rangle= Q_{l,cl}(R)\langle \s (\CS)^{-1}\rangle=\bigg(\prod_{\gp \in \min (R)} Q(\gp)\bigg)\langle \s (\CS)^{-1}\rangle$$
 where $\s: R\ra Q_{l,cl}(R)=\prod_{\gp \in \min (R)} Q(\gp)$ is a direct product of simple Artinian rings $Q(\gp )$. By Corollary \ref{26Dec22}.(1), $\s (\CS )\subseteq \CS_\gp$ for some $\gp \in \min (R)$. By Corollary \ref{26Dec22}.(1),  $$\CS \subseteq \s^{-1} (\CS_\gp )=\CC (\gp ),$$ and so  we must have $\CS = \CC (\gp)$, by the statement (ii).
 
 Now, the theorem follows from the satements (i)--(iii). $\Box $

So, every  maximal left localizable set of a semiprime left Goldie ring that contains the set of regular elements of the ring is a maximal left denominator set, and vice versa.

\begin{corollary}\label{B30Dec22}%\marginpar{B30Dec22}
Let $R$ be a semiprime left Goldie ring which is not a prime ring.   Then $Q_a(R)=\{0\}$. 
\end{corollary}

{\it Proof}. By Theorem \ref{S18Jan19}, $|\max \mrLR |\geq 2$, and the corollary follows from Theorem \ref{XA30Dec22}. $\Box$\\

{\bf Classification of maximal Ore and denominator sets of a semiprime  Goldie ring.}   The sets of localizable left, right or two-sided Ore sets   of $R$ are denoted by $\mLlR$, $\mLrR$ and $\mLR$, respectively. Clearly, $ \mLR = \mLlR \cap \mLrR$.  In order to work with these three sets simultaneously we use the following notation $\mLsR$  where $*\in \{ l,r, \emptyset \}$  and $\emptyset$ 
is the empty set   $(\mL (R) = \mL_\emptyset (R))$. Then $\mL_* (R)=\mrL (R)\cap \Ore_* (R)$ for  $*\in \{ l,r, \emptyset \}$. Ideals of a ring are called {\em incomparable} if none of them is contained in the other. 
 For a ring $R$, $\min (R)$ is the set of its minimal prime ideals.

 The next theorem is an explicit description of maximal Ore sets of a semiprime Goldie ring.
 
\begin{theorem}\label{17Jan19}%\marginpar{17Jan19}
(\cite[Theorem 1.11]{LocSets}) Let $R$ be a semiprime Goldie ring and $\CN_* :=\{ S\in \max  \mL_*(R) \, | \,$ $ \CC_R\subseteq S\}$ where $*\in \{ l,r,\emptyset\}$. Then 
\begin{enumerate}
\item $\max  \Ore (R)= \max \Den (R) = \{ \CC (\gp ) \, | \, \gp \in\min (R)\}=\CN_*$ for all $*\in \{ l,r,\emptyset\}$  where $\CC (\gp ):= \{ c\in R \, | \, c+\gp \in \CC_{R/\gp} \}$. So, every maximal Ore set of $R$ is a maximal denominator set, and vice versa. 
\item For all $S\in \max  \Ore (R)$ the ring $S^{-1}R$ is a simple Artinian ring. 
\item $Q_{cl} (R) \simeq \prod_{S\in \max  \Ore (R)} S^{-1}R$. %(where $Q_{cl} (R) := \CC_R^{-1}R\simeq  R\CC_R^{-1}$ is the classical quotient ring of $R$).
\item $\max  \ass\, \Ore (R) = \ass \, \max  \Ore (R)=\min (R)$ where $\ass \, \max  \Ore (R):=\{ \ass_R(S)\, | \, S\in \max  \Ore (R)\}$. In particular, the ideals in the set $\ass \, \max \, \Ore (R)$ are incomparable. 
\end{enumerate}
\end{theorem}

\begin{corollary}\label{a27Dec22}%\marginpar{a27Dec22}
Let $R$ be a semiprime Goldie ring.  Then $\{ \CS\in \max  \mrLR \, | \, \CC_R\subseteq \CS\}=\max  \Ore (R)= \max \Den (R) = \{ \CC (\gp ) \, | \, \gp \in\min (R)\}=\CN_*$ for all $*\in \{ l,r,\emptyset\}$. 
\end{corollary}

{\it Proof}. The corollary follows from Theorem \ref{S18Jan19} and Theorem \ref{17Jan19}.(1). $\Box $ \\

{\bf Maximal localizable sets of a direct product of division rings.} Let $I$ be a  set and $\CP (I)$ be its {\em power set} of $I$, i.e. the set of all subsets of $I$ (an element of $\CP (I)$ is a  subset of $I$). An element $F\in \CP (I)$ is called a {\em filter} on $I$ if the following conditions hold: $\emptyset \not\in F$; if $\ga \in F$ then $\gb \in F$ for all $\gb\subseteq I$ such that $\ga \subseteq \gb$; and 
for all elements $\ga, \gb \in F$, $\ga \cap \gb \in F$. The set of all filters on $I$ is denoted by $\CF (I)$. A maximal element of $\CF (I)$ is called an {\em ultrafilter} on $I$. The set of all ultrafilters on $I$ is denoted by $\CU                      (I)$. For each element $i\in I$, the set $$F_i:=\{ \ga \subseteq I \, | \, i\in \ga\}\in \CU (I)$$ is called a {\em principal ultrafilter}. An ultrafilter is principal iff it contains a finite set. If the set $I$ is a finite set then all ultrafilters are principal. If $I$ is an infinite set then an ultrafilter is non-principal iff it contains the {\em Frechet filter} of co-finite subsets (a subset $\ga$  of $I$ is called {\em co-finite} if its {\em complement} $C\ga :=I\backslash \ga$ is a finite set). A filter $F\in \CF (I)$ is an ultrafilter iff for each subset $\ga $ of $I$ either $\ga \in F$ or $C\ga \in F$.  The set $(\CF (I),  \subseteq )$ is a poset w.r.t. inclusion where $\CF\subseteq \CF'$ if for all elements $\gf\in \CF$, $\gf\in \CF'$.

Let $D=\prod_{i\in I}D_i$ be  a direct product of  division rings $D_i$. For each subset $J$ of $I$, let $D_J:=\prod_{j\in J}D_j$. We identify the set $D_J$ with the ideal $D_J\times \prod_{k\in CJ}\{0\}$ of $D$ where $CJ:=I\backslash J$ is the complement of $J$ in $I$. If $J\subseteq J'\subseteq I$ then $D_J\subseteq D_{J'}$, and $D_J=D_{J''}$ for some subset $J''$ of $I$ iff $J=J''$. For each element $d=(d_i)_{i\in I}\in D$, where $d_i\in D_i$, the set $\supp (d):=\{ i\in I\, | \, d_i\neq 0\}$ is called the {\em support} of $d$. For all elements $d,e\in D$, 
$$\supp (de)=\supp (d)\cap \supp (e)\;\; {\rm and}\;\; \supp (d\pm e)\subseteq \supp (d)\cup \supp (e).$$
Each element of the ring $D$ is either a unit or a zero divisor. An element $d\in D$ is a unit (resp., a zero divisor) iff $\supp (d)=I$ (resp., $\supp (d) \neq I$).

Proposition \ref{B28Dec22} describes all the localizable sets in the ring $D$ and all the localizations of $D$.

\begin{proposition}\label{B28Dec22}%\marginpar{B28Dec22}
Let $D=\prod_{i\in I}D_i$ be  a direct product of  division rings $D_i$. 
\begin{enumerate}
\item Suppose that $S$ is a multiplicative submonoid of $D$. Then the following statements are equivalent:
\begin{enumerate}
\item $S\in \mrL (D)$.
\item $S$ is a multiplicative subset of $D$.
\item $S\in \Den (D)$.
\end{enumerate}
\item Suppose that $S$ is a multiplicative subset of $D$. Then $D\langle S^{-1}\rangle \simeq S^{-1}D\simeq DS^{-1}\simeq D/\ass_D(S)$ where $\ass_D(S)=\bigcup_{s\in S}D_{C\supp (s)}=\sum_{s\in S}D_{C\supp (s)}$.

\end{enumerate}
\end{proposition}

{\it Proof}. 1.  Straightforward.
 
2. By statement 1, $S\in \Den (D)$, and so $D\langle S^{-1}\rangle \simeq S^{-1}D\simeq DS^{-1}$ and  $$\ass_D(S)=\bigcup_{s\in S}D_{C\supp (s)}=\sum_{s\in S}D_{C\supp (s)}.$$ Since $D/\ass_D(S)\subseteq S^{-1}D$ and for every element $s\in S$, the element  $s+\ass_D(S)$ is a unit of the factor ring $D/\ass_D(S)$, we must have $D/\ass_D(S)= S^{-1}D$. $\Box$

By Proposition \ref{B28Dec22}.(1), every localizable set of the ring $D$ is a denominator set, and vice versa. By Proposition \ref{B28Dec22}.(2), for all nonzero elements $d\in D$, $R\langle d^{-1}\rangle\simeq R/D_{C\supp (s)}$.

Let $S$ be a multiplicative subset of $D$, then the set
%\marginpar{Ssat}
\begin{equation}\label{Ssat}
S_{sat}:=\bigcup_{s\in S} \bigg(\prod_{i\in \supp (s)}D_i^\times \times \prod_{j\in C\supp (s)}D_j\bigg)
\end{equation}
is called the {\em saturation} of $S$. Clearly, $(S_{sat})_{sat}=S_{sat}$. For each subset $S$ of $D$, let $\supp (S):=\{ \supp (s)\, | \, s\in S\}\in \CP (I)$. For each filter $\CF \in \CF (I)$, the set
%\marginpar{Ssat1}
\begin{equation}\label{Ssat1}
\CS (\CF):=\{ s\in D\, | \, \supp (s) \in \CF \}
\end{equation}
is a multiplicative subset of $D$.  For all filters $\CF, \CG\in \CF (I)$ such that $\CF \subseteq \CG$, $\CS (\CF)\subseteq \CS (\CG )$. A  subset $S$ of $D$ is called a {\em saturated subset} if $S=S_{sat}$.    Let $\CP (D)_{sat}$ be the set of all saturated subsets of $D$. The sets $(\CF (I), \subseteq )$ and $(\CP (D)_{sat}, \subseteq )$  are posets.

Proposition \ref{C28Dec22}.(4) shows that the  posets  $\CF (I)$ and $\Sub(D)_{sat}$ are isomorphic. 

\begin{proposition}\label{C28Dec22}%\marginpar{C28Dec22}
Let $D=\prod_{i\in I}D_i$ be  a direct product of  division rings $D_i$, $S$ be a multiplicative subset of $D$, and $\CF\in \CF (I)$. Then
\begin{enumerate}
\item $S_{sat}$ is a multiplicative subset of $D$, $S^{-1}D\simeq S^{-1}_{sat}R$ (an $R$-isomorphism), and $\ass_R(S)=\ass_R(S_{sat})$.
\item $\supp (S_{sat})\in \CF (I)$ and $S_{sat} =\CS (\supp (S_{sat}))$.
\item $\CS (\CF )_{sat}=\CS (\CF )$ and $\supp (\CS (\CF))=\CF$.
\item The map $\CF (I)\ra \CP (D)_{sat}$, $\CF \mapsto \CS (\CF)$ is an isomorphism of posets with inverse $S\mapsto \supp (S)$. 
\end{enumerate}
\end{proposition}

{\it Proof}. 1. By the definition of the set  $S_{sat}$, it is a multiplicative subset of $D$ such that  
 $\ass_R(S)=\ass_R(S_{sat})$. Now, by Proposition \ref{A28Dec22}.(2), $$S^{-1}D\simeq D/\ass_R(S)=D/\ass_R(S_{sat})\simeq S_{sat}^{-1}D.$$
 2. Since the set $S_{sat}$ is a saturated subset of $D$,  $\supp (S_{sat})\in \CF (I)$.   The equality $S_{sat} =\CS (\supp (S_{sat}))$ follows from  the definitions of the sets $S_{sat}$ and $\CS (\supp (S_{sat}))$.  
 
 3. The equality $\CS (\CF )_{sat}=\CS (\CF )$ is obvious. The equality $\supp (\CS (\CF))=\CF$  follows from the definitions of the sets $\CS (\CF)$ and $\supp (\CS (\CF))$.
 
 4. Statement 4 follows from statements 2 and 3. $\Box$

Let $\CF \in \CF (I)$. Then $(\CF , \supseteq )$ is a directed set since for all $\gf , \gq\in \CF$, $f\supseteq \gf\cap \gq$ and $\gq \supseteq \gf\cap\gq$. For each $\gf \in \CF$, let $D_\gf =\prod_{i\in \gf}D_i$, a direct product of division rings $D_i$. For each pair of elements $\gf , \gq\in \CF$ such that $\gf\supseteq \gq$, let $p_{\gq\gf}: D_\gf\ra D_\gq$ be the natural  projection map. Then $(D_\gf, p_{\gq\gf})$ is a directed system.

Let $\CP (D)_{mult}$ be the set of all multiplicative subsets of $D$.  For the   ring $D$, Theorem \ref{A28Dec22} describes  the following sets: $\max \mrL (D)$, $\max \mrLoc (D)$, $\mrL\rad (D)$, $\CC\mrL (D)$, $Q_c(D)$, $\CL\mrL (D)$, and 
 $\CN\CL\mrL (D)$. 
 
\begin{theorem}\label{A28Dec22}%\marginpar{A28Dec22}
Let $D=\prod_{i\in I}D_i$ be  a direct product of  division rings $D_i$. Then 
\begin{enumerate}
\item $\max \mrL (D)=\max \Den_* (D)=\max \Ore_* (D)=\max \CP (D)_{mult}=\{\CS (\CF ) \, | \,\CF \in \CU (I)\}$ where $*\in \{ l,r,\emptyset\}$.
\item For each $\CS = \CS (\CF)\in \max \mrL (D)$, $$D\langle \CS^{-1} \rangle\simeq \varinjlim_{\gf\in \CF}D_\gf \simeq \CS^{-1}D\simeq D\CS^{-1}\simeq D/\ass_D (\CS )$$ is a division ring and $\ass_D (\CS )=\bigcup_{\gf\in \CF}D_{C\gf} =\sum_{\gf\in \CF}D_{C\gf}$ is a maximal ideal of the ring $D$.
\item $\mrL\rad (D)=0$.
\item $\CC\mrL (D)=\CC_D=D^\times$ and $Q_c(D)=Q_{cl}(D)=D$.
\item  $\CL\mrL (D)=D\backslash \{ 0\}$ and $\CN\CL\mrL (D)=\{ 0\}$.
\end{enumerate}
\end{theorem}

{\it Proof}. 1. By Proposition \ref{C28Dec22}.(4), $\max \mrL (D)=\{\CS (\CF ) \, | \,\CF \in \CU (I)\}$. By Proposition \ref{B28Dec22}.(1), $$\max \mrL (D)=\max \Den_* (D)=\max \Ore_* (D)=\max \CP (D)_{mult}\;\; {\rm  where}\;\; *\in \{ l,r,\emptyset\}.$$

2. By Proposition \ref{B28Dec22}.(2), $$D\langle \CS^{-1} \rangle\simeq \CS^{-1}D\simeq D\CS^{-1}\simeq D/\ass_D (\CS )\;\; {\rm   and}\;\;\ass_D (\CS )=\bigcup_{\gf\in \CF}D_{C\gf} =\sum_{\gf\in \CF}D_{C\gf}.$$
 
Let $s+\ass_D(S)$ be a nonzero element of the ring $D/\ass_R(S)$. Then $\supp (s)\in \CF$ and the element $s+D_{C\supp (s) }\in D/D_{C\supp (s)}$ is a unit. The ring $D/\ass_D(S)$ is an epimorphic image of the ring $ D/D_{C\supp (s)}$ since $D_{C\supp (s) }\subseteq  \ass_D(S)$. Therefore the element $$s+\ass_D(S)\in D/\ass_D(S)$$ is a unit (as an image of the unit $s+D_{C\supp (s) }$). Therefore, the ring $D/\ass_D(S)$ is a division ring and the ideal $\ass_D(S)$ is a maximal ideal of $D$.

For each element $\gf \in \CF$, let $$\phi_\gf : D_\gf\simeq D/D_{C\gf}\ra D/\ass_R(\CF)$$ be the natural epimorphism that is determined by the inclusion  $D_{C\gf}\subseteq \ass_R(S)$. For each pair of elements $\gf , \gq\in \CF$ such that $\gf\supseteq \gq$, $\phi_\gf = \phi_\gq p_{\gq \gf}$. Hence, 
there is a ring homomorphism $$\varinjlim_{\gf \in \CF}D_\gf \ra D/\ass_D(S).$$
The ideal $\ass_D(S)=\bigcup_{\gf\in \CF}D_{C\gf}$ belongs to the kernel of  the homomorphism $\phi_I: D=D_I\ra \varinjlim_{\gf \in \CF}D_\gf$. Hence,  there is a homomorphism 
$$D/\ass_D(S)\ra \varinjlim_{\gf \in \CF}D_\gf.$$
Hence, $\varinjlim_{\gf \in \CF}D_\gf \simeq D/\ass_D(S)$. 

3. For each $i\in I$, $F_i$ is a principal ultrafilter on $I$, and $\CS (F_i)=D_i^\times\times \prod_{j\in I\backslash \{ i\} }D_j$. Hence, $\ass_D(\CS (F_i))=\{0\}\times \prod_{j\in I\backslash \{ i\} }D_j$. Now, by statement 1,  $$\mrL\rad (D)=\bigcap_{\CF\in \CU (I)}\ass_D(\CS (\CF))\subseteq \bigcap_{i\in I}\ass_D(\CS (F_i))=0,$$ and so $\mrL\rad (D)=0$.

4. Let $\CC = \CC\mrL (D)$. Recall that for each $i\in I$, $F_i$ is a principal ultrafilter on $I$, and $\CS (F_i)=D_i^\times\times \prod_{j\in I\backslash \{ i\} }D_j$. Hence,
$$D^\times \subseteq \CC=\bigcap_{\CF\in \CU (I)}\CS (\CF)\subseteq \bigcap_{i\in I}\CS (F_i)=D^\times,$$ and so 
$\CC = D^\times =\CC_D$. Since $\CC=D^\times$,  $Q_c(D)=D\langle \CC^{-1}\rangle=D\langle (D^\times)^{-1}\rangle=D\langle D^\times\rangle=D=Q_{cl}(D)$.

5. Statement 5 is obvious.  $\Box $

%%%%%%%%%%%%%%%%%% SECTION 6 %%%%%%%%%%%%%%%%%%%%%%%%

\section{Localization of a module at a localizable set}\label{LOCMOD}%\marginpar{LOCMOD}

The aim of the section is to introduce the concept of localization of a module at a localizable set and to consider its basic properties. Proposition \ref{SA20Jan19}.(2) is a universal property of localization of a module. Theorem \ref{S20Jan19} is a criterion for the localization functor $M\ra S^{-1}M$ to be an exact functor.

{\it Definition.} Let $R$ be a ring, $S\in \mrLRa$,   and $M$ be a left  $R$-module. Then the left $\RSm$-module $$S^{-1}M:= \RSm\t_R M$$ is called the {\em (left) localization} of $M$ at $S$. If $S\in \mrLRa$   and $M$ be a right  $R$-module. Then the right  $\RSm$-module $$MS^{-1}:= M\t_R \RSm$$ is called the {\em (right) localization} of $M$ at $S$.  \\

 By the definition,  $S^{-1}M$ is a left  $\RSm$-module. By applying $-\t_R M$ to the $R$-homomorphism $R\ra \RSm$ we obtain the $R$-homomorphism  $i_M: M\ra S^{-1}M$. 
 Proposition \ref{SA20Jan19}.(2) is the {\em universal property} of the $\RSm$-module $S^{-1} M$. 

\begin{proposition}\label{SA20Jan19}%\marginpar{SA20Jan19}
Let $R$ be a ring, $S\in \mrLRa$, $M$ be an $R$-module, and $i_M : M\ra S^{-1}M$. Then 
\begin{enumerate}
\item $S^{-1}M\simeq \bS^{-1}(M/ \ga M)$ where $\bS := (S+\ga ) / \ga \in \mrL (\bR , 0)$ and $\bR = R/ \ga$ (Corollary \ref{a17Dec22}.(2)). 
\item Let $\CM$ be an $\RSm$-module and $f: M\ra \CM$ be an $R$-homomorphism. Then there is a unique $\RSm$-homomorphism $S^{-1}f: S^{-1}M\ra \CM$ such that $f=S^{-1}f\circ i_M$, 
$$\begin{array}[c]{ccc}
M&\stackrel{i_M}{\ra} &S^{-1}M\\
&\stackrel{f}{\searrow}&{\downarrow}^{\exists !\, S^{-1}f}\\
&  & \CM .
\end{array}$$ 
This property is uniquely characterized the $\RSm$-module $S^{-1}M$ up to isomorphism.
\end{enumerate}
\end{proposition}

{\it Proof}. 1. Let $\bM=M/\ga M$. Then $S^{-1}(\ga M) = \RSm\t_R \ga M = \RSm\ga \t_R M=0$, and so $$S^{-1}M= S^{-1}\bM = \RSm\t_R\bM\simeq \bR \langle \bS^{-1}\rangle \t_{\bR } \bM = \bS^{-1} \bM . $$

2. The $R$-homomorphism $f: M\ra \CM$  determines a  ring homomorphism $$R\ra \End_\Z (\CM ), \;\; r\mapsto (m\mapsto rm).$$ The images of the elements of the set $S$ in $\End_\Z (\CM )$ are units. Now, by  Proposition  \ref{A5Dec22}, there is a unique $\RSm$-homomorphism $S^{-1}f: S^{-1}M\ra \CM$ such that $f=S^{-1}f\circ i_M$. By the uniqueness of the map $S^{-1}f$, this property is uniquely characterized the $\RSm$-module $S^{-1}M$ up to isomorphism.  $\Box$

The $R$-homomorphism $i_M:M\ra S^{-1}M$ is equal to the compositions of $R$-homomorphisms
$$ i_M: M\ra R/\ga \t_R M =M/\ga M\ra S^{-1}M=\bS^{-1}(M/\ga M).$$
Therefore, $\gt_S(M):= \ker (i_M)\supseteq\ga M$. We have the  descending chain of $R$-modules 
$$ M\supseteq \gt_S(M)\supseteq \gt_S^2(M)\supseteq \cdots \supseteq \gt_S^n(M)\supseteq \cdots $$
where $\gt_S^n(M) = \gt_S\gt_S\cdots \gt_S (M)$, $n$ times. 
 An $R$-module $M$ is called $S$-{\em torsion} (resp., $S$-{\em torsionfree}) if $S^{-1}M=0$ (resp., $\gt_S(M)=0$, i.e., the map $i_M: M\ra S^{-1}M$, $m\mapsto 1\t m$ is an $R$-module monomorphism). Let $\gf_S (M)=\im (i_M)$, the image of the map $i_M$, and we have a short exact sequence of $R$-modules 
%\marginpar{xSlcm2}
\begin{equation}\label{xSlcm2}
0\ra \gt_S(M)\ra M\ra \gf_S(M)\ra 0. 
\end{equation}
 We have the  descending chain of $R$-module epimorphisms 
$$ M\ra \gf_S(M)\ra \gf_S^2(M)\ra \cdots \ra \gf_S^n(M)\ra \cdots $$
where $\gf_S^n(M) = \gf_S\gf_S\cdots \gf_S (M)$, $n$ times.

For a ring $R$, let $R$-{\rm mod} be the category of left $R$-modules.  By Proposition \ref{SA20Jan19}.(1), the localization at $S$ functor, $$S^{-1} : R-{\rm mod}\ra \RSm-{\rm mod}, \;\; M\mapsto S^{-1}M,$$ is the composition of two  functors 
%\marginpar{Slcm}
\begin{equation}\label{Slcm}
S^{-1} = \bS^{-1} \circ (R/\ga \t_R-)
\end{equation}
and the second one is a right exact functor.

Applying the right exact functor $-\t_RM$ to the short exact sequence $0\ra \ga \ra R\ra R/\ga \ra 0$ of right $R$-modules yields the  short exact sequence $ \ga \t_RM\ra M\ra R/\ga \t_RM\ra 0$, and so $$ R/\ga \t_RM\simeq M/\ga M,$$  an isomorphism of left $R$-modules. Now, applying the right exact functor $R/\ga \t_R-$ to the short exact sequence   of left $R$-modules, $0\ra M_1\ra M_2\ra M_3\ra 0$,  yields the  short exact sequence of left $R$-modules
 $$0\ra M_1\cap \ga M_2/\ga M_1 \ra M_1/\ga M_1\ra M_2\ga M_2\ra M_3\ga M_3\ra 0.$$

Suppose that  the functor $S^{-1}$ is also a right exact functor for some $S\in \mrLRa$. Then  the sequence of left $\RSm$-modules 
%\marginpar{xSlcm1}
\begin{equation}\label{xSlcm1}
0\ra S^{-1} (M_1\cap \ga M_2/\ga M_1)\simeq  \bS^{-1} (M_1\cap \ga M_2/\ga M_1)\ra S^{-1}M_1\ra S^{-1}M_2\ra  S^{-1}M_3\ra 0
\end{equation}
is exact. Notice that $$M_1\cap \ga M_2/\ga M_1\simeq \bM_1\cap (\ga M_2/\ga M_1).$$

Theorem \ref{S20Jan19} is a criterion for the functor $S^{-1}:M\ra S^{-1}M$ to be exact. 

\begin{theorem}\label{S20Jan19}%\marginpar{S20Jan19}
Let $R$ be a ring, $S\in \mrLRa$, $\bR = R/ \ga $ and $\bS:= (S+\ga ) / \ga $. The functor $S^{-1}$ is exact iff for all $R$-modules $M_1$ and $M_2$ such that $M_1\subseteq M_2$, the $\bR$-modules $M_1\cap \ga M_2/\ga M_1$ is $\bS$-torsion. 
\end{theorem}

{\it Proof}. The theorem follows from the exact sequence (\ref{xSlcm1}). 
 $\Box $

\begin{corollary}\label{Sb20Jan19}%\marginpar{Sb20Jan19}
Let $R$ be a ring, $S\in \mrLRa$, $\bR = R/ \ga $ and $\bS:= (S+\ga ) / \ga $. If the functor $S^{-1}$ is exact then the $\bR$-module $\ga / \ga^2$ is $\bS$-torsion or, equivalently,  $\bS$-torsion.
\end{corollary}

{\it Proof}. Applying Theorem \ref{S20Jan19} to the pair of $R$-modules $M_1=\ga \subseteq M_2= R$, we conclude that the $R$-module (resp., $\bR$-module)  $(M_1\cap \ga M_2) / \ga M_1=\ga/ \ga^2$ is $S$-torsion (resp., $\bS$-torsion). $\Box$

Since $S^{-1}\ga M=\RSm \t_R \ga M=\RSm \ga\t_R  M=0$,  $\ga M \subseteq \gt_S(M)$. By taking the short exact sequence (\ref{xSlcm2}) modulo $\ga M$, we obtain a short exact sequence of $\bR$-modules 
%\marginpar{Slcm3}
\begin{equation}\label{Slcm3}
0\ra \gt_S(M)/\ga M\ra M/\ga M\ra \gf_S(M)\ra 0. 
\end{equation}

{\bf Licence.} For the purpose of open access, the author has applied a Creative Commons Attribution (CC BY) licence to any Author Accepted Manuscript version arising from this submission.

{\bf Disclosure statement.} No potential conflict of interest was reported by the author.

{\bf Data availability statement.} Data sharing not applicable – no new data generated.

\small{

School of Mathematics and Statistics 

University of Sheffield

Hicks Building

Sheffield S3 7RH

UK

email: v.bavula@sheffield.ac.uk
}


\begin{thebibliography}{99}

\bibitem{Bav-intdifline} V. V. Bavula,  The algebra of integro-differential operators on an affine line and its modules, {\em J. Pure Appl. Algebra}, {\bf 217} (2013)  495-529. Arxiv:math.RA: 1011.2997.


  \bibitem{Bav-Crit-S-Simp-lQuot} V. V. Bavula,  New criteria for a ring to have a semisimple left quotient ring,  {\em Journal of Alg. and its Appl.}, {\bf 14} (2015) no. 6, 1550090, 28pp.
 %DOI: 10.1142/S0219498815500905.  %Arxiv:math.RA:1303.0859.

\bibitem{larglquot} V. V. Bavula,  The largest left quotient ring of a ring,
 {\it Comm.  Algebra}, {\bf 44} (2016)
 no. 8, 3219-3261.  Arxiv:math.RA:1101.5107. 
 
 

\bibitem{Bav-LocArtRing} V. V. Bavula,   Left localizations of left Artinian rings,  {\em J. Algebra Appl.}
{\bf 15} (2016) no. 9, 1650165, 38 pages. %, DOI: 10.1142/S0219498816501656.
 Arxiv:math.RA:1405.0214.
 
\bibitem{Crit-lNoeth-lQuot} V. V. Bavula,    Criteria for a ring to have a left Noetherian largest  left quotient ring, {\it Algebras and Repr. Theory},  {\bf 21}  (2018),  no. 2, 359--373. % (published online on  05 July 2017).   
 
 
\bibitem{GWA-di-skew-2020} V. V. Bavula,   Generalized Weyl algebras and diskew polynomial rings, {\it J. Algebra Appl.}, {\bf 19} (2020) no. 10, 2050194, 43pp. (Arxiv:1612.08941).


 
 \bibitem{LocSets} V. V. Bavula, Localizable sets and the localization of a ring at a localizable set,   {\em J. Algebra}, {\bf 60} (2022), no.15, 38--75.  (ArXiv:2112.13447).
 
\bibitem{Loc-groupUnits} V. V. Bavula, Localizations of a ring at localizable sets, their groups of units and saturations, {\it Math. Comp. Sci.}, {\bf 16} (2022), no. 1, Paper No. 10, 15 pp. 

% \bibitem{Bav-FilDimSimCrit-1996} V. V. Bavula,     Filter dimension of algebras and modules, a simplicity
%criterion of   generalized Weyl algebras. {\it Comm. Algebra},  {\bf 24} (1996) no. 6, 1971--1992.


%\bibitem{Bav-GlGWA-1996} V. V. Bavula,  Global dimension of generalized Weyl algebras. {\it Representation   theory of algebras} (Cocoyoc, 1994), 81--107, CMS Conf. Proc., 18, Amer. Math. Soc.,   Providence,  RI, 1996.

%\bibitem{Bav-DimMultHolModCharp-2009} V. V. Bavula, Dimension, multiplicity, holonomic modules, and an analogue
%of the inequality of Bernstein for  rings of differential
%operators in prime characteristic, {\em Representation Theory},
%{\bf 13} (2009) 182-227. (Arxiv:math. RA/0605073).


\bibitem{Jategaonkar-LocNRings} A. V. Jategaonkar, Localization in Noetherian Rings, Londom Math. Soc. LMS 98, Cambridge Univ. Press, 1986.

\bibitem{Jacobson-StrRing} N. Jacobson, ``Structure of rings,''Am.
Math. Soc. Colloq., Vol. XXXVII (rev. ed.), Am. Math. Soc.,
Providence, 1968.

\bibitem{Lam-ModRingbook} T. Y. Lam, Lectures on modules and rings. Graduate Texts in Mathematics, {\bf 189}. Springer-Verlag, New York, 1999. xxiv+557 pp.

\bibitem{Lam-Exbook} T. Y. Lam, Exercises in Modules and Ring, Springer, 2007.

\bibitem{MR} J. C. McConnell and J. C. Robson,  Noncommutative Noetherian rings. With the cooperation of L. W. Small. Revised edition. Graduate Studies in Mathematics, 30. American Mathematical Society, Providence, RI, 2001. 636 pp.

\bibitem{Stenstrom-RingQuot} B. Stenstr\"{o}m, Rings of Quotients, Springer-Verlag, Berlin, Heidelberg, New York, 1975.


\end{thebibliography}
\end{document}